\documentclass[12pt,notitlepage]{amsart}
\usepackage{amsmath}     % [reqno] == equations numbered on right
\usepackage{amscd} 
\usepackage{tikz,tikz-cd}
\usepackage{import}
\usepackage{pdfpages}
\usepackage{xifthen}
\usepackage{hyperref}
\usepackage{transparent}
\usepackage{graphicx}
\usepackage{caption}
\usepackage{subcaption}

\usepackage[colorinlistoftodos]{todonotes}
\advance\voffset by -0.4 truein
\advance\hoffset by -0.4 truein
\newcommand{\F}{\mathcal F}

\newcommand{\Z}{\mathbb Z}

\newcommand{\p}{\mathbb P}

\newcommand{\g}{\mathfrak{V}}
\newcommand{\lp}{\mathbb{LP}}
\newcommand{\h}{\mathfrak{W}}

\newcommand{\vp}{\varphi}

\renewcommand{\deg}{\mathrm{deg}\,}

\theoremstyle{definition}
\newtheorem{theorem}[subsection]{Theorem}
\newtheorem{proposition}[subsection]{Proposition}
\newtheorem{Lem}[subsection]{Lemma}

\newtheorem{Def}[subsection]{Definition}

\newtheorem{Exa}[subsection]{Example}
 
\font\smallrm=cmr8
\font\smallsc=cmcsc10
\font\smallsl=cmsl10

\begin{document}
\author
[{\smallrm Eduardo Esteves, Renan Santos and Eduardo Vital}]
{Eduardo Esteves, Renan Santos and Eduardo Vital}
\title
[{\smallrm Quiver representations arising from degenerations of linear series, II}]
{Quiver representations arising from degenerations of linear series, II}

\date{\today}
\begin{abstract}  
We describe all the schematic limits of families of divisors associated to a given family of rank-$r$ linear series on a one-dimensional family of projective varieties degenerating to a connected reduced projective scheme $X$ defined over any field, under the assumption that the total space of the family is regular along $X$. More precisely, the degenerating family gives rise to a special quiver $Q$, called a \emph{$\mathbb{Z}^n$-quiver}, a special representation $\mathfrak L$ of $Q$ in the category of line bundles over $X$, called a \emph{maximal exact linked net}, and a special subrepresentation $\mathfrak V$ of the representation $H^0(X,\mathfrak L)$ induced from $\mathfrak L$ by taking global sections, called a \emph{pure exact finitely generated linked net} of dimension $r+1$. Given $\mathfrak g=(Q,\mathfrak L,\mathfrak V)$ satisfying these properties, we prove that the quiver Grassmanian $\mathbb{LP}(\mathfrak{V})$ of subrepresentations of $\mathfrak{V}$ of pure dimension 1, called a \emph{linked projective space}, is local complete intersection, reduced and of pure dimension $r$. Furthermore, we prove that there is a morphism $\mathbb{LP}(\mathfrak{V})\to\text{Hilb}_X$, and that its image parameterizes all the schematic limits of divisors along the degenerating family of linear series if $\mathfrak g$ arises from one.
\end{abstract}

\thanks{E.~Esteves was supported by CNPq, Proc.~304623/2015-6 and FAPERJ, Proc.~E-26/202.992/2017, and R.~Santos and E.~Vital were supported by CAPES doctor grants at IMPA}

\maketitle

\vspace{-.5cm}

\tableofcontents
% \emph{This is the first version of a work in progress where we study certain quiver representations arising from degenerations of linear series to singular curves.}

\section{Introduction}

This paper and its prequel \cite{Esteves_et_al_2021} aim to describe all the schematic limits of families of divisors associated to a given family of linear series on a family of projective varieties degenerating to a connected reduced projective scheme $X$ over any field $k$, under the assumption that the total space of the family is regular along $X$. We view these limits as points on the Hilbert scheme of $X$, and describe the subscheme containing them using the quiver Grassmannian of pure dimension~1 of a certain quiver representation. 

A \emph{linear series} is a vector space of global sections of a line bundle over a scheme defined over a field. Linear series are linearizations. For a smooth projective connected curve $C$, the Abel map $\text{Hilb}^d_C\to\text{Pic}^d_C$, associating a finite subscheme $D$ of $C$ of length $d$ to the corresponding line bundle $\mathcal O_C(D)$, is a fibration over an Abelian variety by projective spaces: the fiber over a line bundle $L$ is naturally isomorphic to $\mathbb P(V)$, where $V$ is the (complete) linear series of all global sections of $L$. This is Abel's theorem. Linear series can thus be thought of a certain collection of subschemes (effective divisors) of $C$ of the same length (degree). 

Given a family of linear series on a family of smooth curves degenerating to a singular curve $X$, the family of divisors associated to the linear series has as limit a collection of subschemes of $X$ of the same length. What is this collection? If $X$ is irreducible, it is a subscheme of $\text{Hilb}_X$ isomorphic to a projective space, as follows from work by Altman and Kleiman \cite{AK}. What if $X$ is reducible?

If $X$ is reducible, there are infinitely many linear series on $X$ that arise as limits along the family. These limits were studied by Eisenbud and Harris \cite{Eisenbud1986}, as well as later by Osserman \cite{Osserman2006} for when $X$ is a nodal curve of compact type, with two components in Osserman's case. Whereas Eisenbud and Harris proposed to consider a certain limit for each component of $X$, and called the collection of chosen limits a ``limit linear series," Osserman proposed to consider all limits whose associated line bundles have effective multidegrees, calling this collection a ``limit linear series" as well.

Even though certain notions of what a ``limit linear series" is have been proposed, notably by Eisenbud and Harris and by Osserman, using line bundles and sections, and they are different, there is certainly only one possible notion if one were to consider collections of subschemes. Curiously, we could not find in the literature any mention to the connection between a ``limit linear series" and schematic limits of divisors until work by the first author and Osserman \cite{EO}. There it turned out to be necessary to consider ``limit linear series" as defined by Osserman. However, the limits were considered only for nodal curves with two components and a single node, and only as cycles, in the symmetric product of $X$. 

It was only quite recently that Santana Rocha \cite{Rocha2019} was able to describe the limits in $\text{Hilb}^d_X$, though only for the simple curves considered in \cite{EO}. Remarkably, he made use of no new technique, but only of the linked Grassmannians that had already been introduced by Osserman in \cite{Osserman2006}. 

It must be said that the approach by Eisenbud and Harris, despite the lack of a fundamental connection between ``limit linear series" and schematic limits of divisors, 
yielded many important applications, a few of them listed 
in the introduction to \cite{Esteves_et_al_2021}. So many, in fact, and only using curves of compact type, that Eisenbud and Harris asked in \cite{EHBull}, p.~220, for a generalization of their theory to all nodal curves, writing that “...there is probably a small gold mine awaiting a general insight.” 

It is our goal to answer the question we posed above --- What if $X$ is reducible? --- in full generality, even for higher dimensional varieties, and by doing so, to show a path to answer the question in \cite{EHBull}. 

Let thus $X$ be a connected reduced projective scheme over a field $k$. As explained in \cite{Esteves_et_al_2021}, a rank-$r$ linear series on the general fiber of a \emph{regular smoothing} of $X$ gives rise to two quiver representations: a representation $\mathfrak L$ of a quiver $Q$ in the category of line bundles over $X$ and a subrepresentation $\mathfrak V$ of pure dimension $r+1$ of the induced representation $H^0(X,\mathfrak L)$ in the category of vector spaces over $k$ obtained from $\mathfrak L$ by taking global sections.

We have seen in loc.~cit.~that $Q$ is a special quiver, and $\mathfrak L$ and $\mathfrak V$ are 
special representations of $Q$, to be explained below:
\begin{enumerate}
    \item $Q$ is a \emph{$\mathbb Z^n$-quiver},
    \item $\mathfrak L$ is an \emph{exact} \emph{maximal} \emph{linked net}, 
    \item $\mathfrak V$ is a \emph{pure} \emph{exact} \emph{finitely generated} \emph{linked net}. 
\end{enumerate}

Conversely, let $\mathfrak g=(Q,\mathfrak L,\mathfrak V)$ be the data of a quiver $Q$, a representation $\mathfrak L$ of $Q$ in the category of line bundles over $X$ and a subrepresentation $\mathfrak V$ of a given pure dimension $r+1$ of $H^0(X,\mathfrak L)$ satisfying the special properties listed above. We prove in the current paper that there are a natural scheme structure for the quiver Grassmannian $\lp(\g)$ of subrepresentations of $\mathfrak V$ of pure dimension~$1$ for which $\lp(\g)$ is reduced and a  local complete intersection of pure dimension $r$ with rational irreducible components, and a natural morphism $\lp(\g)\to\text{Hilb}_X$ whose image is the collection of schematic limits of divisors associated to a degenerating family of linear series, if $\mathfrak g$ arises from one; see Theorems~\ref{Main_Theorem} and \ref{main_result}, Proposition~\ref{maplpg} and Theorem~\ref{final}.

We refrain from calling the above data $\mathfrak g$ a ``limit linear series," though we prove here it has every right to be called so!

We give more details now. First of all, a \emph{regular smoothing} of $X$ is the data of a flat projective map $\mathcal X\to B$, where $\mathcal X$ is regular and $B$ is the spectrum of a discrete valuation ring $R$ with residue field $k$, and an isomorphism of the special fiber of the map with $X$. 

Let $(L_\eta,V_\eta)$ be a linear series on the general fiber of a regular smoothing $\mathcal X\to B$ of $X$. Since $\mathcal X$ is regular, there is a line bundle extension $\mathcal L$ of $L_\eta$ to $\mathcal X$. Let $X_0,\dots,X_n$ be the irreducible components of $X$; they are Cartier divisors of $\mathcal X$. Every other line bundle extension of $L_\eta$ is of the form $\mathcal L_u:=\mathcal L(\sum \ell_i X_i)$ for a unique $(n+1)$-tuple $u=(\ell_0,\dots,\ell_n)\in\mathbb Z^{n+1}_{\geq 0}$ with $\min\{\ell_i\}=0$. The vertex set of $Q$ is precisely the set $Q_0$ of those $(n+1)$-tuples. As for the arrows, there is an arrow, and only one, connecting $u$ to $v$ if and only if $\mathcal L_u(X_i)\cong\mathcal L_v$ for some $i$. The representation $\mathfrak L$ associates to the vertices $u\in Q_0$ the line bundles $L_u:=\mathcal L_u|_X$ and to the arrows the restrictions to $X$ of the natural maps $\mathcal L_u\to\mathcal L_u(X_i)$. Finally, the vector space associated to $u\in Q_0$ by 
$\mathfrak V$ is the image in $H^0(X,L_u)$ of $V_\eta\cap H^0(\mathcal X,\mathcal L_u)$.

In \cite{Esteves_et_al_2021}, \S 3 and \S 4, we explained the special properties $\mathfrak g:=(Q,\mathfrak L,\g)$ satisfies, which we summarize here. First, a quiver is a \emph{$\mathbb Z^{n}$-quiver} if it is endowed with a partition of its arrow set in $n+1$ parts, called \emph{arrow types}, such that for each vertex there is a unique arrow of each type leaving it; paths containing the same number of arrows of each type are the circuits; each vertex is connected to each other by a path that does not contain arrows of all types, called an \emph{admissible path}, and two such paths contain the same number of arrows of each type. The quiver $Q$ arising from a degeneration of linear series is a $\mathbb Z^n$-quiver, the arrows partitioned by their association to the components $X_i$.

Second, a representation of a $\mathbb Z^n$-quiver $Q$ in a $k$-linear Abelian category is a \emph{linked net over $Q$} if the compositions of maps along nontrivial circuits are zero, along two admissible paths connecting the same two vertices are the same up to homothety, and along two admissible paths with no arrow type in common have trivially intersecting kernels. The representation $\mathfrak L$ arising from a degeneration of linear series is a linked net, and thus so is the representation $\mathfrak V$.

Third, a representation of a quiver in an Abelian category is \emph{finitely generated} if it is \emph{generated} by a finite set of vertices $H$, that is, if for each vertex $v$ there are paths $\gamma_1,\dots,\gamma_m$ leaving vertices of $H$ and arriving at $v$ such that the associated maps sum to an epimorphism to the object corresponding to $v$. Fourth, it is \emph{pure} if every epimorphism between the objects associated to each two vertices is an isomorphism. The linked net $\mathfrak V$ arising from a degeneration of linear series is clearly pure and is generated by the finite set of vertices corresponding to spaces with at least one section with finite vanishing.

Fifth, a representation of a $\mathbb Z^n$-quiver in an Abelian category is \emph{exact} if the kernel of the map associated to each nontrivial path $\gamma$ containing at most one arrow of each type, called \emph{simple}, is equal to the image of the map associated to a reverse path, a simple path taking the final point of $\gamma$ to its initial point. It is not completely straightforward, but the linked nets $\mathfrak L$ and $\mathfrak V$ arising from a degeneration of linear series are exact.

Finally, a representation of a quiver in the category of line bundles over $X$ is \emph{maximal}, if the map associated to each arrow is generically zero on one and only one irreducible component of $X$. The linked net $\mathfrak L$ arising from a degeneration of linear series is clearly maximal.

Conversely, let $\mathfrak g=(Q,\mathfrak L,\mathfrak V)$ be the data of a quiver $Q$, a representation $\mathfrak L$ of $Q$ in the category of line bundles over $X$ and a subrepresentation $\mathfrak V$ of dimension $r+1$ of $H^0(X,\mathfrak L)$ satisfying Properties (1)-(3) listed above. Here is a brief description of how we obtain that $\lp(\g)$ is a local complete intersection. We  jhjhstart with the important Theorem~\ref{12n}, which says that a point on $\lp(\g)$ is generated by a \emph{polygon}. A \emph{polygon} is a collection of vertices on a nontrivial minimal circuit of the quiver. Polygons appeared in \cite{Esteves_et_al_2021}: Its Prop.~10.1 claims that exact pure linked nets of vector spaces generated by polygons decompose as direct sums of exact pure linked nets of dimension 1, which are generated by vertices by \cite{Esteves_et_al_2021}, Thm.~7.8. Thus, if $\g$ is generated by a polygon, we can simultaneously diagonalize all maps of $\g$. We use this simplification to show that $\lp(\g)$ is a local complete intersection in this case; see Proposition~\ref{prop_triagle}. Finally we argue in the proof of Theorem~\ref{main_result} that for exact pure finitely generated linked nets $\g$ the scheme $\lp(\g)$ is isomorphic in a neighborhood of a point generated by a polygon $H$ to an open subscheme of $\lp(\g_H)$, where $\g_H$ is a certain pure exact linked net generated by $H$, which we define in Section~\ref{shadow}.

Rather than only describing the points on $\lp(\g)$ we describe in Section~\ref{lpg} the reduced subscheme $\lp(\g)^*_v\subseteq\lp(\g)$, parameterizing subrepresentations of $\g$ generated by the vertex $v$ for each $v\in Q_0$, and argue in Section~\ref{lpg_z2} that their closures $\lp(\g)_v$ are the irreducible components of $\lp(\g)$, our Theorem~\ref{Main_Theorem}, concluding that $\lp(\g)$ is generically reduced and of pure dimension $r$. As it is a local complete intersection, thus Cohen--Macaulay, it is reduced. We go beyond this to describe the stratification of $\lp(\g)$ induced by the $\lp(\g)_v$ in terms of minimal generation of subrepresentations; see Proposition~\ref{lemma: non_exact_point_is_in_frontier}.

Finally, in Section~\ref{divisors} we associate to each $\mathfrak W\in\lp(\g)$ the subscheme $Z(\mathfrak W)$ of $X$, the intersection of the zero schemes of all the sections given by $\mathfrak W$. We prove in Proposition~\ref{Zw} that the $Z(\mathfrak W)$ are numerically equivalent, which is enough, since $\lp(\g)$ is reduced, to show that the induced map $\lp(\g)\to\text{Hilb}_X$ to the Hilbert scheme of $X$ is a morphism; see Proposition~\ref{maplpg}. And we use that $\lp(\g)$ is a degeneration of $\mathbb P(V_\eta)$ if $\g$ arises from a degenerating linear series $(L_\eta,V_\eta)$ to show that the image of the morphism is the collection of schematic limits of the divisors 
associated to $(L_\eta,V_\eta)$; see Proposition~\ref{smoothg} and Theorem~\ref{final}.

Once one has explicit data, the objects we study here can be thoroughly described. In Section~\ref{secexample} we study the explicit example of the degeneration of the pencil of lines on a general pencil of curves degenerating to a union of lines in the plane, describing completely 
$\lp(\g)$. We describe as well $\lp(\g^*)$ and the map $X\to\lp(\g^*)$; see below.

There is plenty that we do not do here! As emphasized above, we consider only degenerations to $X$ of linear series along families whose total space is regular, what may not be the case even if $X$ is a nodal curve. In this special case though, one could argue that we could replace $X$ by a semistable model. This is not satisfactory however as we would obtain maps to different Hilbert schemes associated to different degenerations. The theory developed by Amini and the first author in \cite{AE1}, \cite{AE2} and \cite{AE3} might point out to a solution to this problem. 

Second, we consider only quiver Grassmanians of subrepresentations of pure dimension $1$. What about higher dimensions? It is proved by Osserman in \cite{OssermanFlatness}, Thm.~4.2, p.~3387, that quiver Grassmannians of pure subrepresentations of any dimension of pure exact finitely generated linked nets of vector spaces over $\mathbb Z^n$-quivers are Cohen--Macaulay if $n=1$. What about higher $n$?

Third, even if $\mathfrak g=(Q,\mathfrak L,\mathfrak V)$ arises from a degeneration of linear series to a nodal curve $X$, the morphism $\lp(\g)\to\text{Hilb}^d_X$ need not be an embedding! What is its image? Is there a natural resolution of $\text{Hilb}^d_X$ to which the morphism factors as an embedding? Can we obtain an Abel map for this resolution?

Fourth, linear series are useful to describe morphisms to projective spaces. If $\mathfrak g=(Q,\mathfrak L,\mathfrak V)$ arises from a degeneration of linear series to $X$, is there a natural morphism from $X$ to a natural degeneration of projective spaces? Could $\lp(\g^*)$ be the target of this morphism, where $\g^*$ is the dual representation of $\g$? We argue in Proposition~\ref{dual} that there is a natural rational map $X\dashrightarrow\lp(\g^*)$, but we do not argue that $\lp(\g^*)$ is a degeneration of projective spaces. Our theory does not apply to $\g^*$ as $\g^*$ may fail to be a linked net even when $n=1$!

This paper is organised as follows. In Section~\ref{linek_nets} we recall in more detail what $\mathbb Z^n$-quivers and linked nets are and introduce necessary notation. In Section~\ref{ln1} we prove that a finitely generated linked net over a $\mathbb Z^n$-quiver of simple objects in a $k$-linear Abelian category is generated by a polygon; see Theorem~\ref{12n}. In Section~\ref{ln1n2} we illustrate our proof with the classification of these linked nets for $n=2$. 

In Section~\ref{lpg} we define the linked projective space $\lp(\g)$ associated to a linked net $\g$ of vector spaces over a $\mathbb Z^n$-quiver $Q$, define its scheme structure if $\mathfrak V$ is pure and finitely generated, and describe the reduced subschemes $\lp(\g)_v$ which we prove in Section~\ref{lpg_z2} to cover $\lp(\g)$; see Theorem~\ref{Main_Theorem}. In fact, we do more: we describe each stratum in the natural stratification associated to the $\lp(\g)_v$, in particular describing when each selection of $\lp(\g)_v$ intersect nontrivially; see Propositions~\ref{lemma: non_exact_point_is_in_frontier}~and~\ref{cap}.

%In Section~\ref{lpg_z2} we restrict our attention to linked projective spaces of finitely generated linked nets over $\Z^2$-quivers. We prove Proposition~\ref{lemma: non_exact_point_is_in_frontier} that treats the points minimally generated from two and three vertices. We also prove Theorem~\ref{Main_Theorem}, which states that the irreducible components are the non-empty $\lp(\g)_{v}$ and the set of exact points is its nonsigular locus. 

In Section~\ref{shadow} we introduce the shadow partition of a $\mathbb Z^n$-quiver $Q$ associated to a finite set of vertices $H$ of $Q$ which is equal to its hull. Given a pure linked net $\g$ of vector spaces over $Q$ we define a new representation $\g_H$ of $Q$ generated by $H$, and prove that under certain conditions, for instance when $H$ is a polygon, $\g_H$ is a linked net, which is exact if so is $\g$; see Proposition~\ref{gH}.

The results in Section~\ref{shadow} will be crucial in Section~\ref{CH_LPg}, as we have already explained, to show that a pure exact finitely generated linked net $\g$ of vector spaces over a $\mathbb Z^n$-quiver is local complete intersection and reduced; see Theorem~\ref{main_result}.

In Section~\ref{smoothings} we use that $\lp(\g)$ is reduced to show that if $\g$ is smoothable, for instance, if $\g$ arises from a degeneration of linear series, then $\lp(\g)$ is a degeneration of projective spaces. 

We prove in Section~\ref{divisors} that a pure exact finitely generated subrepresentation $\g$ of 
$H^0(X,\mathfrak L)$ for a exact maximal linked net $\mathfrak L$ over a $\mathbb Z^n$-quiver $Q$ of line bundles over $X$ gives rise to a morphism $\lp(\g)\to\text{Hilb}_X$, whose image is the collection of schematic limits of divisors in a degenerating family of linear series, if $\mathfrak g:=(Q,\mathfrak L,\mathfrak V)$ arises from one; see Theorem~\ref{final}. Finally, in Section~\ref{secexample} we give an example.

We thank Omid Amini, Marcos Jardim and Oliver Lorscheid for many discussions on the subject. 

\section{\texorpdfstring{${\mathbb{Z}}^n$}{{\mathbb{Z}}n}-quivers and linked nets}\label{linek_nets}

Throughout the paper, $Q$ will be a fixed quiver. We fix a nontrivial partition of the arrow set of $Q$. Each part $\mathfrak a$ is called an \emph{arrow type}, and we say $a\in\mathfrak a$ has type $\mathfrak a$. The number of parts is $n+1$ with $n\in\mathbb N$. 

Given a path $\gamma$ in $Q$ and an arrow type $\mathfrak a$, we denote by $\mathfrak t_{\gamma}(\mathfrak a)$ the number of arrows of that type the path contains. We call $\mathfrak t_\gamma$ the \emph{type} of $\gamma$ and the collection of arrow types $\{\mathfrak a\,|\,\mathfrak t_{\gamma}(\mathfrak a)>0\}$ its \emph{essential type}. The path $\gamma$ is called \emph{admissible} if $\mathfrak t_{\gamma}(\mathfrak a)=0$ for some $\mathfrak a$ and \emph{simple} if $\mathfrak t_{\gamma}(\mathfrak a)\leq 1$ for every $\mathfrak a$. A simple non-admissible path is called a \emph{minimal circuit}.

We will assume the partition of the arrow set makes $Q$ into a $\mathbb Z^n$-quiver, that is, the following three conditions are satisfied:
\begin{enumerate}
    \item There is exactly one arrow of each type leaving each vertex.
    \item Each vertex is connected to each other by an admissible path.
    \item Two paths $\gamma_1$ and $\gamma_2$ leaving the same vertex arrive at the same vertex if and only if 
    $\mathfrak t_{\gamma_1}-\mathfrak t_{\gamma_2}$ is the constant function.    
    \end{enumerate}
 
 Two distinct vertices connected by a simple path are called \emph{neighbors}. If $v_1$ and $v_2$ are neighbors and $I$ is the essential type of a simple path connecting $v_1$ to $v_2$ we write $v_2=I\cdot v_1$. 
 
 A \emph{polygon} is a nonempty collection $\Delta$ of vertices of $Q$ which are pairwise neighbors. It is finite with at most $n+1$ vertices by \cite{Esteves_et_al_2021}, Prop.~5.9. Letting $m:=\#\Delta$, we call $\Delta$ a \emph{$m$-gon}. (A $2$-gon is a \emph{segment}, a $3$-gon is a \emph{triangle}.) Given a vertex $v\in\Delta$, there is a unique ordering $v_1,\dots,v_m$ of the vertices of $\Delta$ with $v_1=v$ and $v_{i+1}:=I_i\cdot v_i$ for $i=1,\dots,m-1$, where $I_1,\dots,I_{m-1}$ is a sequence of pairwise disjoint collections of arrow types; see \cite{Esteves_et_al_2021}, Prop.~5.9. In this case, we say $v_1,\dots,v_m$ \emph{form an oriented polygon}.
 
 See \cite{Esteves_et_al_2021} for basic properties of $\mathbb Z^n$-quivers.

We fix a field $k$ and call its elements \emph{scalars}.
We will consider representations $\g$ of $Q$ in $k$-linear Abelian categories, for instance, the category of vector spaces (over $k$) or the category of coherent sheaves on an $k$-scheme of finite type. For each vertex $v$ of $Q$, we denote by $V^\g_v$ the associated object, and for each path $\gamma$ in $Q$, we denote by $\varphi^\g_{\gamma}$ the corresponding composition of morphisms of $\g$. If $\g$ is clear from the context, we omit the superscript.

 Given a representation $\g$ of $Q$ in a $k$-linear Abelian category, $\g$ is \emph{pure} if each epimorphism between associated objects is an isomorphism. It is called 
 %\emph{binary} if each associated map is zero or a monomorphism, and 
 \emph{simple} if all associated objects are simple. We say $\g$ is $1$-\emph{generated} by a collection $H$ of vertices if for each vertex $v$ of $Q$ there is $u\in H$ and a path $\gamma$ connecting $u$ to $v$ such that $\varphi_{\gamma}$ is an epimorphism. If $\g$ is not $1$-generated by a smaller collection, we say $H$ is \emph{minimal}. 
 We say $\g$ is \emph{generated} by a collection $H$ of vertices if for each vertex $v$ of $Q$ there are paths $\gamma_1,\dots,\gamma_m$ connecting vertices of $H$ to $v$ such that 
 $V_v=\sum\text{Im}(\varphi_{\gamma_i})$. 

%A collection of four distinct vertices $\{v_0,v_1,v_2,v_3\}$ is called a \emph{paralelogram} if they are pairwise neighbors for which, up to re-ordering, there are non-empty subsets $I,J\subseteq\{0,\dots,n\}$ with $I\cap J=\emptyset$ and $I\cup J\neq\{0,\dots,n\}$ such that
%$$
%v_1=I\cdot v_0,\quad v_2=J\cdot v_0,\quad v_3=J\cdot v_1,\quad v_3=I\cdot v_2.
%$$

We say $\g$ is a \emph{weakly linked net} over $Q$ if $\g$ satisfies the following two conditions, 
 \begin{enumerate}
 \item if $\gamma_1$ and $\gamma_2$ are two paths connecting the same two vertices and $\gamma_2$ is admissible then $\varphi_{\gamma_1}$ is a scalar multiple of $\vp_{\gamma_2}$; 
 \item $\vp_\gamma=0$ for each minimal circuit $\gamma$;
 \end{enumerate}
 and we say it is a \emph{linked net} if in addition a third condition is verified:
 \begin{enumerate}
 \item[(3)] if $\gamma_1$ and $\gamma_2$ are two admissible paths leaving the same vertex with no arrow type in common then $\text{Ker}(\varphi_{\gamma_1})\cap\text{Ker}(\varphi_{\gamma_2})=0$.
 \end{enumerate}
 %We call $\g$ \emph{special} (resp.~\emph{general}) if $\varphi_\gamma$ is zero (resp.~nonzero) for every minimal circuit $\gamma$.

Clearly, if $\g$ is a weakly linked net that is $1$-generated by a finite set, then $\g$ is generated by a finite set. The converse holds as well, by \cite{Esteves_et_al_2021}, Prop.~6.4. In this case, we say $\g$ is \emph{finitely generated}. If $\g$ is finitely generated then $\g$ is \emph{locally finite} by \cite{Esteves_et_al_2021}, Prop.~6.5, that is, for each vertex $v$ of $Q$ there is an integer $\ell$ such that $\vp_\mu=0$ for each path $\mu$ arriving at $v$ with length greater than $\ell$.
 
 For each vector $v$ in a vector space $V$, we will denote by $[v]$ the set of its nonzero scalar multiples. In a $k$-linear category, the set of morphisms between any two objects is a vector space, so given a morphism $\vp$, we may consider the set $[\vp]$. We let $\text{Ker}[\vp]:=\text{Ker}(\vp)$ and $\text{Im}[\vp]:=\text{Im}(\vp)$. Also, since composition is $k$-bilinear, $[\psi][\vp]:=[\psi\vp]$ when $\psi\vp$ is defined. We write $[\vp]=0$ if $\vp=0$ and 
say $[\vp]$ is an isomorphism (resp.~monomorphism, resp.~epimorphism) if so is $\vp$.
 
 Given two vertices $v_1$ and $v_2$ of $Q$, let $\varphi^{v_1}_{v_2}:=[\varphi_{\gamma}]$ for any admissible path $\gamma$ connecting $v_1$ to $v_2$. If $\g$ is a weakly linked net, $\vp^{v_1}_{v_2}$ is well defined. In addition, if $v_1=v_2$, the class $\varphi^{v_1}_{v_2}$ is an isomorphism; otherwise $\vp^{v_2}_{v_1}\varphi^{v_1}_{v_2}=0$, or equivalently, $\text{Im}(\varphi^{v_1}_{v_2})\subseteq\text{Ker}(\varphi^{v_2}_{v_1})$. 
 We say $\g$ is \emph{exact} if  $\text{Im}(\varphi^{v_1}_{v_2})=\text{Ker}(\varphi^{v_2}_{v_1})$ for each two neighbors $v_1$ and $v_2$. 
 %Notice that, since $\g$ has pure dimension, 
 %$$
%\dim\text{Ker}(\varphi^{v_2}_{v_1})-\dim\text{Im}(\varphi^{v_1}_{v_2})=\dim\text{Ker}(\varphi^{v_1}_{v_2})-\dim\text{Im}(\varphi^{v_2}_{v_1});
%$$
%thus it is enough to check only one equality for each pair of neighbors.

If $\g$ is a weakly linked net of vector spaces, then $\g$ is pure if and only if the associated spaces have the same finite dimension, which we call the dimension of $\g$ and denote $\dim\g$. It is simple if in addition $\dim\g=1$.
 
\section{Simple linked nets}\label{ln1}

\begin{Lem}\label{seq} Let $I$ be a nonempty proper collection of arrow types of a $\mathbb Z^n$-quiver $Q$ and $v_1,v_2,v_3$ vertices of $Q$ such that $v_2=I\cdot v_1$ and $v_3=I\cdot v_2$. Let $\g$ be a weakly linked net over $Q$. Then the following statements hold:
\begin{enumerate}
    \item If $\vp^{v_1}_{v_2}$ is an epimorphism then $\vp^{v_2}_{v_1}$ is zero.
    \item If $\vp^{v_2}_{v_1}$ is zero and $\g$ is a linked net then $\vp^{v_2}_{v_3}$ is an monomorphism.
\end{enumerate}
\end{Lem}

\begin{proof} If $\vp^{v_1}_{v_2}$ is an epimorphism, since $\vp^{v_2}_{v_1}\vp^{v_1}_{v_2}=0$, we have that $\vp^{v_2}_{v_1}=0$, proving the first statement. 

Assume now that $\vp^{v_2}_{v_1}=0$. Since $v_2=I\cdot v_1$, there is a simple admissible path $\gamma_1$ connecting $v_2$ to $v_1$ with essential type $T-I$, where $T$ is the set of arrow types of $Q$. Since $\vp^{v_2}_{v_1}=[\vp_{\gamma_1}]$, we have that $\vp_{\gamma_1}=0$. On the other hand, there is a simple admissible path $\gamma_2$ connecting $v_2$ to $v_3$ with essential type $I$. Since $\g$ is a linked net, $\text{Ker}(\vp_{\gamma_1})\cap\text{Ker}(\vp_{\gamma_2})=0$. Since $\vp_{\gamma_1}=0$, it follows that $\vp_{\gamma_2}$ is a monomorphism. Of course 
$\vp^{v_2}_{v_3}=[\vp_{\gamma_2}]$, thus $\vp^{v_2}_{v_3}$ is a monomorphism.
\end{proof}

\begin{Def} Let $\g$ be a weakly linked net over a $\mathbb Z^n$-quiver $Q$. 
Let $v_1$ and $v_2$ be neighboring vertices of $Q$. We call $v_1$ and $v_2$ 
\emph{unrelated neighbors} for $\g$ if $\varphi^{v_1}_{v_2}=0$ and  $\varphi^{v_2}_{v_1}=0$, and call them \emph{related} otherwise. 
\end{Def}

If $\g$ is simple, then $v_1$ and $v_2$ are related if and only if  $\text{Im}(\varphi^{v_1}_{v_2})=\text{Ker}(\varphi^{v_2}_{v_1})$;
%, and thus if and only if $\text{Im}(\varphi^{v_2}_{v_1})=\text{Ker}(\varphi^{v_1}_{v_2})$. 
thus, $\g$ has only related neighbors if and only if $\g$ is exact. 

%Also, notice that, if $v_2=I\cdot v_1$ and $\varphi^{v_2}_{v_1}=0$, then $\varphi^{v_2}_w$ is injective, whence an isomorphism, for each $w\in C^I(v_2)$. 

%Given a vertex $v$ of $Q$ and a subset $I\subseteq\{0,\dots,n\}$, the \emph{cone} $C^I(v)$ is the collection of vertices connected to $v$ by admissible paths $\gamma$ with $\gamma(i)>0$ only if $i\in I$.  

% We say $(V_{v_1}, V_{v_2})$ is \emph{perfect} or a \textit{generating pair} if the maps $\vp^{v_2}_{R+T}$ and
% $\vp^{v_2}_{S+T}$ are both isomorphisms. This is equivalent to say that the restriction of $\g$ to $C^{R+T}(v_2) \cup C^{S+T}(v_2)$ is supported
% on $\{v_2\}$. As the restriction of $\g$ to $C^{R+S}(v_1)$ is supported on $\{v_1\}$, if $(V_{v_1}, V_{v_2})$ is perfect, then
% $\g$ has minimal support on the independent segment $\{v_1,v_2\}$. 

Under certain conditions, as we will see below, unrelated neighbors give rise to more unrelated neighbors.

\begin{Lem}\label{triangle} Let $\g$ be a weakly linked net over a $\mathbb Z^n$-quiver $Q$. Let $v_1,v_2,v_3$ be vertices of $Q$ forming an oriented triangle and $\g$ be a linked net over $Q$. Then:
\begin{enumerate}
    \item $\vp^{v_2}_{v_3}\vp^{v_1}_{v_2}=\vp^{v_1}_{v_3}$.
    \item If $\vp^{v_1}_{v_2}$ is an isomorphism then $v_1$ and $v_3$ are unrelated if and only if $v_2$ and $v_3$ are unrelated.
\end{enumerate}
\end{Lem}

\begin{proof} The first statement follows from the fact that there is an admissible path connecting $v_1$ to $v_3$ through $v_2$. Assume $\vp^{v_1}_{v_2}$ is an isomorphism. Since $\vp^{v_2}_{v_3}\vp^{v_1}_{v_2}=\vp^{v_1}_{v_3}$ and $\vp^{v_1}_{v_2}$ is an epimorphism, we have that  
$\vp^{v_1}_{v_3}=0$ if and only if $\vp^{v_2}_{v_3}=0$. And since $\vp^{v_1}_{v_2}\vp^{v_3}_{v_1}=\vp^{v_3}_{v_2}$ and $\vp^{v_1}_{v_2}$ is a monomorphism, we have that $\vp^{v_3}_{v_1}=0$ if and only if $\vp^{v_3}_{v_2}=0$. Thus $v_1$ and $v_3$ are unrelated if and only if $v_2$ and $v_3$ are unrelated. 
\end{proof}

\begin{Lem}\label{minsup1}
Let $\g$ be a finitely generated weakly linked net over a $\mathbb{Z}^n$-quiver $Q$. Then there is a unique collection $H$ of vertices $1$-generating $\g$ contained in every such collection. Furthermore, $H$ is finite and if $\g$ is simple then $\vp^u_v=0$ 
for all distinct $u,v\in H$. 
\end{Lem}

\begin{proof} By \cite{Esteves_et_al_2021}, Prop.~6.4, there is a finite set of vertices $H'$ that $1$-generates $\g$. It contains a minimal such collection $H$. By loc.~cit., Prop.~6.3, we have that $H$ is contained in every collection of vertices $1$-generating $\g$. The uniqueness of such a $H$ is clear. Furthermore, it follows from loc.~cit., Prop.~6.3 that 
$\vp^v_w$ is not an epimorphism for distinct $v,w\in H$. Thus, if $\g$ is simple, then $\vp^v_w=0$ for distinct $v,w\in H$.
\end{proof}

\begin{Def} Let $\g$ be a simple linked net over a $\mathbb Z^n$-quiver $Q$. A polygon of $Q$ is said to be \emph{unrelated} for $\g$ if each two vertices of it are unrelated for $\g$.
\end{Def}

\begin{theorem}\label{12n} Let $\g$ be a locally finite simple linked net over a $\mathbb Z^n$-quiver. Then $\g$ is (minimally) generated by a polygon. Furthermore, the size of the polygon minimally generating $\g$ is the maximum size of the unrelated polygons for $\g$. 
\end{theorem}

\begin{proof} If $\g$ is exact then $\g$ generated by a vertex (a $1$-gon) by \cite{Esteves_et_al_2021}, Thm.~7.8. Also, there are no unrelated vertices for $\g$. Suppose now $\g$ is not exact. Then there are unrelated neighbors for $\g$. Let $m$ be the maximum number for which there is an unrelated $(m+1)$-gon for $\g$. Then $m\geq 1$. It will be enough to show that there is an unrelated $(m+1)$-gon for $\g$ generating $\g$. 

By \cite{Esteves_et_al_2021}, Prop.~5.9, there is a minimal circuit $a_n\cdots a_0$ such that, denoting by $w_i$ the initial vertex of $a_i$ for each $i$, there are $m+1$ vertices $v_0,\dots,v_m$ among $w_0,\dots,w_n$ which are unrelated for $\g$. Order the $v_j$ such that $v_j=w_{r_j}$ for an increasing sequence of integers $r_0,\dots,r_m$ with $r_0=0$ and $r_m\leq n$. Let $\mathfrak a_i$ be the type of $a_i$ for each $i$. For convenience, put $r_{m+1}:=n+1$ and $a_{n+1}:=a_0$.

The proof consists of a procedure for changing the minimal circuit and the $v_i$ in such a way that at the end $\{v_0,\dots,v_m\}$ generates $\g$. We describe it in steps below.

\vskip0.2cm

{\bf Step 1.} The minimal circuit $a_n\cdots a_0$ and the $v_j$ can be chosen such that $\vp_{a_i}=0$ if and only if $i=r_j-1$ for 
some $j$. 

\vskip0.1cm 

It is enough to prove that $\vp_{a_i}=0$ for exactly $m+1$ values of $i$, as the final vertices of these $a_i$ form a set of unrelated vertices for $\g$, due to $m\geq 1$. For each $j=0,\dots,m$, there is a unique arrow $a$ among $a_{r_j},\dots,a_{r_{j+1}-1}$ such that $\vp_a$ is zero. Indeed, since $v_j$ and $v_{j+1}$ are unrelated, $\vp^{v_j}_{v_{j+1}}=0$ and thus $a$ exists. But if there were $a_i$ and $a_\ell$ with $r_j\leq i<\ell<r_{j+1}$ such that 
$\vp_{a_i}$ and $\vp_{a_\ell}$ are zero, then $v_0,\dots,v_j,w_{i+1},v_{j+1},\dots,v_m$ would form an oriented $(m+2)$-gon of unrelated vertices, contradicting the maximality of $m$. 
\vskip0.2cm

{\bf Step 2.} In addition, for each arrow type $\mathfrak b$, the minimal circuit $a_n\cdots a_0$ and the $v_j$ can be chosen such that $\vp_b=0$ for the arrow $b$ of type $\mathfrak b$ arriving at $v_0$.

\vskip0.1cm 

Let $b$ be an arrow arriving at $v_0$ of type $\mathfrak b$. If $\mathfrak b=\mathfrak a_n$ then $b=a_n$, whence $\vp_b=0$. Assume $\mathfrak b\neq\mathfrak a_n$. Then $\mathfrak b=\mathfrak a_j$ for (a unique) $j<n$. Consider the following minimal circuit:
$$
e_n\cdots e_ja_{j-1}\cdots a_0;
$$
see Figure~\ref{fig:12n}.
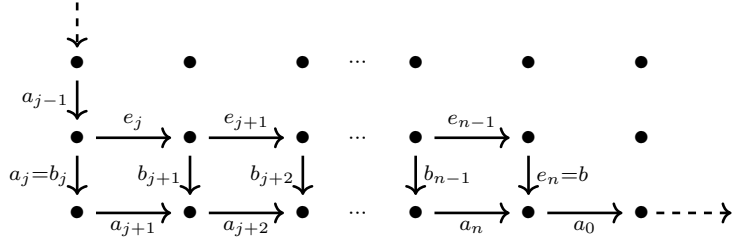
\begin{figure}[ht]
    \centering
    \begin{tikzpicture}
     %the vertices
    \foreach \j in {-1,0,1}
    \foreach \i in {-3,-2,-1,0,1,2}
    \node (\i\j) at (1.5*\i,\j) {$\bullet$};
    % \node (-21) at (-2,1) {$\bullet$};
    % \node (2-1) at (2,-1) {$\bullet$};
      % the extended arrows
    \draw[->, dashed, line width = 1] (-4.5,1.8) -- (-4.5,1.2);
    \draw[->, dashed, line width = 1] (3.2,-1) -- (4.2,-1);
    % the arrows
    \path[commutative diagrams/.cd, every arrow, every label] 
       % arrows horizontal 
    (00)   edge[line width = 1, ->]  node[] {$e_{n-1}$}  (10)
    (0-1)  edge[line width = 1, ->]  node[swap] {$a_{n}$}    (1-1)
    (1-1)  edge[line width = 1, ->]  node[swap] {$a_{0}$}    (2-1)
    (00)   edge[draw = none]         node[sloped,auto=false] {$\cdots$}       (-10)
    (01)   edge[draw = none]         node[sloped,auto=false] {$\cdots$}       (-11)
    (-10)  edge[line width = 1, <-]  node[swap] {$e_{j+1}$} (-20)
    (-1-1) edge[line width = 1, <-]  node {$a_{j+2}$} (-2-1)
    (-1-1) edge[draw = none]         node[sloped,auto=false] {$\cdots$}  (0-1)
    (-30)  edge[line width = 1, ->] node[] {$e_j$}  (-20)
    (-3-1) edge[line width = 1, ->] node[swap] {$a_{j+1}$}  (-2-1)
    % arorws vertical
    (-30) edge[line width = 1, <-] node[] {$a_{j-1}$} (-31)
    (-3-1) edge[line width = 1, <-] node[] {$a_j=b_j$} (-30)
    (-2-1) edge[line width = 1, <-] node[] {$b_{j+1}$} (-20)
    (-1-1) edge[line width = 1, <-] node[] {$b_{j+2}$} (-10)
    (0-1) edge[line width = 1, <-]  node[swap] {$b_{n-1}$}(00)
    (1-1) edge[line width = 1, <-]  node[swap] {$e_{n}=b$} (10)
      ;
    \end{tikzpicture}
    \caption{Proof of Theorem~\ref{12n}}
    \label{fig:12n}
\end{figure}
Here $e_\ell$ is an arrow of type $\mathfrak a_{\ell+1}$ 
for $\ell=j,\dots,n-1$ and $e_n:=b$. The initial vertex of $e_\ell$ is connected to the initial vertex of $a_{\ell+1}$ by an arrow $b_\ell$ of type $\mathfrak b$ for $\ell=j,\dots,n$. Of course, $b_j=a_j$ and $b_n=b$. 

Suppose $\vp_b\neq 0$. We claim that $\vp_{b_\ell}\neq 0$ and that $\vp_{e_\ell}=0$ if and only if $\vp_{a_{\ell+1}}=0$ for each $\ell=j,\dots,n-1$. In particular, $\vp_{e_{n-1}}=0$ and $\vp_{b_j}\neq 0$. Indeed, $\vp_{b_n}\neq 0$ because $b_n=b$. Assume by descending induction on $\ell$ that we have proved that $\vp_{b_\ell}\neq 0$. 
%since $\vp_{a_n}=0$ and $[\vp_{a_n}][\vp_{b_{n-1}}]=[\vp_b][\vp_{e_{n-1}}]$, it follows that $\vp_{e_{n-1}}=0$, and hence $\vp_{b_{n-1}}\neq 0$ by the third property of a linked net, the case $\ell=n-1$ of our claim. Assume by descending induction on $\ell$ that we have proved that $\vp_{b_\ell}\neq 0$ and that $\vp_{e_\ell}= 0$ if and only if $\vp_{a_{\ell+1}}= 0$ for a certain $\ell>j$. 
If 
$\vp_{e_{\ell-1}}=0$ then $\vp_{b_{\ell-1}}\neq 0$ by the third property of a linked net and hence $\vp_{a_{\ell}}=0$ because
\begin{equation}\label{1234}
[\vp_{a_\ell}][\vp_{b_{\ell-1}}]=[\vp_{b_{\ell}}][\vp_{e_{\ell-1}}].
\end{equation} 
And if $\vp_{e_{\ell-1}}\neq 0$, since also $\vp_{b_\ell}\neq 0$, we have that 
$\vp_{a_\ell}\neq 0$ and $\vp_{b_{\ell-1}}\neq 0$ from\eqref{1234}.
The claim is proved.

It follows from the claim that $\vp_b=0$ if $\mathfrak b=\mathfrak a_{r_\ell-1}$ for some $\ell$, because then $b_j=a_{r_\ell-1}$ and thus $\vp_{b_j}=0$.

Now, if $\vp_b\neq 0$ we replace the minimal circuit $a_n\dots a_0$ by the also minimal circuit $e_{n-1}\cdots e_ja_{j-1}\cdots a_0b$. The latter starts at a different vertex, the initial vertex of $b$, but has a pattern similar to that of the former; in particular, the arrows corresponding to zero maps have the same types in both circuits. We will call the latter circuit the \emph{$\mathfrak b$-shift} of the former centered at $v_0$.

Since $\mathfrak b\neq\mathfrak a_n$, we can make a sequence of $\mathfrak b$-shifts centered at initial vertices, as long as the arrow $b$ of type $\mathfrak b$ arriving at the inital vertex of each $\mathfrak b$-shift in the sequence satisfies $\vp_b\neq 0$. The sequence is necessarily finite though, since $\g$ is locally finite. 
%The upshot is that, for a given arrow type $\mathfrak b$, we may assume that the arrow $b$ of that type arriving at $v_0$ satisfies $\vp_b=0$.

\vskip0.2cm

{\bf Step 3.} In addition, the minimal circuit $a_n\cdots a_0$ and the $v_j$ can be chosen such that $\vp_b=0$ for each arrow arriving at $v_0$.

\vskip0.1cm 

Let $\mathfrak b_1,\mathfrak b_2,\dots,\mathfrak b_p$ be pairwise distinct arrow types. 
%, each distinct from the type of $a_{r_i-1}$ for every $i=0,\dots,m$. 
such that the arrow $b_i$ of type $\mathfrak b_i$ arriving at $v_0$ satisfies $\vp_{b_i}=0$ for each $i=1,\dots,p-1$. Let $b_p$ be the arrow of type $\mathfrak b_p$ arriving at $v_0$. Suppose $\vp_{b_p}\neq 0$. Consider the $\mathfrak b_p$-shift of the minimal circuit centered at $v_0$. For each $i=1,\dots, p$, let $b'_i$ be the arrow of type 
$\mathfrak b_i$ arriving at the initial vertex of $b_p$, which is the initial vertex of the new minimal circuit. Since $\vp_{b_p}\vp_{b'_i}$ factors through $\vp_{b_i}$, we must have $\vp_{b'_i}=0$ for each $i=1,\dots,p-1$. If $\vp_{b'_p}\neq 0$, we may then consider the $\mathfrak b_p$-shift of the new minimal circuit centered at its initial vertex, and proceed as above. Again since $\g$ is locally finite, we must arrive at a minimal circuit such that $\vp_{b}=0$ for the arrow $b$ of type $\mathfrak {b_i}$ arriving at the initial vertex for each $i=1,\dots,p$. 

\vskip0.2cm

{\bf Step 4.} In addition, the minimal circuit $a_n\cdots a_0$ and the $v_j$ can be chosen such that $\vp_b=0$ for each arrow arriving at $v_j$ for each $j=0,\dots,m$.

\vskip0.1cm 

Assume that $\vp_b=0$ for each arrow  $b$ arriving at $v_\ell$ for $\ell=j,\dots,m$. Let $b$ be an arrow arriving at $v_0$ and $\mathfrak b$ its type. 
%If $\mathfrak b$ is the type of $a_{r_i-1}$ for some $i$ then $\vp_b=0$, as we have seen. 
Suppose $\vp_b\neq 0$, and consider the $\mathfrak b$-shift of the minimal circuit centered at $v_0$. Let $i$ be such that $\mathfrak a_i=\mathfrak b$. Then $i<n$. But also, $i\geq r_m$, because otherwise, as we have seen in Step 2, the arrow $b'$ of type $\mathfrak b$ arriving at $v_m$ would satisfy $\vp_{b'}\neq 0$. But then the $\mathfrak b$-shift of the minimal circuit does not change the vertices $v_1,\dots,v_m$, only $v_0$ gets replaced by the inital vertex of $b$. We may thus proceed as in Steps 2 and 3, to obtain a minimal circuit for which $\vp_b=0$ for each arrow $b$ arriving at $v_j,\dots,v_m$ and at $v_0$. Reordering the arrows of the minimal circuit, we may assume that $\vp_b=0$ for each arrow  $b$ arriving at $v_{j-1},\dots,v_m$ and repeat the substep.

\vskip0.2cm

{\bf Step 5.} With the minimal circuit $a_n\cdots a_0$ and the $v_j$ chosen such that $\vp_b=0$ for each arrow arriving at $v_j$ for each $j=0,\dots,m$, we have that $\g$ is generated by $\{v_0,\dots,v_m\}$. 

\vskip0.1cm 

Let $v$ be any vertex of $Q$. Choose $w$ among the $w_i$ such that the admissible paths from $w$ to $v$ have the smallest length. Consider such an admissible path $\gamma$. Let $j\in\{0,\dots,m\}$ such that $w=w_i$ for some $i$ satisfying $r_j\leq i<r_{j+1}$. We claim that $\vp^{v_j}_v$ is an isomorphism.

We may choose $\gamma$ such that that there are admissible paths 
$\gamma_1$, $\gamma_2$ and $\gamma_3$ satisfying that  $\gamma=\gamma_2\gamma_1$, and that $\gamma_3\gamma_1$ connects $w$ to $v_{j+1}$, and such that $\gamma_2$ and $\gamma_3$ have no arrow type in common. Then 
$$
\text{Ker}(\vp_{\gamma_2})\cap
\text{Ker}(\vp_{\gamma_3})=0.
$$
As $\gamma_3$ arrives at $v_{j+1}$, we must have $\vp_{\gamma_3}=0$, and hence $\vp_{\gamma_2}$ is an isomorphism. Also $\vp_{\gamma_1}$ is an isomorphism. Indeed, $\gamma_3$ is nontrivial by the choice of $w$. Were $\vp_{\gamma_1}=0$, the path $\gamma_1$ would contain an arrow $b$ with $\vp_b=0$. Letting $z$ be the final point of $b$, we would have 
that $z\neq v_{j+1}$ and that 
$\{v_0,\dots,v_j,z,v_{j+1},\dots,v_m\}$ would be a $(m+2)$-gon of unrelated vertices for $\g$, contradicting the maximality of $m$. Thus 
$\vp_\gamma$ is an isomorphism, and hence so is 
$\vp^{v_j}_v$ because $\gamma a_{i-1}\cdots a_{r_j}$ 
connects $v_j$ to $v$ and $\vp_{a_\ell}$ is an isomorphism for each $\ell=r_j,\dots,r_{j+1}-2$.
\end{proof}

Notice that, since every subset of a polygon is a polygon, it follows from Theorem~\ref{12n} that a simple locally finite linked net is minimally generated by a polygon.

\section{Simple linked nets over \texorpdfstring{$\mathbb Z^2$}{Z2}-quivers}\label{ln1n2}

%An exact pure finitely generated linked net of dimension 1 is generated by a vertex, by \cite{Esteves_et_al_2021}, Thm.~6.4. 
In this section we will illustrate the polygons minimally generating a non-exact simple locally finite linked net over a $\mathbb Z^2$-quiver. 
We will thus assume $n=2$. 

\begin{figure}[ht]
    \centering
    \begin{tikzpicture}
      \node (a) at (0,0) {$\bullet$};
      \node (b) at (0,2) {$\bullet$};
      \node (c) at ({sqrt(3)},1) {$\bullet$};
      \path[commutative diagrams/.cd, every arrow, every label]  
            (b) edge[->] node[fill=white,anchor=center, pos=0.5] {$\mathfrak a_2$} (a)
            (a) edge[->] node[fill=white,anchor=center, pos=0.5] {$\mathfrak a_0$} (c)
            (c) edge[->] node[fill=white,anchor=center, pos=0.5] {$\mathfrak a_1$} (b)
    ;
    \end{tikzpicture}
    \caption{Partition of the arrow set.}
    \label{fig:triangle}
\end{figure}
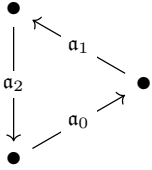

As each two $\mathbb Z^2$-quivers are equivalent by \cite{Esteves_et_al_2021}, Prop.~2.4, 
we may consider a particular representation of $Q$ as a planar quiver, namely: (see Figure~\ref{fig:triangle})
\begin{enumerate}
    \item[(0)] $\mathfrak a_0$ is the set of arrows from South-West to North-East.
    \item[(1)] $\mathfrak a_1$ is the set of arrows from South-East to North-West.
    \item[(2)] $\mathfrak a_2$ is the set of arrows from North to South.
\end{enumerate}
Given the ordering of the partition, we will find it more natural to say that an arrow is of type $\mathfrak a_i$ instead of type $A_i$. 
%So, now the total set of arrow types is $\{0,1,2\}$.

The planar representation will render our statements more descriptive. For starters, if $v_1,v_2,v_3$ form an oriented triangle then there are arrows connecting $v_1$ to $v_2$, $v_2$ to $v_3$ and $v_3$ to $v_1$, one of each type. Up to reordering, we may assume an arrow of type $0$ connects $v_1$ to $v_2$. If an arrow of type 2 connects $v_2$ to $v_3$, we say 
$v_1,v_2,v_3$ form a \emph{clockwise} triangle; otherwise we say they form a \emph{counterclockwise} triangle. The wording is natural, as can be seen in Figure~\ref{fig:triangle}.

\begin{figure}[ht]
\scalebox{0.9}{\begin{minipage}{\textwidth}
  \centering
\begin{subfigure}{.32\textwidth}
\centering
\begin{tikzpicture}
% the grid
\foreach \i in {-1,0,1}
\node (\i-3) at ({\i*sqrt(3)},{-1.5}) {$\bullet$};
\foreach \i in {-1,0}
\node (\i-2) at ({\i*sqrt(3)+.5*sqrt(3)},-1) {$\bullet$};
\foreach \i in {-1,0,1}
\node (\i-1) at ({\i*sqrt(3)},-.5) {$\bullet$};
\foreach \i in {-1,0}
\node (\i0) at ({\i*sqrt(3)+.5*sqrt(3)},0) {$\bullet$};
\foreach \i in {-1,0,1}
\node (\i1) at ({\i*sqrt(3)},.5) {$\bullet$};
\foreach \i in {-1,0}
\node (\i2) at ({\i*sqrt(3)+.5*sqrt(3)},1) {$\bullet$};
\foreach \i in {-1,0,1}
\node (\i3) at ({\i*sqrt(3)},1.5) {$\bullet$};
\foreach \i in {-1,0}
\node (\i4) at ({\i*sqrt(3)+.5*sqrt(3)},2) {$\bullet$};
\filldraw[orange] ({.5*sqrt(3)},1) circle (3pt);
\filldraw[orange] (0,.5) circle (3pt);
% the arrows
\path[commutative diagrams/.cd, every arrow, every label]  
%     % red arrows, the imperfect linkages 
%             % red vertical
            (01) edge[red, line width = 1, <-]  (03)
            (-12) edge[red, line width = 1, <-]  (-14)
%             % red right up
            (-13) edge[red, line width = 1, ->]  (-14)
            (-12) edge[red, line width = 1, ->]  (03)
            (01) edge[red, line width = 1, ->]  (02)
            (0-1) edge[red, line width = 1, ->]  (00)
            (0-3) edge[red, line width = 1, ->]  (0-2)
%             % red left up
            (00) edge[red, line width = 1, ->]  (01)
            (0-2) edge[red, line width = 1, ->]  (0-1)
%     % blue arrows, the perfect linkage
%             % blue vertical 
            (00) edge[blue, line width = 1, ->]  (0-2)
            (01) edge[blue, line width = 1, ->]  (0-1)
            (02) edge[blue, line width = 1, ->]  (00)
            (0-1) edge[blue, line width = 1, ->]  (0-3)
%             %blue left up
            (03) edge[blue, line width = 1, ->]  (-14)
            (02) edge[blue, line width = 1, ->]  (03)
            (01) edge[blue, line width = 1, ->]  (-12)
            (-12) edge[blue, line width = 1, ->]  (-13)
%             % blue right up 
            (02) edge[blue, line width = 1, ->]  (13)
            ;
\end{tikzpicture}
    \caption{Type I.}
    \label{fig: type_I}
\end{subfigure}
\hfill 
\begin{subfigure}{.32\textwidth}
\centering
\begin{tikzpicture}
% the grid
\foreach \i in {-1,0,1}
\node (\i-4) at ({\i*sqrt(3)+.5*sqrt(3)},-2) {$\bullet$};
\foreach \i in {0,1}
\node (\i-3) at ({\i*sqrt(3)},{-1.5}) {$\bullet$};
\foreach \i in {-1,0,1}
\node (\i-2) at ({\i*sqrt(3)+.5*sqrt(3)},-1) {$\bullet$};
\foreach \i in {0,1}
\node (\i-1) at ({\i*sqrt(3)},-.5) {$\bullet$};
\foreach \i in {-1,0,1}
\node (\i0) at ({\i*sqrt(3)+.5*sqrt(3)},0) {$\bullet$};
\foreach \i in {0,1}
\node (\i1) at ({\i*sqrt(3)},.5) {$\bullet$};
\foreach \i in {-1,0,1}
\node (\i2) at ({\i*sqrt(3)+.5*sqrt(3)},1) {$\bullet$};
\foreach \i in {0,1}
\node (\i3) at ({\i*sqrt(3)},1.5) {$\bullet$};
\filldraw[orange] ({.5*sqrt(3)},0) circle (3pt);
\filldraw[orange] (0,.5) circle (3pt);
% the arrows
\path[commutative diagrams/.cd, every arrow, every label]  
%     % red arrows, the imperfect linkages 
%             % red vertical
            (00) edge[red, line width = 1, <-]  (02)
            (11) edge[red, line width = 1, <-]  (13)
%             % red right up
            (0-1) edge[red, line width = 1, ->]  (00)
            (0-3) edge[red, line width = 1, ->]  (0-2)
%             % red left up
            (00) edge[red, line width = 1, ->]  (01)
            (0-2) edge[red, line width = 1, ->]  (0-1)
            (0-4) edge[red, line width = 1, ->]  (0-3)
            (11) edge[red, line width = 1, ->]  (02)
            (12) edge[red, line width = 1, ->]  (13)
%     % blue arrows, the perfect linkage
%             % blue vertical 
            (00) edge[blue, line width = 1, ->]  (0-2)
            (0-2) edge[blue, line width = 1, ->]  (0-4)
            (01) edge[blue, line width = 1, ->]  (0-1)
            (0-1) edge[blue, line width = 1, ->]  (0-3)
%             %blue left up
            (01) edge[blue, line width = 1, ->]  (-12)
%             % blue right up 
            (00) edge[blue, line width = 1, ->]  (11)
            (11) edge[blue, line width = 1, ->]  (12)
            (01) edge[blue, line width = 1, ->]  (02)
            (02) edge[blue, line width = 1, ->]  (13)
            ;
\end{tikzpicture}
    \caption{Type II.}
    \label{fig: type_II}
\end{subfigure}
\hfill
\begin{subfigure}{.32\textwidth}
\centering
\begin{tikzpicture}
% the grid
\foreach \i in {-1,0}
\node (\i-4) at ({\i*sqrt(3)+.5*sqrt(3)},-2) {$\bullet$};
\foreach \i in {-1,0,1}
\node (\i-3) at ({\i*sqrt(3)},{-1.5}) {$\bullet$};
\foreach \i in {-1,0}
\node (\i-2) at ({\i*sqrt(3)+.5*sqrt(3)},-1) {$\bullet$};
\foreach \i in {-1,0,1}
\node (\i-1) at ({\i*sqrt(3)},-.5) {$\bullet$};
\foreach \i in {-1,0}
\node (\i0) at ({\i*sqrt(3)+.5*sqrt(3)},0) {$\bullet$};
\foreach \i in {-1,0,1}
\node (\i1) at ({\i*sqrt(3)},.5) {$\bullet$};
\foreach \i in {-1,0}
\node (\i2) at ({\i*sqrt(3)+.5*sqrt(3)},1) {$\bullet$};
\foreach \i in {-1,0,1}
\node (\i3) at ({\i*sqrt(3)},1.5) {$\bullet$};
\filldraw[orange] (0,-.5) circle (3pt);
\filldraw[orange] (0,.5)  circle (3pt);
% the arrows
\path[commutative diagrams/.cd, every arrow, every label]  
%     % red arrows, the imperfect linkages 
%             % red vertical
            (00) edge[red, line width = 1, <-]  (02)
            (11) edge[red, line width = 1, <-]  (13)
            (0-1) edge[red, line width = 1, <-]  (01)
            (-10) edge[red, line width = 1, <-]  (-12)
            (-11) edge[red, line width = 1, <-]  (-13)
%             % red right up
            (-11) edge[red, line width = 1, ->]  (-12)
            (-10) edge[red, line width = 1, ->]  (01)
%             % red left up
            (00) edge[red, line width = 1, ->]  (01)
            (11) edge[red, line width = 1, ->]  (02)
%     % blue arrows, the perfect linkage
%             % blue vertical 
            (0-1) edge[blue, line width = 1, ->]  (0-3)
%             %blue left up
            (01) edge[blue, line width = 1, ->]  (-12)
            (-12) edge[blue, line width = 1, ->]  (-13)
            (0-1) edge[blue, line width = 1, ->]  (-10)
            (-10) edge[blue, line width = 1, ->]  (-11)
%             % blue right up 
            (0-1) edge[blue, line width = 1, ->]  (00)
            (00) edge[blue, line width = 1, ->]  (11)
            (01) edge[blue, line width = 1, ->]  (02)
            (02) edge[blue, line width = 1, ->]  (13)
            ;
\end{tikzpicture}
     \caption{Type III.}
    \label{fig:type_III}
\end{subfigure}
\\
\begin{subfigure}{.45\textwidth}
\centering
\begin{tikzpicture}
% the grid
\foreach \i in {-1,0,1}
\node (\i-3) at ({\i*sqrt(3)},{-1.5}) {$\bullet$};
\foreach \i in {-1,0,1}
\node (\i-2) at ({\i*sqrt(3)+.5*sqrt(3)},-1) {$\bullet$};
\foreach \i in {-1,0,1}
\node (\i-1) at ({\i*sqrt(3)},-.5) {$\bullet$};
\foreach \i in {-1,0,1}
\node (\i0) at ({\i*sqrt(3)+.5*sqrt(3)},0) {$\bullet$};
\foreach \i in {-1,0,1}
\node (\i1) at ({\i*sqrt(3)},.5) {$\bullet$};
\foreach \i in {-1,0,1}
\node (\i2) at ({\i*sqrt(3)+.5*sqrt(3)},1) {$\bullet$};
\foreach \i in {-1,0,1}
\node (\i3) at ({\i*sqrt(3)},1.5) {$\bullet$};
\foreach \i in {-1,0,1}
\node (\i4) at ({\i*sqrt(3)+.5*sqrt(3)},2) {$\bullet$};
\foreach \i in {-1,0,1}
\node (\i5) at ({\i*sqrt(3)},2.5) {$\bullet$};
\filldraw[orange] ({.5*sqrt(3)},1) circle (3pt);
\filldraw[orange] (0,.5) circle (3pt);
\filldraw[orange] (0,1.5) circle (3pt);
% the arrows
\path[commutative diagrams/.cd, every arrow, every label]  
%     % red arrows, the imperfect linkages 
%             % red vertical
            (13) edge[red, line width = 1, <-]  (15)
            (-13) edge[red, line width = 1, <-]  (-15)
            (02) edge[red, line width = 1, <-]  (04)
            (01) edge[red, line width = 1, <-]  (03)
            (-12) edge[red, line width = 1, <-]  (-14)
%             % red right up
            (-13) edge[red, line width = 1, ->]  (-14)
            (-12) edge[red, line width = 1, ->]  (03)
            (01) edge[red, line width = 1, ->]  (02)
            (0-1) edge[red, line width = 1, ->]  (00)
            (0-3) edge[red, line width = 1, ->]  (0-2)
%             % red left up
            (14) edge[red, line width = 1, ->]  (15)
            (13) edge[red, line width = 1, ->]  (04)
            (02) edge[red, line width = 1, ->]  (03)
            (00) edge[red, line width = 1, ->]  (01)
            (0-2) edge[red, line width = 1, ->]  (0-1)
%     % blue arrows, the perfect linkage
%             % blue vertical 
            (00) edge[blue, line width = 1, ->]  (0-2)
            (01) edge[blue, line width = 1, ->]  (0-1)
            (02) edge[blue, line width = 1, ->]  (00)
            (0-1) edge[blue, line width = 1, ->]  (0-3)
%             %blue left up
            (03) edge[blue, line width = 1, ->]  (-14)
            (-14) edge[blue, line width = 1, ->]  (-15)
            (01) edge[blue, line width = 1, ->]  (-12)
            (-12) edge[blue, line width = 1, ->]  (-13)
%             % blue right up 
            (04) edge[blue, line width = 1, ->]  (15)
            (03) edge[blue, line width = 1, ->]  (04)
            (13) edge[blue, line width = 1, ->]  (14)
            (02) edge[blue, line width = 1, ->]  (13)
            ;
\end{tikzpicture}
    \caption{Type IV.}
    \label{fig:type_IV}
\end{subfigure}
\hfill
\begin{subfigure}{.45\textwidth}
\centering
\begin{tikzpicture}
% the grid
\foreach \i in {-1,0,1}
\node (\i-2) at ({\i*sqrt(3)+.5*sqrt(3)},-1) {$\bullet$};
\foreach \i in {-1,0,1}
\node (\i-1) at ({\i*sqrt(3)},-.5) {$\bullet$};
\foreach \i in {-1,0,1}
\node (\i0) at ({\i*sqrt(3)+.5*sqrt(3)},0) {$\bullet$};
\foreach \i in {-1,0,1}
\node (\i1) at ({\i*sqrt(3)},.5) {$\bullet$};
\foreach \i in {-1,0,1}
\node (\i2) at ({\i*sqrt(3)+.5*sqrt(3)},1) {$\bullet$};
\foreach \i in {-1,0,1}
\node (\i3) at ({\i*sqrt(3)},1.5) {$\bullet$};
\foreach \i in {-1,0,1}
\node (\i4) at ({\i*sqrt(3)+.5*sqrt(3)},2) {$\bullet$};
\foreach \i in {-1,0,1}
\node (\i5) at ({\i*sqrt(3)},2.5) {$\bullet$};
\foreach \i in {-1,0,1}
\node (\i6) at ({\i*sqrt(3)+.5*sqrt(3)},3) {$\bullet$};
\filldraw[orange] ({.5*sqrt(3)},1) circle (3pt);
\filldraw[orange] ({.5*sqrt(3)},2) circle (3pt);
\filldraw[orange] (0,1.5) circle (3pt);
% the arrows
\path[commutative diagrams/.cd, every arrow, every label]  
%     % red arrows, the imperfect linkages 
%             % red vertical
            (16) edge[red, line width = 1, ->]  (14)
            (13) edge[red, line width = 1, <-]  (15)
            (05) edge[red, line width = 1, ->]  (03)
            (02) edge[red, line width = 1, <-]  (04)
            (-16) edge[red, line width = 1, ->]  (-14)
%             % red right up
            (-15) edge[red, line width = 1, ->]  (-16)
            (-14) edge[red, line width = 1, ->]  (05)
            (01) edge[red, line width = 1, ->]  (02)
            (0-1) edge[red, line width = 1, ->]  (00)
            (03) edge[red, line width = 1, ->]  (04)
%             % red left up
            (14) edge[red, line width = 1, ->]  (15)
            (13) edge[red, line width = 1, ->]  (04)
            (02) edge[red, line width = 1, ->]  (03)
            (00) edge[red, line width = 1, ->]  (01)
            (0-2) edge[red, line width = 1, ->]  (0-1)
%     % blue arrows, the perfect linkage
%             % blue vertical 
            (00) edge[blue, line width = 1, ->]  (0-2)
            (01) edge[blue, line width = 1, ->]  (0-1)
            (02) edge[blue, line width = 1, ->]  (00)
            (03) edge[blue, line width = 1, ->]  (01)
%             %blue left up
            (03) edge[blue, line width = 1, ->]  (-14)
            (-14) edge[blue, line width = 1, ->]  (-15)
            (05) edge[blue, line width = 1, ->]  (-16)
            (04) edge[blue, line width = 1, ->]  (05)
%             % blue right up 
            (04) edge[blue, line width = 1, ->]  (15)
            (02) edge[blue, line width = 1, ->]  (13)
            (13) edge[blue, line width = 1, ->]  (14)
            (15) edge[blue, line width = 1, ->]  (16)
            ;
\end{tikzpicture}
     \caption{Type V.}
    \label{fig:type_V}
\end{subfigure}
\end{minipage}}
\caption{Non-exact finitely generated linked nets of dimension 1.}
\end{figure}
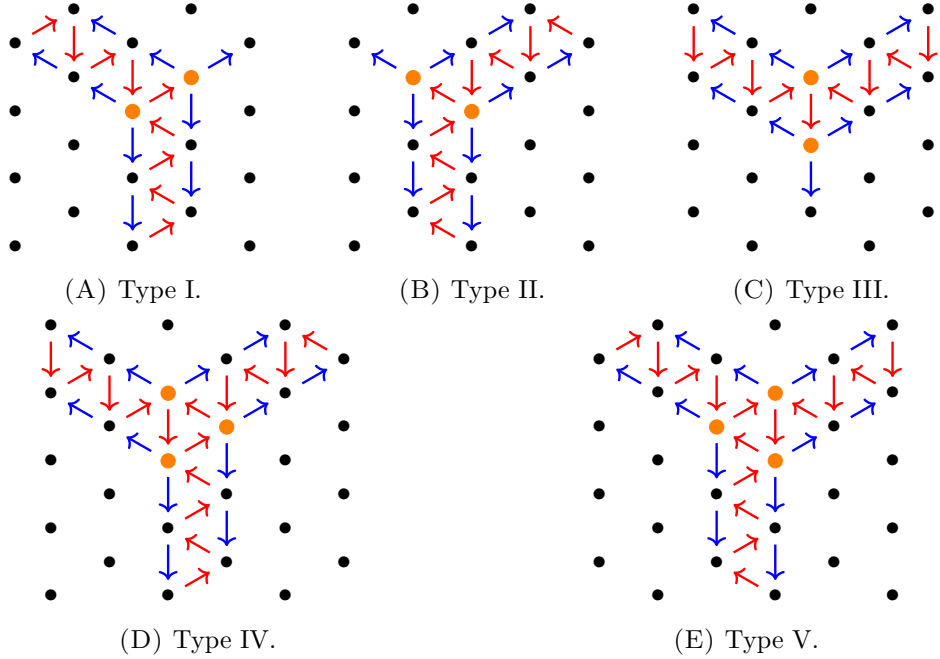

%%THE FIVE TYPES !!!!%%%%%
We show in Figures~\ref{fig: type_I}~to~\ref{fig:type_V} five types of non-exact locally finite simple linked nets over the $\Z^2$-quiver $Q$. An arrow is colored red if the associated map is zero and blue otherwise. To say that a linked net admits a configuration of one of the types shown is to say that there is a finite collection of arrows in the quiver corresponding to maps as indicated by the type. Of course, types I, II and III are the same up to reordering the partition of the arrow set of $Q$, and the same goes for types IV and V. The orange colored vertices are explained below.

\begin{theorem} \label{thm: exact_is_simple_and_non_exact_is_five_types_rank0}
A non-exact locally finite simple linked net over the $\mathbb{Z}^2$-quiver $Q$ admits a configuration of type I, II, III, IV or V. Furthermore, it is minimally generated by the collection of orange vertices indicated in each type. In particular, no linked net admits configurations of two different types.
\end{theorem}

\begin{proof} The third statement follows from the second, as there is a unique minimal collection of vertices $1$-generating the linked net, by Lemma~\ref{minsup1}, and two different types have different collections of orange vertices. 

As for the second statement, observe that for type I, II or III, there are two strips of blue and red arrows meeting at the orange vertices, whereas for type IV or V there are three of them. 
The strips have finite length, but it follows from 
Lemmas~\ref{seq}~and~\ref{triangle} that each of them extends indefinitely away from the orange vertices. Thus, by removing all the red arrows in all the extended strips, we get two connected subquivers for type I, II or III, and three for type IV or V, spanning the whole set of vertices of $Q$. Each subquiver contains a unique orange vertex. 

To prove that the collection of orange vertices $1$-generates the linked net for each type it is enough to prove:
\vskip0.2cm

{\bf Claim:} The restriction of the linked net to each subquiver is generated by the orange vertex it contains. 

Indeed, observe that each orange vertex from which a red arrow $a$ leaves can be connected to each vertex $v$ of the corresponding subquiver by a path $\gamma$ whose essential type is contained in $T-\{\mathfrak a\}$, where $\mathfrak a$ is the type of $a$. Since we have a linked net, $\text{Ker}(\vp_a)\cap\text{Ker}(\vp_\gamma)=0$, and since $\vp_a=0$ we have that $\vp_\gamma$ is a monomorphism, whence an isomorphism. Our claim is proved for the orange vertices we considered, in particular for type IV or V.

On the other hand, given an orange vertex $v$ from which no red arrow leaves, or equivalently whose all arrows leaving it correspond to isomorphisms, let $a_1,a_2,a_3$ denote these arrows, $\mathfrak a_1,\mathfrak a_2,\mathfrak a_3$ their respective types and $v_1,v_2,v_3$ their respective final vertices. Put $I_i:=\{\mathfrak a_1,\mathfrak a_2,\mathfrak a_3\}-\{\mathfrak a_i\}$ for $i=1,2,3$. Notice that we are in type I, II or III, and thus there is one and only additional orange vertex $w$. Up to reordering we may assume the arrow $a$ connecting $w$ to $v$ has type $\mathfrak a_3$. Notice then that there are a red arrow leaving $v_1$ with type $\mathfrak a_2$ and a red arrow leaving $v_2$ with type $\mathfrak a_1$. As before, since we have a linked net, $\vp^{v_1}_u$ is an isomorphism for each $u\in C_{I_2}(v_1)$ and $\vp^{v_2}_u$ is an isomorphism for each $u\in C_{I_1}(v_2)$, whence $\vp^v_u$ is an isomorphism for each $u\in C_{I_2}(v_1)\cup C_{I_1}(v_2)$. Finally, since $\vp^v_w=0$, also $\vp^v_u$ is an isomorphism for each $u\in C_{\mathfrak a_3}(v)$. As the union of the cones $C_{I_2}(v_1)$, $C_{I_1}(v_2)$ and $C_{\mathfrak a_3}(v)$ is the full set of vertices of the subquiver corresponding to $v$, our claim is proved. 

That the collection of orange vertices minimally generates the linked net follows from the fact, which can be ascertained for each type, that $\vp^{v_1}_{v_2}=0$ for each two orange vertices 
$v_1,v_2$.

Finally, we prove the first statement, following the proof we gave to Theorem~\ref{12n}. There are two cases to analyze. Either there is a triangle in $Q$ whose every two vertices are unrelated or not. We consider the first case first.

We claim we have a configuration of type IV or V. By symmetry, we may assume there is a counterclockwise triangle of unrelated vertices. We will show that we have a configuration of type IV. The triangle vertices are depicted in Figure~\ref{fig:proofcase1} in orange and the triangle arrows in red. Since we have a linked net, the map associated to each of the two arrows not in the triangle leaving each of these three orange vertices is a monomorphism, whence an isomorphism. These arrows are depicted in Figure~\ref{fig:proofcase1} in blue. By Lemma~\ref{seq}(1), the arrows of the same type that follow each of these six blue arrows correspond to isomorphisms as well. We depicted them in dashed blue in Figure~\ref{fig:proofcase1}. By the same lemma, the red arrows follow arrows of the same type corresponding to zero maps, depicted in yellow in Figure~\ref{fig:proofcase1}. Using Lemma~\ref{triangle} we obtain that the arrows depicted in dashed red in Figure~\ref{fig:proofcase1} connect unrelated neighbors. We have obtained the configuration of type IV. Notice that the red arrows form indeed a minimal circuit such that all arrows arriving at each vertex of the circuit correspond to zero maps, and hence the orange vertices generate the linked net, as seen in the proof of Theorem~\ref{12n}. 

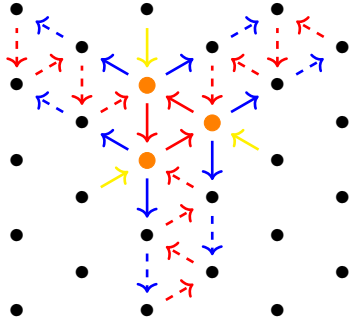
\begin{figure}[ht]
%\begin{subfigure}{.45\textwidth}
\centering
\begin{tikzpicture}
% the grid
\foreach \i in {-1,0,1}
\node (\i-3) at ({\i*sqrt(3)},{-1.5}) {$\bullet$};
\foreach \i in {-1,0,1}
\node (\i-2) at ({\i*sqrt(3)+.5*sqrt(3)},-1) {$\bullet$};
\foreach \i in {-1,0,1}
\node (\i-1) at ({\i*sqrt(3)},-.5) {$\bullet$};
\foreach \i in {-1,0,1}
\node (\i0) at ({\i*sqrt(3)+.5*sqrt(3)},0) {$\bullet$};
\foreach \i in {-1,0,1}
\node (\i1) at ({\i*sqrt(3)},.5) {$\bullet$};
\foreach \i in {-1,0,1}
\node (\i2) at ({\i*sqrt(3)+.5*sqrt(3)},1) {$\bullet$};
\foreach \i in {-1,0,1}
\node (\i3) at ({\i*sqrt(3)},1.5) {$\bullet$};
\foreach \i in {-1,0,1}
\node (\i4) at ({\i*sqrt(3)+.5*sqrt(3)},2) {$\bullet$};
\foreach \i in {-1,0,1}
\node (\i5) at ({\i*sqrt(3)},2.5) {$\bullet$};
\filldraw[orange] ({.5*sqrt(3)},1) circle (3pt);
\filldraw[orange] (0,.5) circle (3pt);
\filldraw[orange] (0,1.5) circle (3pt);
% the arrows
\path[commutative diagrams/.cd, every arrow, every label]  
%     % red arrows, the imperfect linkages 
%             % red vertical
            (05) edge[yellow, line width = 1, ->]  (03)
            (13) edge[red, line width = 1, <-, dashed]  (15)
            (-13) edge[red, line width = 1, <-, dashed]  (-15)
            (02) edge[red, line width = 1, <-, dashed]  (04)
            (01) edge[red, line width = 1, <-]  (03)
            (-12) edge[red, line width = 1, <-, dashed]  (-14)
%             % red right up
            (11) edge[yellow, line width = 1, ->]  (02)
            (-13) edge[red, line width = 1, ->, dashed]  (-14)
            (-12) edge[red, line width = 1, ->, dashed]  (03)
            (01) edge[red, line width = 1, ->]  (02)
            (0-1) edge[red, line width = 1, ->, dashed]  (00)
            (0-3) edge[red, line width = 1, ->, dashed]  (0-2)
%             % red left up
            (-10) edge[yellow, line width = 1, ->]  (01)
            (14) edge[red, line width = 1, ->, dashed]  (15)
            (13) edge[red, line width = 1, ->, dashed]  (04)
            (02) edge[red, line width = 1, ->]  (03)
            (00) edge[red, line width = 1, ->, dashed]  (01)
            (0-2) edge[red, line width = 1, ->, dashed]  (0-1)
%     % blue arrows, the perfect linkage
%             % blue vertical 
            (00) edge[blue, line width = 1, ->, dashed]  (0-2)
            (01) edge[blue, line width = 1, ->]  (0-1)
            (02) edge[blue, line width = 1, ->]  (00)
            (0-1) edge[blue, line width = 1, ->, dashed]  (0-3)
%             %blue left up
            (03) edge[blue, line width = 1, ->]  (-14)
            (-14) edge[blue, line width = 1, ->, dashed]  (-15)
            (01) edge[blue, line width = 1, ->]  (-12)
            (-12) edge[blue, line width = 1, ->, dashed]  (-13)
%             % blue right up 
            (04) edge[blue, line width = 1, ->, dashed]  (15)
            (03) edge[blue, line width = 1, ->]  (04)
            (13) edge[blue, line width = 1, ->, dashed]  (14)
            (02) edge[blue, line width = 1, ->]  (13)
            ;
\end{tikzpicture}
    \caption{Proof of Theorem~\ref{thm: exact_is_simple_and_non_exact_is_five_types_rank0}, Case 1}
    \label{fig:proofcase1}
%\end{subfigure}
\end{figure}

We may now assume there is no triangle of unrelated vertices. We claim we have a configuration of type I, II or III. By hypothesis, in each triangle there is at least one arrow that corresponds to an isomorphism. There might be triangles with two arrows corresponding to isomorphisms, but there is at least one triangle with only one arrow corresponding to an isomorphism, as there are unrelated vertices. Consider such a triangle. Without loss of generality we may assume that the orientation of the triangle is clockwise and the arrow corresponding to an isomorphism is of type 2. We will show we have a configuration of type I or II. The triangle in question is the clockwise triangle depicted in Figure~\ref{fig:proofcase2} containing the two orange vertices. The arrow of type 2 is depicted in blue and the other arrows in red.

Since we have a linked net, the arrow of type 2 leaving from the other orange vertex corresponds to an isomorphism. By Lemma~\ref{seq}(1), so does each arrow in the path of essential type 2 leaving from each orange vertex. These arrows are depicted in dashed blue in Figure~\ref{fig:proofcase2}. By Lemma~\ref{triangle}, the arrows that connect one vertex from one path to the other correspond to zero maps. These arrows are depicted in dashed red in Figure~\ref{fig:proofcase2}.

Notice that each of the triangles with dashed red and blue arrows lies on a strip below the initial triangle with red and blue arrows. Each triangle lying below the inital triangle has a vertex with a blue arrow, of type 2, arriving at it. If we do a 2-shift to that triangle, centered at that vertex, as explained in the proof of Theorem~\ref{12n}, we end up with the triangle with the same orientation above it. We have seen in that proof that one cannot do 2-shifts indefinitely. In the case at hand, if there were a triangle on the extended doubly infinite strip, with any orientation, above the initial triangle, of the same sort, that is, with the arrow of type 2 being the unique one corresponding to an isomorphism, then we would likewise conclude that all triangles below it on the same strip would be of the same sort. It is not possible however to have all the triangles on the whole doubly infinite strip of the same sort as the initial triangle, since $\g$ is locally finite. Thus there is a topmost triangle on that strip of the sort we are considering. Assume it is a clockwise triangle, as depicted in Figure~\ref{fig:proofcase2}. As it is the topmost such triangle, the yellow arrows correspond to zero maps. And since there is no triangle of unrelated vertices, the green arrow corresponds to an isomorphism. We will show we have a configuration of type I. 

Indeed, since we have a linked net, the arrow of type 1 leaving from the leftmost orange vertex corresponds to an isomorphism. By Lemma~\ref{seq}(1), so does each arrow in the path of essential type 1 leaving from each orange vertex. These arrows are depicted in dashed green in Figure~\ref{fig:proofcase2}. By Lemma~\ref{triangle}, the arrows that connect one vertex from one path to the other correspond to zero maps. These arrows are depicted in dashed yellow in Figure~\ref{fig:proofcase2}. Finally, the dashed black arrow corresponds to an isomorphism by Lemma~\ref{seq}(2), as the orange vertices are unrelated. We have obtained the configuration of type I. 

Observe that, since we have a linked net, the arrow of type 0 leaving the final vertex of the blue arrow in Figure~\ref{fig:proofcase2} corresponds to an isomorphism, and thus the arrow of type 1 following it must correspond to the zero map. It follows that the arrows in the clockwise triangle containing the orange vertices form a minimal circuit such that all arrows arriving at each orange vertex correspond to zero maps, and hence the orange vertices generate the linked net, as seen in the proof of Theorem~\ref{12n}. 

\begin{figure}[ht]
\centering
\begin{tikzpicture}
% the grid
\foreach \i in {-1,0,1}
\node (\i-3) at ({\i*sqrt(3)},{-1.5}) {$\bullet$};
\foreach \i in {-1,0}
\node (\i-2) at ({\i*sqrt(3)+.5*sqrt(3)},-1) {$\bullet$};
\foreach \i in {-1,0,1}
\node (\i-1) at ({\i*sqrt(3)},-.5) {$\bullet$};
\foreach \i in {-1,0}
\node (\i0) at ({\i*sqrt(3)+.5*sqrt(3)},0) {$\bullet$};
\foreach \i in {-1,0,1}
\node (\i1) at ({\i*sqrt(3)},.5) {$\bullet$};
\foreach \i in {-1,0}
\node (\i2) at ({\i*sqrt(3)+.5*sqrt(3)},1) {$\bullet$};
\foreach \i in {-1,0,1}
\node (\i3) at ({\i*sqrt(3)},1.5) {$\bullet$};
\foreach \i in {-1,0}
\node (\i4) at ({\i*sqrt(3)+.5*sqrt(3)},2) {$\bullet$};
\filldraw[orange] ({.5*sqrt(3)},1) circle (3pt);
\filldraw[orange] (0,.5) circle (3pt);
% the arrows
\path[commutative diagrams/.cd, every arrow, every label]  
%     % red arrows, the imperfect linkages 
%             % red vertical
            (04) edge[yellow, line width = 1, ->] (02)
            (01) edge[yellow, line width = 1, <-]  (03)
            (-12) edge[yellow, line width = 1, <-, dashed]  (-14)
%             % red right up
            (-13) edge[yellow, line width = 1, ->, dashed]  (-14)
            (-12) edge[yellow, line width = 1, ->, dashed]  (03)
            (01) edge[red, line width = 1, ->]  (02)
            (0-1) edge[red, line width = 1, ->, dashed]  (00)
            (0-3) edge[red, line width = 1, ->, dashed]  (0-2)
%             % red left up
            (00) edge[red, line width = 1, ->]  (01)
            (0-2) edge[red, line width = 1, ->, dashed]  (0-1)
%     % blue arrows, the perfect linkage
%             % blue vertical 
            (00) edge[blue, line width = 1, ->, dashed]  (0-2)
            (01) edge[blue, line width = 1, ->, dashed]  (0-1)
            (02) edge[blue, line width = 1, ->]  (00)
            (0-1) edge[blue, line width = 1, ->, dashed]  (0-3)
%             %blue left up
            (03) edge[green, line width = 1, ->, dashed]  (-14)
            (02) edge[green, line width = 1, ->]  (03)
            (01) edge[green, line width = 1, ->, dashed]  (-12)
            (-12) edge[green, line width = 1, ->, dashed]  (-13)
%             % blue right up 
            (02) edge[line width = 1, ->, dashed]  (13)
            ;
\end{tikzpicture}
    \caption{Proof of Theorem~\ref{thm: exact_is_simple_and_non_exact_is_five_types_rank0}, Case 2}
    \label{fig:proofcase2}
\end{figure}

\end{proof}

%\begin{theorem}\label{123} A linked net of vector spaces of dimension $1$ with finite support over a $\mathbb Z^2$-quiver is generated by a triangle.
%\end{theorem}
%\begin{proof} 
%Since the linked net of the statement has finite support it has one of the five configuration from I to V of Theorem~\ref{thm: exact_is_simple_and_non_exact_is_five_types_rank0}. That is, the linked net is generated by a triangle. 
%\end{proof}

\section{The linked projective space}\label{lpg}

%This section is dedicated to the study of $\lp(\g)$, where $\mathfrak{g}$ is a linked net of vector spaces over a $\mathbb{Z}^2$-quiver. We determine the structure of $\lp(\g)$, in the same spirit of the work by Santana in \cite{Rocha2019}. In summary, in Theorem \ref{Main_Theorem} we prove that the exact points of $\lp(\g)$ form an open set $\lp(\g)^{*}$, which is equal to the nonsingular locus. Therefore, the non-exact points are all singular and we know how to classify them: They are precisely the five types described in Theorem~\ref{thm: exact_is_simple_and_non_exact_is_five_types_rank0}.

Recall that for each vector space $V$ and vector $s\in V$, we denote by $[s]$ the class of its nonzero scalar multiples, that is, $[s]:=\{cs\,|\,c\in k^*\}$. If $s=0$, we write $[s]=0$. 
If $\vp\colon V\to W$ is a map of vector spaces, we let 
$[\vp][s]:=[\vp(s)]$. Given another vector $t\in V$ we let $[s]\wedge [t]:=[s\wedge t]$. When we write $[s]\in\mathbb P(V)$ we assume implicitly that $s\neq 0$.

%Given a linked net $\g$ of vector spaces over $k$, the net is pure if and only if the associated vector spaces have the same finite dimension.

%\begin{Def} The \emph{dimension} of a pure linked net of vector spaces $\g$, denoted $\dim\g$, is the dimension over $k$ of any associated vector space.
%\end{Def}

Let $\g$ be a representation of a $\mathbb Z^n$-quiver $Q$ in the category of nontrivial finite-dimensional vector spaces. Let $\lp(\g)$ be the quiver Grassmannian of subrepresentations of pure dimension~1 of $\g$. It is a set. Assume $\g$ is a weakly linked net. For each finite collection $H$ of vertices of $Q$, let $\lp_H(\g)$ be the subscheme of $\prod_{v\in H}\mathbb P(V_v)$ defined by
$$
\lp_H(\g):=\Big\{(s_v\,|\,v\in H)\in\prod_{v\in H}\mathbb P(V_v) \,\Big|\, 
\vp^v_w(s_v)\wedge s_w=0 \text{ for all }v,w\in H\Big\}.
$$
There is a natural map $\Psi^\g_H\colon\lp(\g)\to\lp_H(\g)$, induced by restriction. If $\g$ is $1$-generated by $H$ then $\Psi^\g_H$ is clearly injective. It is bijective if in addition $\g$ is pure and $P(H)=H$, as we will see below.
%If $H$ is well chosen, $\Psi_H$ is a bijection, as stated by the next proposition.

Recall from \cite{Esteves_et_al_2021} that the \emph{hull} of a set $H$ of vertices of a $\mathbb Z^n$-quiver $Q$ is the set $P(H)$ of all vertices $v$ of $Q$ such that for each arrow type there are 
$z\in H$ and a path $\gamma$ connecting $z$ to $v$ not containing any arrow of that type. Of course, $H\subseteq P(H)$. By \cite{Esteves_et_al_2021}, Prop.~5.6, if $H$ is finite, so is $P(H)$. Also, $P(P(H))=P(H)$. Hence, every finitely generated linked net is $1$-generated by a finite set of vertices equal to its hull.

By \cite{Esteves_et_al_2021}, Prop.~5.7, for each vertex $v$ of $Q$ there is $w_v\in P(H)$ such that for each $z\in H$ there is an admissible path connecting $z$ to $v$ through $w_v$. Furthermore, $w_v$ is unique if $P(H)=H$.

\begin{Def} Let $Q$ be a $\mathbb Z^n$-quiver and $H$ a nonempty collection of vertices of $Q$ such that $P(H)=H$. For each vertex $v $ of $Q$ we call the unique vertex $w_v\in H$ for which there is an admissible path connecting $z$ to $v$ through $w_v$ for each $z\in H$ the \emph{shadow} of $v$ in $H$. 
\end{Def}

\begin{proposition}\label{HHPPHH} Let $\g$ be a pure nontrivial weakly linked net over a $\mathbb Z^n$-quiver $Q$ of vector spaces over $k$. Let $H_1,H_2,H_3$ be finite collections of vertices of $Q$ $1$-generating $\g$ with $P(H_1)=H_1$. Let 
$$
\widetilde\Psi^{H_1}_{H_2}\colon\prod_{v\in H_1}\mathbb P(V_v)\longrightarrow \prod_{v\in H_2}\mathbb P(V_v)
$$
sending $(s_v\,|\,v\in H_1)$ to $(t_v\,|\,v\in H_2)$ satisfying $t_v=\varphi^{w_v}_v(s_{w_v})$ for each $v\in H_2$, where $w_v$ is the shadow of $v$ in $H_1$. 
Then $\widetilde\Psi^{H_1}_{H_2}$ is a well-defined scheme morphism and restricts to a morphism $\Psi^{H_1}_{H_2}\colon\lp_{H_1}(\g)\to\lp_{H_2}(\g)$. Furthermore,  
\begin{enumerate}
    \item $\Psi^\g_{H_1}$ is bijective,
    \item $\Psi^\g_{H_2}=\Psi^{H_1}_{H_2}\Psi^\g_{H_1}$, 
    \item $\Psi^{H_2}_{H_3}\Psi^{H_1}_{H_2}=\Psi^{H_1}_{H_3}$ if $P(H_2)=H_2$,
    \item $\Psi^{H_1}_{H_2}$ is an isomorphism of schemes if $P(H_2)=H_2$.
\end{enumerate}
\end{proposition}

\begin{proof} For each vertex $v$ of $Q$ let $w_v\in H_1$ be its shadow. %Notice that $w_v=v$ if $v\in H_1$. 
Since $\g$ is $1$-generated by $H_1$ there is
$z\in H_1$ such that $\vp^z_v$ is an isomorphism. Since there is an admissible path from $z$ to $v$ through $w_v$, we have that $\vp^z_v= \vp^{w_v}_v \vp^z_{w_v}$, and hence 
$\vp^{w_v}_v$ is an epimorphism, thus an isomorphism because $\g$ is pure. It follows that $\widetilde\Psi^{H_1}_{H_2}$ is a well-defined scheme morphism.

Furthermore, we claim $\widetilde\Psi^{H_1}_{H_2}$ takes 
$\lp_{H_1}(\g)$ to $\lp_{H_2}(\g)$. Indeed, for each two vertices $u$ and $v$ of $Q$, there is an admissible path connecting $w_u$ to $v$ through $w_v$. Thus
$\vp^{w_u}_{v}=\vp^{w_v}_{v} \vp^{w_u}_{w_v}$. Furthermore, if $\vp^u_v$ is nonzero, then so is $\vp^u_v \vp^{w_u}_u$, and hence $\vp^{w_u}_{v}=\vp^u_v \vp^{w_u}_u$ as well. Thus the equation 
$$
\vp^u_v(\vp^{w_u}_u(s_{w_u}))\wedge \vp^{w_v}_{v}(s_{w_v})=0
$$
holds trivially if $\vp^u_v=0$ and follows otherwise from 
$$
\vp^{w_u}_{w_v}(s_{w_u})\wedge s_{w_v}=0
$$
by applying $\vp^{w_v}_{v}$ to both sides. It follows that $\widetilde\Psi^{H_1}_{H_2}$ restricts to a morphism  $\Psi^{H_1}_{H_2}\colon\lp_{H_1}(\g)\to\lp_{H_2}(\g)$. 

Furthermore, the last argument, as it applies to all vertices $u,v$ of $Q$, shows as well that $\Psi^\g_{H_1}$ is surjective. Since $\Psi^\g_{H_1}$ is injective, it follows that $\Psi^\g_{H_1}$ is bijective and 
$\Psi^\g_{H_2}=\Psi^{H_1}_{H_2}\Psi^\g_{H_1}$, proving Statements~(1)~and (2).

Clearly, $\widetilde\Psi^{H_1}_{H_1}$ is the identity. 
Then Statement~(4) follows from Statement~(3).

It remains to prove Statement~(3). For each vertex $v\in H_3$, let $u_v$ be its shadow in $H_2$, and $w'_v$ the shadow of $u_v$ in $H_1$. Let $w_v$ be the shadow of $v$ in $H_1$. Then, as we have seen, $\vp^{u_v}_v$, $\vp^{w'_v}_{u_v}$ and $\vp^{w_v}_v$ are isomorphisms. Then $\vp^{u_v}_v\vp^{w'_v}_{u_v}$ is an isomorphism, hence nonzero since $\g$ is pure and nontrivial. So 
\begin{equation}\label{wvuv}
\vp^{w'_v}_v=\vp^{u_v}_v\vp^{w'_v}_{u_v}.
\end{equation}
On the other hand, there is an admissible path from $w'_v$ to $v$ through $w_v$, whence \begin{equation}\label{wvwv}
\vp^{w'_v}_v=\vp^{w_v}_v\vp^{w'_v}_{w_v}.
\end{equation}

The map $\widetilde\Psi^{H_1}_{H_3}$ takes $(s_z\,|\,z\in H_1)$ to $(\vp^{w_v}_v(s_{w_v})\,|\,v\in H_3)$, whereas it follows from \eqref{wvuv} that
$\widetilde\Psi^{H_2}_{H_3}\widetilde\Psi^{H_1}_{H_2}$ takes 
$(s_z\,|\,z\in H_1)$ to $(\vp^{w'_v}_v(s_{w'_v})\,|\,v\in H_3)$. 
Since $\vp^{w'_v}_{w_v}(s_{w'_v})\wedge s_{w_v}=0$ on $\lp_{H_1}(\g)$, applying 
$\vp^{w_v}_v$ to both sides and using \eqref{wvwv} we obtain
$$
\vp^{w'_v}_{v}(s_{w'_v})\wedge\vp^{w_v}_v(s_{w_v})=0
$$
for each $v\in H_3$, and thus 
$\Psi^{H_2}_{H_3}\Psi^{H_1}_{H_2}=\Psi^{H_1}_{H_3}$.
\end{proof}

\begin{Def} Let $\g$ be a pure nontrivial finitely generated weakly linked net of vector spaces over a $\mathbb Z^n$-quiver $Q$. Give $\lp(\g)$ the scheme structure induced from the bijection 
$\psi^\g_H\colon\lp(\g)\to\lp_H(\g)$ for a finite set $H$ of vertices of $Q$ generating $\g$ and satisfying $P(H)=H$. We call $\lp(\g)$ the 
\emph{linked projective space} associated to $\g$. We say it \emph{has the Hilbert polynomial of the diagonal} if so has $\lp_H(\g)$.
\end{Def}

It follows from Proposition~\ref{HHPPHH} that $\psi_H$ is a scheme morphism for each finite subset of vertices $H$ that $1$-generates $\g$ and that the scheme structure on $\lp(\g)$ does not depend on the choice of $H$. Rather, the choice of $H$ with $P(H)=H$ gives us an embedding $\lp(\g)\hookrightarrow\prod_{v\in H}\mathbb P(V_v)$. But even the extrinsic structures on $\lp(\g)$ given by different $H$ are somewhat comparable, because the isomorphisms between the $\lp_H(\g)$ are restrictions of linear maps on their ambient spaces. So, for instance, the multivariate Hilbert polynomial of $\lp_{H_1}(\g)$ is that of the (small) diagonal if and only if so is the multivariate Hilbert polynomial of $\lp_{H_2}(\g)$. Indeed, we may assume that $H_2\supseteq H_1$, and in this case, the multivariate Hilbert polynomial of the latter, $\text{Hilb}_{\lp_{H_2}(\g)}(n_v\,|\,v\in H_2)$, is obtained from that of the former, $\text{Hilb}_{\lp_{H_1}(\g)}(n_u\,|\,u\in H_1)$, by replacing each $n_u$ for $u\in H_1$ by the sum of the $n_v$ for all $v\in H_2$ such that $w_v=u$. 

%We will often 
% describe a point of $\lp(\g)$ as a collection $(U_v\,|\,v\in H)$ of subspaces $U_v$ of dimension 1 of the $V_v$, instead of $([s_v]\,|\,v\in H)$, the relation being that $U_v$ is generated by $s_v$ for each $v$. We will also 

Given that a point on $\lp(\g)$ corresponds to a weakly linked net $\h$ of vector spaces over $Q$, we attribute to the point adjectives we attribute to $\h$. For instance, the point is exact if $\h$ is exact. We will also write $\h\in\lp(\g)$.

\begin{Def} Let $\g$ be a pure nontrivial finitely generated weakly linked net of vector spaces over a $\mathbb Z^n$-quiver $Q$. For each vertex $v$ of $Q$, let
\begin{displaymath}
\lp(\g)_{v}^{*} := \{ \mathfrak W\subseteq\mathfrak V\,|\,\h \text{ is generated by }v\}
\end{displaymath}
and put $\lp(\g)_v := \overline{\lp(\g)_v^{*}}$.
\end{Def}

\begin{proposition}\label{birLPv} Let $\g$ be a pure nontrivial finitely generated weakly linked net of vector spaces over a $\mathbb Z^n$-quiver $Q$. Let $v$ be a vertex of $Q$. Then $\lp(\g)_v^*$ is a nonsingular open subscheme of $\lp(\g)$. If nonempty, there is a birational map
$$
\mathbb P(V_v)\longrightarrow \lp(\g)_v.
$$
In particular, each nonempty $\lp(\g)_v$ is irreducible of dimension~$\dim\g-1$ and rational.
\end{proposition}

\begin{proof} Let $H$ be a finite set of vertices containing $v$ and $1$-generating $\g$. Then $\lp(\g)_{v}^{*}=(\psi_H^\g)^{-1}(U)$, where $U$ is the set of $(s_u\,|\,u\in H)$ such that $\vp^v_u(s_v)\neq 0$ for each $u\in H$, thus open in $\lp_H(\g)$. The rational map is naturally defined by taking $[s]\in \mathbb{P}(V_v)$ to
the subnet $\h\subseteq\g$ generated by $ks$. It is defined on the open set $U'$ of $\mathbb P(V_v)$ parameterizing the $[s]$ for which $\vp^v_u(s)\neq 0$ for each $u\in H$, with image $\lp(\g)^*_v$. The composition $U'\to\lp(\g)_v^*\to U$ has a natural inverse, induced by projection. Taking $H$ such that $H=P(H)$, we get an isomorphism $U'\to\lp(\g)_v^*$.
\end{proof}

%All $\h\in\lp(\g)^*_v$ are exact by \cite[Thm.~7.6]{Esteves_et_al_2021}. And in fact, it follows from loc.~cit., Thm.~7.8,  
% Theorem labeled as ``\label{simple1}" in article I
%that every exact point on $\lp(\g)$ belongs to $\lp(\g)_v^{*}$ for some $v$.

%\begin{Rem} 
%If $\vp^v_w\colon V_v\to V_w$ is zero then $\lp(\g)_v = \emptyset$. 
%\end{Rem}

\begin{proposition}\label{clsdzero} Let $\g$ be a pure nontrivial finitely generated weakly linked net of vector spaces over a $\mathbb Z^n$-quiver $Q$. Let $\h\in\lp(\g)$. Let $H$ be a set of vertices of $Q$ that $1$-generates $\h$ and $v$ be a vertex of $Q$. If $\h\in\lp(\g)_v$ then $v\in H$ and $\varphi^u_v(V^\h_u)=0$ for each other vertex $u$ of $Q$. In particular, $\lp(\g)_u^{*} \cap \lp(\g)_v^{*}=\emptyset$.
\end{proposition}

\begin{proof} If $\h\in\lp(\g)^*_v$ then $\vp^v_u(V^\h_v)=V^\h_u$ and hence 
$\h\not\in\lp(\g)^*_u$ for each vertex $u$ distinct from $v$ because $\vp^u_v(V^\h_u)=0$. The latter is a closed condition, and thus holds as well if $\h\in\lp(\g)_v$. Since $\vp^z_v(V^\h_z)$ is nonzero for some $z\in H$, it follows that $v\in H$.
\end{proof}

\section{Linked projective spaces of exact linked nets}\label{lpg_z2}

%Assume in this section that $n=2$, that is, that $Q$ is a $\Z^2$-quiver.

\begin{Def} Let $\g$ be a pure nontrivial finitely generated weakly linked net of vector spaces over a $\mathbb{Z}^n$-quiver $Q$. For each finite subset of vertices $H$ of $Q$, put
$$
\lp(\g)_H:=\bigcap_{v\in H}\lp(\g)_v\quad\text{and}\quad
\lp(\g)^*_H:=\lp(\g)_H-\bigcup_{v\not\in H}\lp(\g)_v.
$$
\end{Def}

\begin{proposition}\label{lemma: non_exact_point_is_in_frontier}
  Let $\g$ be a pure nontrivial exact finitely generated linked net of vector spaces over a $\mathbb{Z}^n$-quiver $Q$. Let $H$ be a finite set of vertices of $Q$. Then 
  \begin{equation}\label{lpHmin}
  \lp(\g)^*_H=\Big\{\h\in\lp(\g)\,\Big|\,\h\text{ is minimally $1$-generated by }H\Big\}.
  \end{equation}
 Furthermore, $\lp(\g)^*_H$ is open and dense in $\lp(\g)_H$.
\end{proposition}

\begin{proof} %Let $H$ be a finite set of vertices $1$-generating $\g$. 
Let $\h\in\lp(\g)$. By Theorem~\ref{12n}, there are vertices $u_1,\dots,u_r$ forming an oriented polygon $\Delta$ 
minimally $1$-generating $\h$. For each $i=1,\dots,r$, let 
$s_i$ be a generator of $V^\h_{u_i}$. For each $i,j\in\{1,\dots,r\}$, let $\psi^i_j$ be a map representing 
$\vp^{u_i}_{u_j}$. For convenience, put $u_{r+1}:=u_1$ and $s_{r+1}:=s_1$, and let
$\psi^{r+1}_j:=\psi^1_j$ for each $j$.
%For convenience, put $u_{qr+j}:=u_j$ and $s_{qr+j}:=s_j$ for each $q\in\mathbb Z$ and $j=1,\dots,r$. Likewise, put $\psi^{q_1r+i}_{q_2r+_j}:=\psi^i_j$ for for each $q_1,q_2\in\mathbb Z$ and $i,j\in\{1,\dots,r\}$. 
By Lemma~\ref{minsup1}, we have $\psi^{i}_j(s_i) = 0$ for each distinct $i,j$. Since $\g$ is exact, there is $s^i\in V_{u_i}^\g$ for each $i=1,\dots,r$ such that  
$s_{i} = \psi^{i+1}_{i}(s^{i+1})$ for each $i$, where for convenience we put $s^{r+1}:=s_1$.

We claim that $\h\in\lp(\g)_\Delta$, or equivalently, $\h\in\lp(\g)_{u_i}$ for each $i$. Reordering the vertices $u_j$ if necessary, it is enough to show that $\h\in\lp(\g)_{u_1}$. For each $t \in k$, let $\h^t$ be the subnet of $\g$ generated by  
$$
s^t:=\sum_{\ell=2}^{r+1}t^{\ell-2}\psi^{\ell}_1(s^\ell)\in V^\g_{u_1}.
$$
That $\h^t$ is indeed a subnet of $\g$ follows from \cite{Esteves_et_al_2021}, Prop.~7.1.
To prove the claim we will show that $\h^t\in\lp(\g)^*_{u_1}$ for a general $t$ and $\h=\lim_{t\to 0}\h^t$. Indeed, since $\g$ is pure and finitely generated, it is enough to show that for each vertex $z$, the space $V^{\h^t}_z$ is nonzero for a general $t\in k$ and 
$V^{\h}_z=\lim_{t\to 0}V^{\h^t}_z$. Now, for each vertex $z$ there is 
$j\in\{1,\dots,r\}$ such that $u_j$ is the shadow of $z$ in $\Delta$. Then $V^\h_z=\vp^{u_j}_z(ks_j)$ and $V^{\h^t}_z=\vp^{u_j}_z(k\psi^1_j(s^t))$. 
But
$$
%\begin{align*}
\psi^1_j(s^t)=\sum_{\ell=2}^{r+1}t^{\ell-2}\psi^1_j\psi^{\ell}_1(s^\ell)
=\sum_{\ell=j+1}^{r+1}t^{\ell-2}\psi^{\ell}_j(s^\ell)
=t^{j-1}s_j+\sum_{\ell=j+2}^{r+1}t^{\ell-2}\psi^{\ell}_j(s^\ell).
%\end{align*}
$$
Since $\vp^{u_j}_z(ks_j)\neq 0$, we have $\vp^{u_j}_z(k\psi^1_j(s^t))\neq 0$ for general $t$, and thus $V^{\h^t}_z\neq 0$, as wished. Furthermore, $V^{\h}_z=\lim_{t\to 0}V^{\h^t}_z$, finishing the proof of the claim.
%In particular,
%\begin{equation}\label{uu}
%\h\in\lp(\g)_{u_1}\cap\cdots\cap\lp(\g)_{u_r}.
%\end{equation}

Moreover, it follows from Proposition~\ref{clsdzero} that 
$\h\in\lp(\g)^*_\Delta$. Thus, if $\h$ is minimally $1$-generated by $H$ then $H=\Delta$ by Lemma~\ref{minsup1}, and hence 
$\h\in\lp(\g)^*_H$. 

%If $\h\in\lp(\g)_H$, then Proposition~\ref{clsdzero} yields $H\subseteq\Delta$. Under this condition, $\h$ is minimally $1$-generated by $H$ if and only if $\h$ is $1$-generated by $H$.
%, or equivalently, for each vertex $a$ of $Q$, there is $v\in H$ such that $\vp^{v}_a(s_v)\neq 0$. As we need only check this for the vertices $a$ of a finite set $1$-generating $\g$, it follows that 
%those $\h\in\lp(\g)_H$ (minimally) $1$-generated by $H$ form an open subset. 

Conversely, if $\h\in\lp(\g)_H$ then 
$H\subseteq\Delta$ by Proposition~\ref{clsdzero} again. Moreover, if $\h\in\lp(\g)_H^*$ then $\h\not\in\lp(\g)_z$ for each 
$z\in\Delta-H$. But $\h\in\lp(\g)_\Delta$ by our claim, whence 
$\Delta=H$. Thus $\h$ is minimally $1$-generated by $H$. We have proved the first statement of the proposition. 

As for the second statement, since $\g$ is finitely generated, only finitely many $\lp(\g)_v$ are nonempty, and thus $\lp(\g)^*_H$ is clearly open in $\lp(\g)_H$. We have to prove it is dense. 

Assume $\h\in\lp(\g)_H$. We will prove that $\h$ is in the closure of $\lp(\g)^*_H$. Now, $H\subseteq\Delta$ by Proposition~\ref{clsdzero}. So there is a subsequence $p_1,\dots,p_m$ of $1,\dots,r$ such that $H=\{u_{p_1},\dots,u_{p_m}\}$. Up to reordering the $u_i$ we may assume $p_1=1$. For convenience, put $p_{m+1}:=r+1$.
%For convenience put $v_i:=u_{\ell_i}$ for each $i$. 
%Since $\h$ is minimally generated by $\{u_1,\dots,u_r\}$ we have that $\vp^{u_i}_{u_j}(s_{u_i})=0$ for each distinct $i,j$. Since $\g$ is exact, for each $i=1,\dots,m$ and each $j=\ell_i+1,\dots,\ell_{i+1}-1$ there is $s^i_j\in V^\g_{v_i}$ such that $\vp^{v_i}_{u_j}(s^i_j)=s_{u_j}$. 
For each $t\in K$, let $\mathfrak U^t$ be the subnet of $\g$ generated by 
$s^t_1,\dots,s^t_{m}$, where 
$$
s^t_{i}:=\sum_{p_i<\ell\leq p_{i+1}}t^{\ell-p_i-1}\psi^\ell_{p_i}(s^\ell)\in V^\g_{u_{p_i}}
$$ 
for $i=1,\dots,m$. That $\mathfrak U^t$ is indeed a subnet of $\g$ follows from \cite{Esteves_et_al_2021}, Prop.~7.1. 

Now, $H$ is a polygon because $H\subseteq\Delta$, and thus $P(H)=H$ by \cite{Esteves_et_al_2021}, Prop.~5.10. Since $\mathfrak U^t$ is generated by $H$ it follows from \cite{Esteves_et_al_2021}, Prop.~6.4, that $\mathfrak U^t$ is 1-generated by $H$ for each $t$. It is minimally so, because $\psi^{p_i}_{p_j}(s^t_{i})=0$ for each distinct $i,j$. If we show that $\mathfrak U^t\in\lp(\g)$ for general $t$ and 
$\h=\lim_{t\to 0}\mathfrak U^t$, it will thus follow, as we have seen for $\h$, that $\mathfrak U^t\in\lp(\g)^*_H$ and hence that 
$\h$ lies in the closure of $\lp(\g)^*_H$, as wished.

As before, for each vertex $z$ of the quiver, there is 
$j\in\{1,\dots,r\}$ such that $u_j$ is the shadow of $z$ in $\Delta$. Then 
$$
V^{\h}_z=\vp^{u_j}_z(ks_j)\quad\text{and}\quad 
V^{\mathfrak U^t}_z=\sum_{i=1}^m\vp^{u_j}_z(k\psi^{p_i}_j(s^t_i)).
$$
Let $q\in\{1,\dots,m\}$ such that $p_q\leq j<p_{q+1}$. Then 
$\psi^{p_i}_j\psi^\ell_{p_i}=0$ for each $i=1,\dots,m$ and $\ell\in (p_i,p_{i+1}]$, unless $i=q$ and $\ell>j$, in which case 
$\psi^{p_i}_j\psi^\ell_{p_i}=\psi^\ell_j$. It follows that 
$\psi^{p_i}_j(s^t_i)=0$ for each $i=1,\dots,m$, unless $i=q$. Also,
$$
\psi^{p_q}_j(s^t_q)=\sum_{j<\ell\leq p_{q+1}}t^{\ell-p_q-1}\psi^\ell_j(s^\ell)=t^{j-p_q}s_j+\sum_{j+1<\ell\leq p_{q+1}}t^{\ell-p_q-1}\psi^\ell_j(s^\ell).
$$
Since $\vp^{u_j}_z(ks_j)\neq 0$, it follows that $V^{\mathfrak U^t}_z$ has dimension 1 for general $t$. Furthermore, 
$V^{\h}_z=\lim_{t\to 0}V^{\mathfrak U^t}_z$. Since $\g$ is pure and finitely generated, $\mathfrak U^t\in\lp(\g)$ for general $t$ and $\h=\lim_{t\to 0}\mathfrak U^t$, as wished.
\end{proof}

\begin{proposition}\label{cap} Let $\g$ be a pure nontrivial exact finitely generated linked net of vector spaces over a $\mathbb{Z}^n$-quiver $Q$. Let 
$H$ be a finite collection of vertices $1$-generating $\g$. Let $v_1,\dots,v_m$ be distinct vertices of $Q$. Then the intersection $\lp(\g)_{v_1} \cap\cdots\cap\lp(\g)_{v_m}$ is nonempty only if $\{v_1,\dots,v_m\}$ is a polygon contained in $H$. 
\end{proposition}

\begin{proof} By Proposition~\ref{lemma: non_exact_point_is_in_frontier}, the 
$\h\in\lp(\g)$ minimally $1$-generated by $\{v_1,\dots,v_m\}$ form a dense subset of 
$\lp(\g)_{v_1} \cap\cdots\cap\lp(\g)_{v_m}$. If the intersection is nonempty, so is the subset, and hence $\{v_1,\dots,v_m\}$ is a polygon by Theorem~\ref{12n}. As the $\h$ are also $1$-generated by $H$ we must have 
$\{v_1,\dots,v_m\}\subseteq H$ by Lemma~\ref{minsup1}.
\end{proof}

\begin{theorem} \label{Main_Theorem}
  Let $\g$ be a pure nontrivial exact finitely generated linked net of vector spaces over a $\mathbb{Z}^n$-quiver. Then $\lp(\g)$ is generically nonsingular of pure dimension 
  $\dim(\g)-1$, its irreducible components are rational and equal to the nonempty $\lp(\g)_v$, and 
  the set of exact points on $\lp(\g)$ is its nonsingular locus.
\end{theorem}

\begin{proof}  It follows from \cite[Thm.~7.8]{Esteves_et_al_2021}
% Theorem labeled as ``\label{simple1}" in article I
that the exact points on $\lp(\g)$ lie on the union $\bigcup \lp(\g)_{v}^{*}$, which is
contained in the nonsingular locus of $\lp(\g)$ by Proposition~\ref{birLPv}. Furthermore, the nonexact points are minimally $1$-generated by at least two vertices, by \cite{Esteves_et_al_2021}, Prop.~7.6, and thus lie on the intersection of at least two of the $\lp(\g)_v$ by Proposition~\ref{lemma: non_exact_point_is_in_frontier}. Then
\begin{equation}\label{cuplpg}
\lp(\g)=\bigcup_{v\in H}\lp(\g)_{v},
\end{equation}
where $H$ is a finite set $1$-generating $\g$, and the nonexact points are singular points on $\lp(\g)$, in particular, not on $\bigcup \lp(\g)_{v}^{*}$. It follows that the nonsingular locus of $\lp(\g)$ is $\bigcup \lp(\g)_{v}^{*}$, which is also the set of exact points. The remaining statements follow from \eqref{cuplpg} and Proposition~\ref{birLPv}.
\end{proof}

\section{The shadow partition}\label{shadow}

\begin{Def}
Let $H$ be a nonempty set of vertices of a $\mathbb Z^n$-quiver $Q$ such that $P(H)=H$. 
%By Proposition~\ref{PH123}, for each vertex $v$ of $Q$ there is a unique vertex $w_v\in H$ such that for each $z\in H$ there is an admissible path connecting $z$ to $v$ through $w_v$. We call $w_v$ the \emph{shadow} of $v$ in $H$. Clearly, $w_v=v$ if $v\in H$.
For each $w\in H$, let $R_w$ be the set of vertices $v$ of $Q$ having shadow $w$ in $H$. We call it the \emph{shadow region} of $w$. The collection of shadow regions is called the \emph{shadow partition} associated to $H$.
\end{Def}

The following proposition justifies the definition.

\begin{proposition} Let $H$ be a non-empty set of vertices of a $\mathbb Z^n$-quiver $Q$ such that $P(H)=H$. 
Then the shadow regions $R_w$ for $w\in H$ form a nontrivial partition of the vertex set of $Q$. Furthermore, $R_w\cap H=\{w\}$ for each $w\in H$
\end{proposition}

\begin{proof} The first statement is simply a rephrasing of a consequence of 
\cite{Esteves_et_al_2021}, Prop.~5.7, the fact that each vertex has a unique shadow in $H$. As for the second statement, it follows from the fact that the shadow of $v$ in $H$ is $v$ for each $v\in H$. 
\end{proof}

\begin{figure}[ht]
  \scalebox{0.85}{\begin{minipage}{\textwidth}
  \centering
  \begin{tikzpicture}[commutative diagrams/every diagram]
  \node (p0) at (0:0cm)         {$\textcolor{red}{v_0}$};
  \node (p1) at (90:1.5cm)      {$v_8$};
  \node (p2) at (90+60:1.5cm)   {$v_9$};
  \node (p12) at (90+60:2.2cm)   {$R_0$};
  \node (p3) at (90+2*60:1.5cm) {$v_{10}$};
  \node (p4) at (90+3*60:1.5cm) {$\textcolor{red}{v_2}$};
  \node (p5) at (90+4*60:1.5cm) {$\textcolor{red}{v_1}$};
  \node (p6) at (90+5*60:1.5cm) {$v_7$};
  \node (p7) at (90+3*60:3cm)   {$v_3$};
  \node (p14) at (90+3*60:3.7cm)   {$R_2$};
  \node (p8) at (0:2.6cm)           {$v_6$};
  \node (p9) at (90+4*60:3cm)       {$v_5$};
  \node (p13) at (4:3.3cm)       {$R_1$};
  \node (p10) at (120+3*60:2.6cm)   {$v_4$};
  \node (p11) at (120+2*60:2.6cm)   {$v_{11}$};
  \draw [fill=green!50, draw=none,fill opacity=0.4] (0:0cm) -- (90+2*60:3cm) -- 
  (90+60:3cm) -- (90:3cm) -- (0:0cm) -- cycle;
  \draw [fill=green!50, draw=none,fill opacity=0.4] (90+3*60:1.5cm) -- (110+2*60:3.8cm) -- (90+3*60:4.5cm) -- (71+4*60:4cm) --  (90+3*60:1.5cm) -- cycle;
  \draw [fill=green!50, draw=none,fill opacity=0.4] (90+4*60:1.5cm) -- (90+4*60:4.5cm) -- (12:4cm) -- (64:2.9cm) --  (90+4*60:1.5cm) -- cycle;
\path[commutative diagrams/.cd, every arrow, every label]  
%   (p0) edge node[--,blue] {} (p13)    (p9) edge node[] {} (p8)
  (p8) edge node[] {} (p5)    (p9) edge node[] {} (p8)
  (p10) edge node[] {} (p5)   (p9) edge node[] {} (p10)
  (p4) edge node[] {} (p10)   (p10) edge node[] {} (p7)
  (p0) edge node[] {} (p1)    (p2) edge node[] {} (p0)
  (p0) edge node[] {} (p3)    (p4) edge[red] node[] {} (p0) 
  (p0) edge[red] node[] {} (p5)    (p6) edge node[] {} (p0) 
  (p6) edge node[] {} (p8)    (p5) edge node[] {} (p9)
  (p7) edge node[] {} (p4)    (p11) edge node[] {} (p3) 
  (p4) edge node[] {} (p11)   (p11) edge node[] {} (p7) 
  (p1) edge node[] {} (p2)    (p1) edge node[] {} (p6) 
  (p5) edge node[] {} (p6)    (p5) edge[red] node[] {} (p4) 
  (p3) edge node[] {} (p4)    (p3) edge node[] {} (p2) 
  ;
\end{tikzpicture}
\end{minipage}}
\caption{The shadow partition of a triangle.}
\label{fig3}
\end{figure}
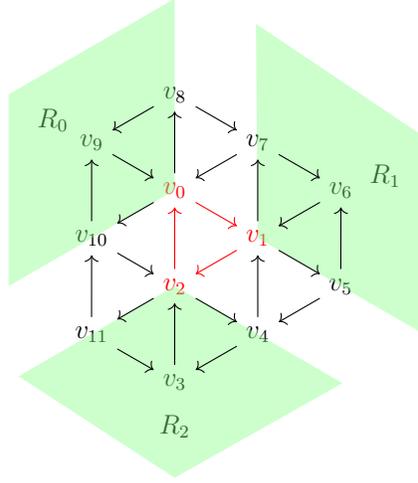

See Figure~\ref{fig3} for the case where $n=2$ and $H$ is a triangle.

Recall that a sequence of vertices 
$w_1,\dots,w_m$ is said to \emph{form an oriented polygon} if there is a sequence $I_1,\dots,I_m$ of pairwise disjoint collections of arrow types such that $I_1\cup\cdots\cup I_m$ is the complete arrow type set and $w_{i+1}=I_i\cdot w_i$ for each $i=1,\dots,m-1$. Actually, it was required that the $I_i$ be nonempty. We drop this requirement here, and say that  $w_1,\dots,w_m$ \emph{form an irredundant oriented polygon} if all the $I_i$ are non-empty.

\begin{Lem}\label{shadowpolygon} Let $H$ be a non-empty set of vertices of a $\mathbb Z^n$-quiver $Q$ such that $P(H)=H$. Let $v_0,\dots,v_m$ be vertices of $Q$ forming an oriented polygon. Let $w_0,\dots,w_m$ be the sequence of their respective shadows in $H$. Then $w_0,\dots,w_m$ form an oriented polygon.
\end{Lem}

\begin{proof} We may assume that $v_0,\dots,v_m$ form an irredundant oriented polygon. By adding the intermediate vertices in picked admissible paths between the vertices $v_i$, we may in addition suppose that $m=n$. 

For each $i=0,\dots,n-1$, let $a_i$ be the arrow connecting $v_i$ to $v_{i+1}$, and $a_n$ that connecting $v_n$ to $v_0$. For each $i=0,\dots,n$, let $\gamma_i$ be an admissible path connecting $w_i$ to $v_i$. For each 
$i=0,\dots,n-1$, let $\rho_i$ be an admissible path connecting $w_i$ to $w_{i+1}$, and $\rho_n$ one connecting $w_n$ to $w_0$. 

Observe that $\gamma_{i+1}\rho_i$ is admissible for each $i=0,\dots,n$, by the defining property of the shadow $w_{i+1}$ of $v_{i+1}$, where we put $\gamma_{n+1}:=\gamma_0$, $w_{n+1}:=w_0$ and 
$v_{n+1}:=v_0$. 

If all the $\rho_i$ are trivial, then all the $w_i$ coincide, and then clearly $w_0,\dots,w_n$ form an oriented polygon. We may now suppose one of the $\rho_i$ is nontrivial. Up to shifting, we may suppose $\rho_0$ is nontrivial. We will show now that $\rho_n\cdots\rho_1\rho_0$ is a minimal circuit, which will end the proof.

Let $\mathfrak a_0$ be the arrow type of $a_0$. Of course, $\rho_n\cdots\rho_1\rho_0$ is a nontrivial circuit. It will thus be enough to show that $\rho_0$ contains at most one arrow of type $\mathfrak a_0$ and that $\rho_i$ does not contain any for any $i>0$. 

Indeed, if $\gamma_0$ contained an arrow of type $\mathfrak a_0$, then $w_0$ would be the shadow of $v_1$ in $H$, contradicting $w_1\neq w_0$. Similarly if $\gamma_1$ contained an arrow of type $\mathfrak a_0$, then there would be a path $\nu$ connecting 
$w_1$ to $v_0$ such that $a_0\nu$ has the same type as $\gamma_1$, and hence $w_1$ would be the shadow of $v_0$ in $H$, contradicting $w_1\neq w_0$. Then $a_0\gamma_0$ contains at most one arrow of type $\mathfrak a_0$, and thus so does $\rho_0$, as $a_0\gamma_0$ and $\gamma_1\rho_0$ connect the same vertices but the latter is admissible. 

Suppose $\gamma_i$ contains no arrow of type $\mathfrak a_0$ for a certain $i\in\{1,\dots,n\}$. Then $a_i\gamma_i$ is admissible, and thus of the same type as $\gamma_{i+1}\rho_i$. It follows that neither 
$\rho_i$ nor $\gamma_{i+1}$ contains an arrow of type $\mathfrak a_0$. By induction, $\gamma_i$ contains no arrow of type $\mathfrak a_0$ for any $i$, and neither does $\rho_i$ for $i=1,\dots,n$.
\end{proof}

\begin{Lem}\label{shadowpath} Let $H$ be a non-empty set of vertices of a $\mathbb Z^n$-quiver $Q$ such that $P(H)=H$. Let $\mu:=\alpha_m\cdots\alpha_1$ for paths $\alpha_i$ of $Q$. For each $i=1,\dots,m$ let $v_i$ be the initial vertex of $\alpha_i$ and $v_{m+1}$ be the final vertex of $\alpha_m$. For each $i=1,\dots,m+1$, let $w_i$ be the shadow of $v_i$ in $H$ and $\gamma_i$ an admissible path connecting $w_i$ to $v_i$. For each $i=1,\dots,m$ let $\rho_i$ be an admissible path connecting $w_i$ to $w_{i+1}$. Put $\rho:=\rho_m\cdots\rho_1$. Then the length of $\mu\gamma_1$ is at least that of $\gamma_{m+1}\rho$, with equality if and only if 
$\alpha_i\gamma_i$ is admissible for each $i$. In particular, if $\mu\gamma_1$ is admissible then equality holds and $\rho$ is admissible. Conversely, if $\rho$ is admissible, and equality holds, then $\mu\gamma_1$ is admissible.
\end{Lem}

\begin{proof} By the defining property of the shadow, the concatenation $\gamma_{i+1}\rho_i$ is admissible for each $i=1,\dots,m$. Since it connects the same vertices as $\alpha_i\gamma_i$, the length of the latter is at least that of the former, being equal if and only if $\alpha_i\gamma_i$ is admissible. Then the lengths of the following concatenations form a increasing sequence, 
$$
\gamma_{m+1}\rho,\quad 
\alpha_m\gamma_m\rho_{m-1}\cdots\rho_1,\quad 
\alpha_m\alpha_{m-1}\gamma_{m-1}\rho_{m-2}\cdots\rho_1,\quad
\mu\gamma_1,
$$
which is constant, or equivalently, the length of $\mu\gamma_1$ is equal to that of $\gamma_{m+1}\rho$, if and only if $\alpha_i\gamma_i$ is admissible for each $i$. The first statement is proved.

If $\mu\gamma_1$ is admissible, since $\gamma_{m+1}\rho$ does not have bigger length, $\gamma_{m+1}\rho$ is admissible as well, and hence has the same length as $\mu\gamma_1$ and $\rho$ is admissible. On the other hand, if $\rho$ is admissible then so is $\gamma_{m+1}\rho$, by the defining property of the shadow. If in addition $\mu\gamma_1$ and $\gamma_{m+1}\rho$ have equal lengths, then $\mu\gamma_1$ is admissible.
\end{proof}

Let $\g$ be a weakly linked net of objects in a $k$-linear Abelian category $\mathcal A$ over a $\mathbb Z^n$-quiver $Q$. Let $H$ be a nonempty set of vertices of $Q$ such that $P(H)=H$. We define a representation $\g_H$ of $Q$ in $\mathcal A$ associated to $\g$ and $H$ as follows. First, for each vertex $v\in V$, set $V^{\g_H}_v:=V^\g_{w_v}$, where $w_v$ is the shadow of $v$ in $H$. Second, given an arrow $a$ of $Q$, let $v_1$ and $v_2$ be its initial and final vertices, and $w_1$ and $w_2$ their respective shadows in $H$. If there is no admissible path from $w_1$ to $v_2$ through $v_1$, put $\vp^{\g_H}_a:=0$; otherwise, let $\vp^{\g_H}_a$ be any map with class $[\vp^\g_\rho]$, where $\rho$ is an admissible path from $w_1$ to $w_2$. Notice that, at any rate, there is an admissible map from $w_1$ to $v_2$ through $w_2$, by the defining property of a shadow.

%\begin{proposition}\label{VHV} Let $\g$ be a weakly linked net over a $\mathbb Z^n$-quiver $Q$ and $H$ a nonempty set of vertices $1$-generating $\g$ satisfying $P(H)=H$. Then $[\vp^\g_a]\vp^{w_1}_{v_1}=\vp^{w_2}_{v_2}[\vp^{\g_H}_a]$ for each arrow $a$ of $Q$, where $v_1$ and $v_2$ are its initial and final vertices, and $w_i$ is the shadow of $v_i$ in $H$ for $i=1,2$. If in addition $\g$ is pure, so is $\g_H$ and the $\vp^{w_i}_{v_i}$ are isomorphisms.
%\end{proposition}

%\begin{proof} If $\g$ is $1$-generated by $H$, for each vertex $v$ of $Q$ the class $\vp^{w_v}_v$ is an epimorphism, where $w_v$ is the shadow of $v$ in $H$. If $a$ is an arrow of $Q$, let $v_1$ and $v_2$ be its initial and final vertices, and $w_1$ and $w_2$ their respective shadows in $H$. If there is no admissible path from $w_1$ to $v_2$ through $v_1$,  then $[\vp^\g_a]\vp^{w_1}_{v_1}$ is zero, and thus coincides with $\vp^{w_2}_{v_2}[\vp^{\g_H}_a]$. Otherwise, $[\vp^{\g_H}_a]=\vp^{w_1}_{w_2}$ by definition and thus $[\vp^\g_a]\vp^{w_1}_{v_1}=\vp^{w_1}_{v_2}=\vp^{w_2}_{v_2}[\vp^{\g_H}_a]$. The last statement is clear.
%\end{proof}

\begin{Lem}\label{shadowphi} Let $\g$ be a weakly linked net over a $\mathbb Z^n$-quiver $Q$ and $H$ a nonempty set of vertices of $Q$ such that $P(H)=H$. Let $\mu$ be a path in $Q$. Let $u$ and $v$ be its initial and final vertices, and $w$ and $z$ their respective shadows in $H$. Let $\gamma$ 
(resp.~$\epsilon$) 
be an admissible path connecting $w$ to $u$ (resp.~$w$ to $z$). 
Then:
\begin{enumerate}
    \item If $\vp^{\g_H}_\mu\neq 0$ then $\mu\gamma$ is admissible.
    \item If $\mu\gamma$ is admissible then 
    $[\vp^{\g_H}_\mu]=[\vp^{\g}_{\epsilon}]$.
\end{enumerate}
\end{Lem}

\begin{proof} Write $\mu=\alpha_m\cdots \alpha_1$ for arrows $\alpha_i$ of $Q$, and keep the notation as in Lemma~\ref{shadowpath}. If $\vp^{\g_H}_\mu\neq 0$ then $\vp^{\g_H}_{\alpha_i}\neq 0$ for each $i$ and thus, by definition of $\g_H$, the concatenation $\alpha_i\gamma_i$ is admissible for each $i$. It follows from Lemma~\ref{shadowpath} that $\mu\gamma$ has the same length as $\gamma_{m+1}\rho$. Since the latter is admissible, so is $\mu\gamma$. 

On the other hand, if $\mu\gamma$ is admissible, then $\alpha_i\gamma_i$ is admissible for each $i$ and $\rho_m\cdots\rho_1$ is admissible, again by Lemma~\ref{shadowpath}. Then, by definition of $\g_H$, the map 
$\vp^{\g_H}_\mu$ has the same class as 
$\vp^\g_{\rho_m}\cdots\vp^\g_{\rho_1}$, whose class is equal to that of $\vp^\g_{\epsilon}$ because $\epsilon$ and $\rho_m\cdots\rho_1$ are admissible and connect the same two vertices.
\end{proof}

\begin{Def} Two distinct vertices $v_1$ and $v_2$ of a $\mathbb Z^n$-quiver $Q$ are said to be \emph{weakly neighbors} if there are a vertex $v$ of $Q$ and simple admissible paths $\gamma_1$ and $\gamma_2$ connecting 
$v$ to $v_1$ and $v_2$, respectively, with no arrow type in common. We call $v$ a \emph{bridge} of $v_1$ and $v_2$. 
\end{Def}

\begin{proposition}\label{bridges} Two neighbors are weakly neighbors, and their bridges are themselves. In particular, the set of bridges of all pairs of vertices in a polygon is the polygon itself.
\end{proposition}

\begin{proof} Let $v_1$ and $v_2$ be neighbors. Then there are simple admissible paths $\mu$ and $\nu$ connecting $v_1$ to $v_2$ and $v_2$ to $v_1$, respectively, whose essential types $I_\mu$ and $I_\nu$ form a partition of the set of arrow types of the quiver. Clearly, $v_1$ and $v_2$ are bridges of $v_1$ and $v_2$. If there were a bridge $v$ distinct from $v_1$ and $v_2$, then a simple admissible path $\gamma_1$ (resp.~$\gamma_2$) connecting $v$ to $v_1$ (resp.~$v_2$) would have essential type $I_{\gamma_1}$ (resp.~$I_{\gamma_2}$) containing $I_\nu$ (resp.~$I_\mu$), because $I_{\gamma_1}\cap I_{\gamma_2}=\emptyset$ and thus neither $\mu\gamma_1$ nor $\nu\gamma_2$ could be admissible. But since $I_{\gamma_1}\cap I_{\gamma_2}=\emptyset$, we would  have that $I_{\gamma_1}=I_\nu$ and $I_{\gamma_2}=I_\mu$, and hence $v=v_2$ and $v=v_1$, a contradiction. 
\end{proof}

%\begin{proof} {\color{red} Mal escrito. Prova que o bridge de dois strictly weakly neighbors é único, desnecessário?} The first statement is clear. As for the second statement, let $v_1,v_2$ be weakly neighbors and $v,v'$ their bridges. Let $\gamma_1$ and $\gamma_2$ be simple admissible paths connecting $v$ to $v_1$ and $v_2$, and $\mu_1$ and $\mu_2$ reverse paths, respectively. Since $\gamma_1$ and $\gamma_2$ dot not have an arrow type in common, we may choose the $\mu_i$ such that $\mu_1=\mu\mu'_1$ and $\mu_2=\mu\mu'_2$, where $\mu'_1$ an $\mu'_2$ have the same type as $\gamma_2$ and $\gamma_1$, respectively. Let $\gamma'_1$ and $\gamma_'2$ be simple admissible paths connecting $v'$ to $v_1$ and $v_2$, respectively. Then $\mu_1\gamma'_1$ and $\mu_2\gamma'_2$ connect $v'$ to $v$, and hence $\mu'_1\gamma'_1$ and $\mu'_2\gamma'_2$ connect the same two vertices. Since $\gamma'_1$ and $\gamma'_2$ dot not have an arrow type in common, it follows that the essential type of $\gamma'_1$ is contained in that of $\mu'_2$, and hence in that of $\gamma_1$. The reverse inclusion is also true, by switching $v$ with $v'$, and since they are simple paths, $\gamma_1$ and $\gamma'_1$ have the same type. Since they have the same end vertex, $v_1$, they must have also the same initial vertex, thus $v=v'$.
%\end{proof}

\begin{proposition}\label{gH} Let $\g$ be a weakly linked net over a $\mathbb Z^n$-quiver $Q$ and $H$ a non-empty set of vertices of $Q$ such that $P(H)=H$. Then:
\begin{enumerate}
    \item $\g_H$ is a weakly linked net over $Q$ generated by $H$.
    \item If $\g$ is pure (resp.~exact), so is $\g_H$.
    \item If $\g$ is a linked net, and every bridge between weakly neighbors of $H$ is in $H$, then $\g_H$ is a linked net as well.
\end{enumerate}
\end{proposition}

\begin{proof} First, let $v_1$ and $v_2$ be vertices of $Q$ and 
$w_1$ and $w_2$ their respective shadows in $H$. Let $\gamma$ be an admissible path connecting $w_1$ to $v_1$. Given two paths $\mu_1,\mu_2$ connecting $v_1$ to $v_2$ with $\mu_2$ admissible, if $\vp^{\g_H}_{\mu_1}$ is nonzero then it follows from Lemma~\ref{shadowphi} that $\mu_1\gamma$ is admissible and 
$[\vp^{\g_H}_{\mu_1}]=\vp^{w_1}_{w_2}$. In particular, also 
$\mu_1$ is admissible and thus has the same type as $\mu_2$. Then $\mu_2\gamma$ is admissible and thus 
$[\vp^{\g_H}_{\mu_2}]=\vp^{w_1}_{w_2}$, by the same lemma. Thus $\vp^{\g_H}_{\mu_1}$ is a scalar 
multiple of $\vp^{\g_H}_{\mu_2}$.

Second, if $\mu$ is a minimal circuit, then $\mu$ is not admissible, and thus $\vp^{\g_H}_{\mu}=0$ by Lemma~\ref{shadowphi}. It follows that $\g_H$ is a weakly linked net.

Third, for each vertex $v$ of $Q$, let $w$ be its shadow in $H$ and let $\gamma$ be an admissible path connecting $w$ to $v$. Then each vertex on $\gamma$ has $w$ as shadow in $H$. It follows directly from the definition that $\vp^{\g_H}_{\gamma}$ is the identity map, whence an isomorphism. So $\g_H$ is $1$-generated by $H$. Statement~(1) is proved.

As for Statement~(2), if $\g$ is pure, so is $\g_H$ because the objects associated to $\g_H$ are among those associated to $\g$. Suppose $\g$ is exact. Let $v_1$ and $v_2$ be neighboring vertices of $Q$. Let $\mu_1$ be an admissible simple path connecting $v_1$ to $v_2$ and $\mu_2$ a reverse path. For each $i=1,2$, let $w_i$ be the shadow of $v_i$ and $\gamma_i$ an admissible path connecting $w_i$ to $v_i$. Since $\g_H$ is a weakly linked net, 
\begin{equation}\label{KI}
\text{Ker}(\vp^{\g_H}_{\mu_2})\supseteq\text{Im}(\vp^{\g_H}_{\mu_1}).
\end{equation}
We need to show equality.  

If both $\mu_1\gamma_1$ and $\mu_2\gamma_2$ are admissible, then $\vp^{\g_H}_{\mu_1}$ has class 
$\vp^{w_1}_{w_2}$ whereas 
$\vp^{\g_H}_{\mu_2}$ has class 
$\vp^{w_2}_{w_1}$. Now, $w_1,w_2$ form an oriented polygon by 
Lemma~\ref{shadowpolygon}. Also, $w_1\neq w_2$ by \cite{Esteves_et_al_2021}, Lem.~5.2, since $v_1$ and $v_2$ are neighbors and $\mu_1\gamma_1$ and $\mu_2\gamma_2$ are admissible. Since $\g$ is exact, $\text{Im}(\vp^{w_1}_{w_2})=\text{Ker}(\vp^{w_2}_{w_1})$, from which \eqref{KI} follows.

If $\mu_2\gamma_2$ is not admissible, then 
the essential type of $\mu_1$ is contained in that of $\gamma_2$. Since $\mu_1$ is simple, there is an admissible path $\gamma'_2$ connecting $w_2$ to $v_1$ such that $\mu_1\gamma'_2$ has the same type as $\gamma_2$. But then $w_2$ is the shadow of $v_1$ in $H$ and thus $w_1=w_2$. The same conclusion is reached if $\mu_1\gamma_1$ is not admissible, the remaining case to argue. So we may assume $w_1=w_2$.

Equality in \eqref{KI} follows now because either 
$\vp^{\g_H}_{\mu_1}$ or $\vp^{\g_H}_{\mu_2}$ is the identity map. Indeed, since $w_1=w_2$, and since $v_1$ and $v_2$ are neighbors, \cite{Esteves_et_al_2021}, Lem.~5.2, implies that either 
$\mu_1\gamma_1$ is admissible or $\mu_2\gamma_2$ is admissible.

If $\mu_1\gamma_1$ is admissible, it has the same type as $\gamma_2$, which connects $w_2$, the shadow of $v_2$ in $H$, to $v_2$. It follows that all vertices on $\mu_1\gamma_1$ have shadow $w_2$ in $H$ as well. But $w_2=w_1$. From the definition of $\g_H$, the map $\vp^{\g_H}_{\mu_1}$ is the identity map. Similarly, if $\mu_2\gamma_2$ is admissible then $\vp^{\g_H}_{\mu_2}$ is the identity map. Statement~(2) is proved. 

Finally, assume $\g$ is a linked net. Let $\mu_1$ and $\mu_2$ be simple admissible paths leaving the same vertex $v$ of $Q$ with no arrow type in common. By \cite{Esteves_et_al_2021}, Lem.~6.6, to prove that $\g_H$ is a linked net, it is enough to show that 
$$
\text{Ker}(\vp^{\g_H}_{\mu_1})\cap
\text{Ker}(\vp^{\g_H}_{\mu_2})=0.
$$

Let $v_1$ and $v_2$ denote the respective final vertices of $\mu_1$ and $\mu_2$. Let $w,w_1,w_2$ denote the respective shadows of $v,v_1,v_2$ in $H$. Let $\gamma$ be an admissible path connecting $w$ to $v$. Since $\mu_1$ and $\mu_2$ have no arrow type in common, $\mu_1\gamma$ or $\mu_2\gamma$ is admissible. 

Suppose first that one of them is not, say $\mu_1\gamma$ is not admissible. From Lemma~\ref{shadowphi}, we get $\vp^{\g_H}_{\mu_1}=0$. But also, the essential type of $\mu_2$ is contained in that of $\gamma$. Then $w$ is the shadow of $v_2$ in $H$, that is $w=w_2$ and hence $\vp^{\g_H}_{\mu_2}$ is the identity, by definition of $\g_H$. So $\text{Ker}(\vp^{\g_H}_{\mu_1})\cap\text{Ker}(\vp^{\g_H}_{\mu_2})=0$. 

Suppose now that both $\mu_1\gamma$ and $\mu_2\gamma$ are admissible. By Lemma~\ref{shadowphi},
$$
\text{Ker}(\vp^{\g_H}_{\mu_1})\cap\text{Ker}(\vp^{\g_H}_{\mu_2})=
\text{Ker}(\vp^{w}_{w_1})\cap\text{Ker}(\vp^{w}_{w_2}).
$$
If $w_1=w$ or $w_2=w$ then 
the above intersection is clearly zero. 

Suppose $w_1\neq w$ and $w_2\neq w$. Let $\rho_i$ be an admissible path connecting $w$ to $w_i$ for $i=1,2$. The paths $\rho_i$ are simple by Lemma~\ref{shadowpolygon}. Since $\mu_1\gamma$ and $\mu_2\gamma$ are admissible, and $\mu_1$ and $\mu_2$ have no arrow types in common, it follows from Lemma~\ref{shadowpath} that the intersection of the essential types of $\rho_1$ and $\rho_2$ is contained in the essential type of $\gamma$. Thus we may choose the $\rho_i$ and $\gamma$ such that 
$\rho_i=\rho'_i\gamma'$ and $\gamma=\gamma''\gamma'$ for a path $\gamma'$ such that $\rho'_1$ and $\rho'_2$ have no arrow type in common. But then the final vertex $z$ of $\gamma'$ is a bridge of $w_1$ and $w_2$. If $H$ is closed under adding bridges then 
$z\in H$. But lies on $\gamma$, which connects $w$, the shadow of $v$ in $H$, to $v$. Thus $z=w$, that is, $\gamma'$ is trivial and hence $\rho_1$ and $\rho_2$ have no arrow type in common. Since $\g$ is a linked net, 
$$
\text{Ker}(\vp^{w}_{w_1})\cap\text{Ker}(\vp^{w}_{w_2})=
\text{Ker}(\vp_{\rho_1})\cap\text{Ker}(\vp_{\rho_2})=0.
$$
Statement~(3) is proved. 
\end{proof}

%Given a subrepresentation $\mathfrak W$ of $\g$, the associated representation $\mathfrak W_H$ is a subrepresentation of $\g_H$ given by the subspaces $V^{\mathfrak W}_{w_v}\subseteq V^{\g_H}_{v}$. Indeed, if $a$ is an arrow of $Q$ from vertex $v_1$ to $v_2$, and $\vp^{\g_H}_a$ is nonzero then $\vp^{\g_H}_a(V^{\mathfrak W}_{w_{v_1}})=\vp^{w_{v_1}}_{w_{v_2}}(V^{\mathfrak W}_{w_{v_1}})\subseteq V^{\mathfrak W}_{w_{v_2}}$.

%Given a representation $\g'$ of the full subquiver of $Q$ with vertices in $H$, there is a representation $\g$ of $Q$ such that $\g_H$ is projectively equivalent to $\g'$.

\section{\texorpdfstring{$\lp(\g)$}{LP(V)} is a local complete intersection}\label{CH_LPg}

%Assume from now on that $n=2$, that is, that $Q$ is a $\Z^2$-quiver. 

%{\color{red} Alterei o lema abaixo e sua prova para $\g$ suportados em um polígono. A estrutura da prova ficou mais simples.}

\begin{Lem}\label{prop_triagle} Let $\g$ be a pure nontrivial exact linked net of vector spaces over a $\Z^n$-quiver $Q$. If $\g$ is generated by a polygon 
then $\lp(\g)$ is a local complete intersection.
\end{Lem}

\begin{proof} Let $v_0,\dots,v_n$ be vertices of $Q$ forming an oriented $(n+1)$-gon $\Delta$ generating $\g$. For convenience, put $v_{n+1}:=v_0$. Let $r:=\dim(\g)$.
%; see \eqref{quiver_Q_1_3_cha_5} for the case $n=2$.
%\begin{equation}\label{quiver_Q_1_3_cha_5}
%  \g:\hspace{.5cm}
%  \begin{tikzpicture}[commutative diagrams/every diagram]
%  \node (102) at (0:0cm)          {$V_{v_2}$.};  
%  \node (021) at (90+60:2.5cm)    {$V_{v_0}$}; 
%  \node (210) at (90+2*60:2.5cm)  {$V_{v_1}$};    
%  \path[commutative diagrams/.cd, every arrow, every label]    
%  (021) edge node[fill=white,anchor=center, pos=0.5] {$ \vp^{v_0}_{v_1} $} (210)
%  (210) edge node[fill=white,anchor=center, pos=0.5] {$ \vp^{v_1}_{v_2} $} (102)
%  (102) edge node[fill=white,anchor=center, pos=0.5] {$ \vp^{v_2}_{v_0} $} (021)
%  ;     
%   \end{tikzpicture}
%\end{equation}
%
%As 
% $$
% \lp(\g)=\bigcup_{i=0}^{2}\lp(\g)_{v_i},
% $$
% and each $\lp(\g)_{v_i}^*$ is reduced, it follows that $\lp(\g)$ is generically reduced. Then, by \cite[Prop. 14.124]{Ulrich_1}, it is enough to show that $\lp(\g)$ is Cohen--Macaulay.

As $\g$ is exact, it follows from \cite{Esteves_et_al_2021}, Prop.~9.1, that $\g$ is the direct sum of locally finite exact linked nets $\h_1,\dots,\h_r$ of vector spaces of dimension 1. Since $\g$ is $1$-generated by $\Delta$, so are the $\h_j$. By \cite{Esteves_et_al_2021}, Prop.~7.6, the net $\h_j$ is generated by $v_{\ell_j}$ for some $\ell_j\in\{0,\dots,n\}$ for each $j=1,\dots,r$; let $s_j\in V^{\h_j}_{v_{\ell_j}}$ be a generator. For each $\ell=0,\dots,n$, let $r_\ell:=\#\{j\,|\,\ell_j=\ell\}$. Clearly, $\sum r_\ell=r$. The $s_j$ induce a basis for $V^\g_{v_i}$ for each $i=0,\dots,n$, and thus a decomposition $V^\g_{v_i}=V_{i,0}\oplus\cdots\oplus V_{i,n}$, 
where $V_{i,\ell}$ is the subspace generated by the $\vp^{v_\ell}_{v_i}(s_j)$
 for all $j$ with $\ell_j=\ell$, for $\ell=0,\dots,n$; each $V_{i,\ell}$ has dimension $r_\ell$.

For each $i=0,\dots,n$, the map $\vp^{v_i}_{v_{i+1}}$ can then be represented by a diagonal matrix $M_i$. For $i=0,\dots,n-1$, all of its entries are 1 but those in positions $r_0+\cdots+r_i+j$ for $j=1,\dots,r_{i+1}$, which are 0. The matrix $M_n$ has all of its entries 1 but those in positions $1,\dots,r_0$.

Let $\mathbb G$ be the product of $n+1$ copies of 
$\mathbb P^{r-1}$. As $\Delta$ is equal to its hull, we have that $\lp(\g)$ is isomorphic to the subscheme $X$ of points $([x_0],\dots,[x_n]) \in \mathbb G$ satisfying the equations
\begin{equation}\label{eqslpg}
M_0x_0\wedge x_1=0,\quad
\dots,\quad
M_{n-1}x_{n-1}\wedge x_n=0,\quad
 M_nx_n\wedge x_0=0.
\end{equation}
We need only prove $X$ is a local complete intersection.
%Let $x_j$, $y_j$, and $z_j$ be the homogeneous coordinates of %$\prod_{i=0}^{2}\p(V_{v_i})$, with $j=0,\dots,r$.

Since $\mathbb G$ is smooth with
$$
\dim\mathbb G=(n+1)(r-1)=(nr-n)+(r-1),
$$
and since $X$ has pure dimension $r-1$, because so has $\lp(\g)$ by 
Theorem~\ref{Main_Theorem}, we need only prove that $X$ is locally given by $nr-n$ equations.

%\begin{proposition}\label{prop_triagle}
%  The linked projective space of an exact linked net over a $\Z^2$-quiver with support on a triangle quiver is Cohen-Macaulay, reduced, and connected.
%\end{proposition}
%\begin{proof}
 
%Let the linked net of vector spaces $\g$ be the one in \eqref{quiver_Q_1_3_cha_5}. 
Write $x_i=(x_{i,0},\dots,x_{i,n})$ for each $i=0,\dots,n$, and 
$x_{i,\ell}=(x^1_{i,\ell},\dots,x^{r_\ell}_{i,\ell})$ for each $\ell=0,\dots,n$. For convenience, put $x^j_{n+1,\ell}:=x^j_{0,\ell}$ and $x_{n+1,\ell}:=x_{0,\ell}$ for each $\ell$ and $j$ and set $x_{n+1}:=x_0$.
% \begin{eqnarray}
 %  (s_{00}+s_{02})\wedge(s_{10}+s_{12}) & = & 0\label{f01_2}\\
  % (s_{00}+s_{02})\wedge s_{11} & = & 0\label{f01_1}\\
   %(s_{10}+s_{11})\wedge(s_{20}+s_{21}) & = & 0\label{f12_2}\\
   %(s_{10}+s_{11})\wedge s_{22} & = & 0\label{f12_1}\\
   %(s_{21}+s_{22})\wedge(s_{01}+s_{02}) & = & 0\label{f20_2}\\
   %(s_{21}+s_{22})\wedge s_{00} & = & 0         \label{f20_1}.
% \end{eqnarray}
%As Cohen--Macaulayness is a local property, we prove that $\lp(\g)$ is Cohen--Macaulay 
%using an open cover  
For each $i,\ell\in \{0,\dots,n\}$ and $j=1,\dots,r_\ell$, let $D^j_{i,\ell}$ be 
the open subset of $\mathbb G$ where $x^j_{i,\ell}\neq 0$. Put 
$D_{i,\ell}:=\bigcup D^j_{i,\ell}$ for each $i,\ell$. 
%For $\h=([x_0],[x_1],[x_2])\in D_{0,0}$, Equation~(10) is equivalent to $(x_{2,1},x_{2,2})=0$, which implies $x_{2,0}\neq 0$ and Equations~(8)~and~(9). Since $x_{2,0}\neq 0$ and $x_{2,1}=0$, Equation~(7) is equivalent to $x_{1,1}=0$, which implies Equation~(6), and $x_{1,0}\wedge x_{2,0}=0$. Since $x_{0,0}$ and $x_{2,0}$ are nonzero, we have
%$$
%r_1+r_2+r_1+r_0-1+r_0+r_2-1=2(r_0+r_1+r_2-1)=2r
%$$
%equations. By symmetry, we conclude that $\lp(\g)\cap D_{i,i}$ is Cohen--Macaulay for $i=0,1,2$. 
We claim that 
$$
X\subseteq D_{0,1}\cup D_{1,2}\cup \cdots\cup D_{n-1,n}\cup D_{n,0}.
$$
Indeed, were $([x_0],\dots,[x_n])\in X$
such that $x_{i,i+1}=0$ for $i=0,\dots,n$, then $M_ix_i$ would be nonzero, and thus a nonzero scalar multiple of $x_{i+1}$ for each $i=0,\dots,n$. But then $M_n\cdots M_1M_0x_0$ would be nonzero, an absurd.

By symmetry, we need only prove $X$ is locally given by $nr-n$ equations on $D^{j_0}_{0,1}$ for each $j_0$. Now, the equation $M_nx_n\wedge x_0=0$ on $D^{j_0}_{0,1}$ is equivalent to
$x_{0,1}^{j_0}x_{n,\ell}=x_{n,1}^{j_0}x_{0,\ell}$ for $\ell=1,\dots,n$, a total of $r_1+\cdots+r_n-1$ equations, and $x_{n,1}^{j_0}x_{0,0}=0$. Furthermore, they imply that $x_{n,0}\neq 0$ or $x_{n,1}\neq 0$, so we need only prove that $X$ is locally given by $nr-n$ equations on 
$D_{0,1}^{j_0}\cap D_{n,p_n}^{j_n}$ for each $p_n\in\{0,1\}$ and each integers $j_0,j_n$. 

Suppose by descending induction on $i$ that 
we need only prove $X$ is locally given by $nr-n$ equations on 
$$
D_{0,1}^{j_0}\cap D_{n,p_n}^{j_n}\cap\cdots\cap 
D_{i+1,p_{i+1}}^{j_{i+1}}
$$
for an integer $i\in\{1,\dots,n-1\}$, each $p_n\in\{0,1\}$ and $p_s\in\{s+1,p_{s+1}\}$ for $s=i+1,\dots,n-1$, and all integers $j_0,j_{i+1},\dots,j_n$. Notice that $p_s\in\{0,1,s+1,\cdots,n\}$, hence $p_s\neq s$ for each $s$. Now, the equation $M_ix_i\wedge x_{i+1}=0$ on that open set is equivalent to 
$x_{i+1,p_{i+1}}^{j_{i+1}}x_{i,\ell}=x_{i,p_{i+1}}^{j_{i+1}}x_{i+1,\ell}$ for $\ell=0,\dots,i,i+2,\dots,n$, a total of $r-r_{i+1}-1$ equations, and $x_{i,p_{i+1}}^{j_{i+1}}x_{i+1,i+1}=0$. They imply that $x_{i,i+1}\neq 0$ or $x_{i,p_{i+1}}\neq 0$, so we need only prove that $X$ is locally given by $nr-n$ equations on 
$$
D_{0,1}^{j_0}\cap D_{n,p_n}^{j_n}\cap\cdots\cap 
D_{i+1,p_{i+1}}^{j_{i+1}}\cap D_{i,p_i}^{j_i}
$$
for each $p_n\in\{0,1\}$ and $p_s\in\{s+1,p_{s+1}\}$ for $s=i,\dots,n-1$, and all integers $j_0,j_i,\dots,j_n$.

By induction, it follows that we need only prove that $X$ is 
given by $nr-n$ equations on
$$
D_{0,1}^{j_0}\cap D_{n,p_n}^{j_n}\cap\cdots\cap 
D_{i,p_i}^{j_i}\cap\cdots\cap D_{1,p_1}^{j_1}
$$ 
for each $p_n\in\{0,1\}$ and $p_s\in\{s+1,p_{s+1}\}$ for $s=1,\dots,n-1$, and all integers $j_0,j_1,\dots,j_n$. Put $p_{n+1}:=p_0:=1$ and  $j_{n+1}:=j_0$ for convenience. Put 
$$
y_i:=\frac{x_{i,p_{i+1}}^{j_{i+1}}}{x_{i,p_i}^{j_i}}\quad\text{for each }i=0,\dots,n.
$$
As we have seen, Equations \eqref{eqslpg} on that open set are equivalent to
\begin{equation}\label{batch1}
\begin{cases}
x_{i,\ell}=y_ix_{i+1,\ell} &\text{for }i=0,\dots,n-1 \text{ and all }\ell\neq i+1,\\
x_{n,\ell}=y_nx_{0,\ell} &\text{for all }\ell\neq 0,
\end{cases}
\end{equation}
a total of $nr-n-1$ equations, and 
\begin{equation}\label{batch2}
y_ix_{i+1,i+1}=0 \text{ for } i=0,\dots,n.
\end{equation}
But Equations \eqref{batch1} imply
$$
y_ix_{i+1,i+1}=y_iy_{i+1}x_{i+2,i+1}=\cdots=y_i\cdots y_nx_{0,i+1}=
y_i\cdots y_ny_0\cdots y_{i-1}x_{i,i+1}
$$
for each $i=0,\dots,n$. Since $x_{0,1}\neq 0$, it follows that Equations \eqref{batch2} are all equivalent to a single equation: $y_0\cdots y_n=0$. 
\end{proof}

%Now we can proof the following Theorem 

\begin{theorem}\label{main_result} Let $\g$ be a pure nontrivial exact finitely generated linked net of vector spaces over a $\mathbb Z^n$-quiver $Q$. Then $\lp(\g)$ is local complete intersection and reduced.
%, and connected.
\end{theorem}

\begin{proof} %That $\lp(\g)$ is connected follows from Proposition \ref{lemma: non_exact_point_is_in_frontier}. In addition, s
Since $\lp(\g)$ is generically nonsingular by Theorem~\ref{Main_Theorem}, if $\lp(\g)$ is a local complete intersection, thus Cohen--Macaulay, then $\lp(\g)$ is reduced by \cite[Prop.~14.126]{Ulrich_1}. It is thus enough to show that $\lp(\g)$ is a local complete intersection.

%Let $H$ be a finite collection of vertices of $Q$ supporting $\g$ such that $P(H)=H$. 

Let $\h\in\lp(\g)$. By Theorem~\ref{12n}, there is a polygon $\Delta$ generating $\h$. %Then there is a minimal circuit $a_n\cdots a_0$ such that $v_i$ is the initial vertex of $a_i$ for $i=0,\dots,n$. Let $\mathfrak a_i$ be the type of $a_i$ and $R_i:=C_{T-\{\mathfrak a_i\}}(v_i)$ for each $i=0,\dots,n$; see Figure~\ref{fig3} for the case where $n=2$.
Since $P(\Delta)=\Delta$, we may consider the associated representation $\g_\Delta$. It follows from Proposition~\ref{gH} that $\g_\Delta$ is a pure exact nontrivial linked net of vector spaces over $Q$ generated by $\Delta$, and hence $\lp(\g_\Delta)$ is a local complete intersection by Lemma~\ref{prop_triagle}.

Let $H$ be a finite set of vertices containing $\Delta$, generating $\g$ and satisfying $P(H)=H$. Given a vertex $v$ of $Q$ denote by $w_v$ its shadow in $\Delta$. The $\mathfrak X\in\lp(\g)$ generated by $\Delta$ form an open subscheme $U$ given by 
$\varphi^{w_v}_v(V^{\mathfrak X}_{w_v})\neq 0$ for each $v\in H$. For each such $\mathfrak X$ there is a corresponding subnet $\mathfrak Y$ of $\g_\Delta$ generated by all $V^{\mathfrak X}_u$ for $u\in\Delta$. Then $\mathfrak Y\in\lp(\g_\Delta)$. Indeed, given a vertex $v$ of $Q$, since $\vp^w_{w_v}(V^{\mathfrak X}_w)\subseteq V^{\mathfrak X}_{w_v}$ for each $w\in\Delta$, it follows that $V^{\mathfrak Y}_v=\vp^{\g_\Delta}_\mu(V^{\mathfrak X}_{w_v})$ for an admissible path connecting $w_v$ to $v$. Since $\vp^{\g_\Delta}_\mu$ is the identity, $\mathfrak Y$ is a pure subnet of $\g$ of dimension 1, whence $\mathfrak Y\in\lp(\g_\Delta)$.

Let
$$
\Theta\colon U \longrightarrow \lp(\g_\Delta)
$$
be the map taking $\mathfrak X\in U$ to $\mathfrak Y$, as above. It is a scheme morphism because its composition with the embedding $\psi^{\g_\Delta}_\Delta$ is the composition of $\psi_H^\g$ with the projection map. Of course, $\mathfrak Y$ determines $\mathfrak X$ for $\mathfrak X\in U$. Also, the image of $\Theta$ is in the open subset $U'$ of $\lp(\g_\Delta)$ given by 
$\varphi^{w_v}_v(V^{\mathfrak Y}_{w_v})\neq 0$ for each $v\in H$. We claim the induced map $\Theta\colon U\to U'$ is an isomorphism. 

Indeed, given $\mathfrak Y\in U'$, we let $\mathfrak X$ be the subnet of $\g$ generated by all $V^{\mathfrak Y}_u$ for $u\in\Delta$. As before, 
for each vertex $v$ of $Q$, we have $V^{\mathfrak X}_v=\vp^{w_v}_v(V^{\mathfrak Y}_{w_v})$, which is of dimension 1 because $\mathfrak Y\in U'$. Thus $\mathfrak X\in\lp(\g)$. The assignment $\mathfrak Y\mapsto\mathfrak X$ is clearly a scheme morphism and the inverse to the morphism $U\to U'$. 

Since $\lp(\g_\Delta)$ is a local complete intersection, so are $U'$ and hence $U$. As 
$U$ is a neighborhood of $\h$, we have that $\lp(\g)$ is a local complete intersection around $\h$. As $\h\in\lp(\g)$ was arbitrary, $\lp(\g)$ is a local complete intersection.
\end{proof}

\section{Smoothings}\label{smoothings}

\begin{Def}
A \emph{general linked net} over a $\mathbb Z^n$-quiver $Q$ of objects in a $k$-linear Abelian category $\mathcal A$ is a representation $\mathfrak V$ of $Q$ in $\mathcal A$ such that 
\begin{enumerate}
    \item $\vp^\g_a$ is an isomorphism for each arrow $a$ of $Q$;
    \item for each two paths $\gamma_1$ and $\gamma_2$ connecting the same two vertices, $\vp^\g_{\gamma_1}$ is a scalar multiple of $\vp^\g_{\gamma_2}$. 
\end{enumerate}
\end{Def}

As before, for each two vertices $u$ and $v$ of $Q$ we may define 
$\vp^u_v:=[\vp^g_\mu]$ for any path $\mu$ connecting $u$ to $v$. 

Given a pure nontrivial general linked net $\mathfrak V$ over $Q$, for each vertex $v$ of $Q$, the natural map $\psi_v\colon\lp(\mathfrak V)\to\mathbb P(V_v)$ is a bijection. 
Indeed, given a one-dimensional subspace $W\subseteq V_v$, put 
$W_w:=\vp_{\gamma}(V_v)$ for each vertex $w$ of $Q$, where $\gamma$ is any path connecting $v$ to $w$. The $W_w$ are one-dimensional by Property~(1). They are well-defined and form a subrepresentation $\h\subseteq\g$ by Property~(2). We have thus a well-defined map $W\mapsto\mathfrak W$, which is the inverse to $\psi_v$. Furthermore, given another vertex $u$, we have 
that $\psi_u=\psi^v_u\psi_v$, where $\psi^v_u$ is the isomorphism given by $\vp^v_u$. Thus $\psi_v$ induces a scheme structure on $\lp(\mathfrak V)$ which is independent of the choice of $v$. In addition, given a finite set of vertices $H$, the natural map
$$
\psi_H=\prod_{v\in H}\psi_v\colon\lp(\mathfrak V)\longrightarrow\prod_{v\in H}
\mathbb P(V_v)
$$
is an isomorphism onto a small diagonal. 

Let $R$ be a discrete valuation ring with residue field $k$ and field of fractions $K$. Let 
$\mathfrak M$ be a representation of $Q$ in the category of free modules of a given rank $n$ over $R$. For convenience, put $M_v:=V^{\mathfrak M}_v$ for each vertex $v$ of $Q$. Assume the induced representation by vector spaces over $K$ is a general linked net and that over $k$ is a weakly linked net $\g$. We call $\mathfrak M$ a \emph{smoothing} of $\g$ over $R$ and say $\g$ is \emph{smoothable}.

Let $H$ be a finite set of vertices of $Q$. Let $B:=\text{Spec}(R)$. Define $\lp_H(\mathfrak M)$ as the $B$-subscheme 
$$
\lp_H(\mathfrak M)\subseteq\prod_{v\in H}\text{Proj}(\text{Symm}(M_v))
$$
of the $B$-product given by the vanishing of the maps of vector bundles
$$
\begin{CD}
\mathcal O_v(-1)\otimes\mathcal O_w(-1) @>(\rho_\mu,1)>> 
\bigwedge^2 \widetilde M_w,
\end{CD}
$$
for all $v,w\in H$, where $\widetilde M_v$ is the pullback of the locally free sheaf associated to $M_v$ on $B$ and $\mathcal O_v(-1)$ is the pullback of the tautological subsheaf on the scheme $\text{Proj}(\text{Symm}(M_v))$ for each vertex $v$ of $Q$, and 
$\rho_\mu\colon\widetilde M_v\to\widetilde M_w$ is the map induced by $\vp^{\mathfrak M}_\mu$ for any path $\mu$ connecting $v$ to $w$. 

\begin{theorem}\label{smoothg} Let $\g$ be a finitely generated exact pure nontrivial linked net over $Q$ of 
vector spaces over $k$. Let $H$ be a finite set of vertices of $Q$ generating $\g$ with $P(H)=H$. Let $\mathfrak M$ be a smoothing of $\g$ over a discrete valuation ring $R$ with residue field $k$. Then $\lp_H(\mathfrak M)$ is reduced and flat over $B:=\text{\rm Spec}(R)$ and $\lp_H(\g)$ is a degeneration of the small diagonal in 
$\prod_{v\in H}\mathbb P(V_v)$.
\end{theorem}

\begin{proof} Let $K$ be the field of fractions of $R$. Since $\g$ is the representation by vector spaces over $k$ induced by $\mathfrak M$, it follows that the special fiber of $\lp_H(\mathfrak M)$ over $B$ is 
$\lp_H(\g)$. Also, the general fiber is isomorphic to each factor $\text{Proj}(\text{Symm}(M_v)\otimes K)$
under the projection, and is thus a small diagonal. 

It remains to show $\lp_H(\mathfrak M)$ is reduced and $B$-flat. Since $P(H)=H$, the special fiber is isomorphic to $\lp(\g)$, and is thus geometrically reduced by Theorem~\ref{main_result}. In addition, no topological component of $\lp_H(\mathfrak M)$ is contained in the special fiber. Indeed, it is enough to show that each general point on $\lp_H(\g)$ is on a section of $\lp_H(\mathfrak M)$ over $B$. Now, a general point on $\lp(\g)$ corresponds to an exact linked subnet $\h\subseteq\g$, that is, to a vertex $w$ of $Q$ and an element $s\in M_w\otimes k$ such that $\vp^\g_\mu(s)$ is nonzero for each admissible path $\mu$ leaving $w$. Lift $s$ to an element $\widetilde s$ of $M_w$. Then $\vp^{\mathfrak M}_\mu(\widetilde s)$ lifts $\vp^{\mathfrak V}_\mu(s)$ for each admissible path $\mu$. 

Given a vertex $v$ of $Q$, and two paths $\mu_1$ and $\mu_2$ connecting $w$ to $v$, since $\mathfrak M$ restricts to a general linked net over $K$, there are $x,y\in R-\{0\}$ with no common factor such that $y\vp^{\mathfrak M}_{\mu_1}=x\vp^{\mathfrak M}_{\mu_2}$. But then $\overline y\vp_{\mu_1}=\overline x\vp_{\mu_2}$, where $\overline x$ and $\overline y$ are the residue classes of $x$ and $y$. Assume $\mu_2$ is admissible. Since $\mathfrak M$ restricts to a weakly linked net over $k$, there is $c\in k$ such that $\vp_{\mu_1}=c\vp_{\mu_2}$, and since $\vp_{\mu_2}(s)\neq 0$, we must have $\overline y c=\overline x$. Since $x$ and $y$ have no common factor, it follows that $y\in R^*$, and hence 
$\vp^{\mathfrak M}_{\mu_1}=(x/y)\vp^{\mathfrak M}_{\mu_2}$. 

Thus, for each vertex $v$ of $Q$, let $\mu$ be any admissible path connecting $w$ to $v$, and consider the $R$-submodule of $M_v$ generated by $\vp^{\mathfrak M}_\mu(\widetilde s)$. It does not depend on the choice of $\mu$. It is of rank 1 with free quotient since $\vp^{\mathfrak V}_\mu(s)\neq 0$, and thus gives rise to a section of $\text{Proj}(\text{Symm}(M_v))$ over $B$. Putting together all these sections for $v\in H$, we have a section of $\prod_{v\in H}\text{Proj}(\text{Symm}(M_v))$ contained in $\lp_H(\mathfrak M)$. 

As in \cite{OssermanFlatness}, Lem.~4.3, p.~3388, we conclude that the reduced induced subscheme associated to $\lp_H(\mathfrak M)$ is flat over $B$. And as in loc.~cit., we conclude that $\lp_H(\mathfrak M)$ is reduced, whence flat, from \cite{Osserman2006}, Lem.~6.13, p.~1191.
%The last statement follows from the invariance of the Hilbert polynomial in flat families.
\end{proof}

\section{Divisors}\label{divisors}

\begin{proposition}\label{LDeltaAb} Let $\mathfrak V$ be an exact weakly linked net of objects in an $k$-linear Abelian category $\mathcal A$ over a $\mathbb Z^n$-quiver $Q$. Let $v_1,\dots,v_m$ be vertices of $Q$ forming an oriented polygon $\Delta$ with $m\geq 2$. Then the class of 
$$
V_\Delta:=\text{Ker}(\vp^{v_1}_{v_2})\oplus\cdots\oplus\text{Ker}(\vp^{v_{m-1}}_{v_m})\oplus \text{Ker}(\vp^{v_m}_{v_1})
$$
in the Grothendieck group of $\mathcal A$ is 
equal to that of $V_v$ for every vertex $v$ of $Q$.
\end{proposition}

\begin{proof} For each $i=1,\dots,m$, let $M_i:=\text{Im}(\vp^{v_1}_{v_{i}})\subseteq V_{v_{i}}$. Then $M_1=V_{v_1}$. Also, the map $\vp^{v_{i}}_{v_{i+1}}$ restricts to a surjection $M_i\to M_{i+1}$ for $i=1,\dots,m-1$.  Furthermore, since $\mathfrak V$ is exact, 
 $\text{Ker}(\vp^{v_{i}}_{v_{i+1}})=\text{Im}(\vp_{v_{i}}^{v_{i+1}})$ for each $i=1,\dots,m$, and since $v_1,\dots,v_m$ form an oriented polygon, $M_m=\text{Ker}(\vp^{v_m}_{v_{1}})$ and $\text{Im}(\vp_{v_{i}}^{v_{i+1}})\subseteq M_i$ for $i=1,\dots,m-1$. We obtain an exact sequence 
\begin{equation}\label{MMsimp}
0 \to \text{Ker}(\vp^{v_i}_{v_{i+1}}) \to
M_{i} \to M_{i+1} \to 0
\end{equation}
for each $i=1,\dots,m-1$. Since $M_1=V_{v_1}$ and $M_m=\text{Ker}(\vp^{v_m}_{v_{1}})$, it follows that the class of $V_\Delta$ in the Grothendieck group of $\mathcal A$ is that of $V_{v_1}$, and hence that of $V_v$ for every vertex $v$ of $Q$ by \cite{Esteves_et_al_2021}, Prop.~9.3.
\end{proof}

%If a linked net $\g$ of vector spaces over a $\mathbb Z^n$-quiver $Q$ arises from a degeneration of linear series, then there is a representation $\mathfrak L$ of $Q$ in the category of line bundles such that $\g$ is a subrepresentation of the representation $H^0(\mathfrak L)$, induced by taking global sections; see \cite{Esteves_et_al_2021}. We will take this below as a definition of a ``limit linear series."

%\begin{proposition} Let $F$ be a simple torsion-free, rank-one sheaf on $X$ and let $U\subseteq\mathbb P(H^0(X,F))$ be the open subscheme parameterizing classes of sections $s$ of $F$ that do not vanish generically anywhere. Then the map $U\to\text{Hilb}^d_X$ taking $[s]$ to $Z(s)$ is an embedding.
%\end{proposition}

For each reduced scheme $X$ of finite type over $k$, each weakly linked net $\mathfrak L$ of invertible sheaves on $X$ over a $\mathbb Z^n$-quiver $Q$, and each path $\gamma$ in $Q$ we denote by $X^{\mathfrak L}_{\gamma}$ the union of the irreducible components of $X$ over which $\vp^{\mathfrak L}_\gamma$ is generically zero. We omit the superscript if $\mathfrak L$ is clear from the context. Notice that if $\gamma=\gamma_2\gamma_1$ then 
$X^{\mathfrak L}_{\gamma}=X^{\mathfrak L}_{\gamma_2}\cup X^{\mathfrak L}_{\gamma_1}$. Recall from \cite{Esteves_et_al_2021} that $\mathfrak L$ is called \emph{maximal} if $X^{\mathfrak L}_a$ is an irreducible component of $X$ for each arrow $a$ of $Q$.

\begin{proposition}\label{maxlink} Let $X$ be a reduced scheme of finite type over $k$ and $\mathfrak L$ be a maximal linked net of invertible sheaves on $X$  over a $\mathbb Z^n$-quiver $Q$. For each arrow type $\mathfrak a$ of $Q$, put $X^{\mathfrak L}_{\mathfrak a}:=X^{\mathfrak L}_a$, where $a$ is an arrow of type $\mathfrak a$. Then the assignment 
$\mathfrak a\mapsto X^{\mathfrak L}_{\mathfrak a}$ is well defined and a bijection between the set of arrow types of $Q$ and the set of irreducible components of $X$.
\end{proposition}

\begin{proof} We divide the proof in two steps:

\vskip0.2cm

{\bf Step 1:} Let $b_1,b_2$ be two arrows leaving the same vertex and $b$ an arrow of the same type as $b_2$ leaving the final vertex of $b_1$. Then 
$X_{b}=X_{b_2}$. If in addition $b_1\neq b_2$ then $X_b\neq X_{b_1}$.

\vskip0.1cm

Indeed, if $b_1=b_2$ then 
$\text{Ker}(\vp_{bb_1})=
\text{Ker}(\vp_{b_1})$ by \cite{Esteves_et_al_2021}, Lem.~6.6, hence $X_{b}\cup X_{b_1}=X_{b_1}$, or equivalently, $X_{b}\subseteq X_{b_1}$. Equality follows as both sides are irreducible components of $X$. 

If $b_1\neq b_2$, then $b_1$ and $b_2$ have different types, hence $\text{Ker}(\vp_{b_1})\cap\text{Ker}(\vp_{b_2})=0$. It follows that $X_{b_1}\neq X_{b_2}$. Let $b'$ be an arrow of the same type as $b_1$ leaving the final vertex of $b_2$. Then $bb_1$ and $b'b_2$ connect the same two vertices and thus 
$\vp_{bb_1}$ is a nonzero scalar multiple of 
$\vp_{b'b_2}$. It follows that
\begin{equation}\label{XbXb}
X_{b}\cup X_{b_1}=X_{b'}\cup X_{b_2}.
\end{equation}
But $X_{b_1}\neq X_{b_2}$ and the subcurves of $X$ in \eqref{XbXb} are irreducible. Then $X_{b}=X_{b_2}$.

\vskip0.2cm

{\bf Step 2.} Let $a_1$ and $a_2$ be arrows of $Q$. Then $X_{a_1}=X_{a_2}$ if and only if $a_1$ and $a_2$ have the same type. 

\vskip0.1cm

The assertion is clear if $a_1=a_2$. Assume $a_1\neq a_2$. Let $v_i$ be the inital vertex of $a_i$ for $i=1,2$. Let $\gamma$ be an admissible path connecting $v_1$ to $v_2$. We prove the claim by induction on the length of $\gamma$.

If $v_1=v_2$ then $a_1$ and $a_2$ do not have the same type and $X_{a_1}\neq X_{a_2}$ by Step~1.
Assume $v_1\neq v_2$. Write $\gamma=b\gamma'$, where $b$ is an arrow, the last of $\gamma$. Let $a'_2$ be an arrow of the same type as $a_2$ leaving the initial vertex of $b$. Then $X_{a_2}=X_{a'_2}$ by Step~1. By induction, 
$a_1$ and $a'_2$ have the same type if and only if $X_{a_1}=X_{a'_2}$. But then $a_1$ and $a_2$ have the same type if and only if $X_{a_1}=X_{a_2}$, as claimed.

\vskip0.2cm

It follows that the assignment 
$\mathfrak a\mapsto X_{\mathfrak a}$ is a well-defined injection. It is surjective because $\vp_{\gamma}=0$ for each minimal circuit $\gamma$, and hence $X=\bigcup X_a$ where the union runs through the arrows $a$ of $\gamma$.
\end{proof}

\begin{Def} 
Let $X$ be a reduced scheme of finite type over $k$. A closed subscheme 
$Z\subseteq X$ is said to be of
\emph{pure codimension one} if the intersection of $Z$ with each irreducible component of $X$ has all irreducible components of codimension one in that component. A coherent sheaf $F$ 
is said to have \emph{rank one} if $F$ is generically invertible everywhere, and \emph{depth one} if its associated points are the generic points of $X$.
A global section $s$ of $F$ defines a closed subscheme of $X$, denoted $Z(s)$, whose sheaf of ideals is the image of the induced map $F^X\to\mathcal O_X$, and which we call the \emph{zero scheme} of $s$, where $F^X:=Hom_{\mathcal O_X}(F,\mathcal O_X)$. 
\end{Def}

An invertible sheaf has rank one and depth one. If $F$ is a rank-1, depth-1 sheaf on $X$ and $s$ does not vanish generically anywhere on $X$, then $F^X\to\mathcal O_X$ is injective, and hence the sheaf of ideals of $Z(s)$ is isomorphic to $F^X$. Furthermore, $Z(s)$ has pure codimension one in $X$.

\begin{Def} Let $X$ be a reduced scheme of finite type over $k$ and $\mathfrak L$ be a weakly linked net of coherent sheaves on $X$ over a $\mathbb Z^n$-quiver $Q$. We denote by $H^0(X,\mathfrak L)$ the representation obtained by taking global sections. Given a subrepresentation $\h\subseteq H^0(X,\mathfrak L)$ of pure dimension 1, we let $Z(\h)$ denote the intersection of the zero schemes of the elements of $V^\h_v$ at all vertices $v$ of $Q$, viewed as sections of the corresponding coherent sheaves.
\end{Def}

The representation $H^0(X,\mathfrak L)$ is a weakly linked net of vector spaces. It is a linked net if so is $\mathfrak L$. It may not be pure though, nor finitely generated. 
%In fact, if $\mathfrak L$ is an exact maximal linked net of invertible sheaves on $X$, then the set of dimensions of the spaces associated to $H^0(X,\mathfrak L)$ is unbounded: for each $b\in\mathbb N$ and each irreducible component $Y$ of $X$ there is a vertex $v$ of $Q$ such that the sheaf associated to $v$ by $\mathfrak L$ has restriction to $Y$ of degree greater than $b$.

\begin{proposition}\label{Zw} Let $X$ be a reduced scheme of finite type over $k$. Let $\mathfrak L$ be an exact maximal linked net of invertible sheaves on $X$ over a $\mathbb Z^n$-quiver $Q$ and $\h\subseteq H^0(X,\mathfrak L)$ a finitely generated pure subrepresentation of dimension 1. Then 
$Z(\h)$ has pure codimension one in $X$ and $[Z(\h)]=c_1(L)\cap [X]$ for each invertible sheaf $L$ associated to $\mathfrak L$.
\end{proposition}

\begin{proof} Since $\mathfrak L$ is exact, each two invertible sheaves $L$ and $M$ associated to $\mathfrak L$ have the same class in the Grothendieck group of coherent sheaves on $X$ by \cite{Esteves_et_al_2021}, Prop.~9.3, and thus $c_1(L)=c_1(M)$. 

Since $\h$ is finitely generated, by Theorem~\ref{12n} there are vertices $v_1,\dots,v_m$ of $Q$ forming an oriented polygon minimally generating  $\h$. It follows that 
$$
Z(\h)=\bigcap_{i=1}^m Z(s_i),
$$ 
where $s_i$ is a generator of $V^\h_{v_i}$ for $i=1,\dots,m$. 

For each $i=1,\dots,m$, let $L_i$ be the invertible sheaf on $X$ associated to $v_i$ by $\mathfrak L$, and let 
$\psi_i\colon L_i\to L_{i+1}$ be the associated map and $Y_i$ the union of the irreducible components of $X$ where $\psi_i$ vanishes generically. (For convenience, we put $v_{m+1}:=v_1$ and $L_{m+1}:=L_1$.) 

If $m=1$, as there are arrows of each type leaving $v_1$, it follows from Proposition~\ref{maxlink} that $s_1$ is generically nonzero on each irreducible component of $X$, and hence $Z(\h)$ has pure codimension one in $X$ and $[Z(\h)]=c_1(L_1)\cap [X]$.

Assume $m>1$. Since $\{v_1,\dots,v_m\}$ minimally $1$-generates $\h$, the vertices are unrelated for $\h$, and thus $s_i$ is a global section of the subsheaf $\text{Ker}(\psi_i)$ for each $i$. The subsheaf is a coherent sheaf on $Y_i$ which has rank one and depth one. Furthermore, the section $s_i$ is nonzero generically on $Y_i$, because there is no unrelated polygon for $\h$ with more than $m$ vertices by Theorem~\ref{12n}. Let $(s_1,\dots,s_m)$ denote the corresponding section of the sum
$$
M:=\bigoplus_{i=1}^m \text{Ker}(\psi_i).
$$
The sum is a torsion-free, rank-one sheaf on $X$. It has the same class in the Grothendieck group of coherent sheaves on $X$ as $L_1$ by Proposition~\ref{LDeltaAb}. The section $(s_1,\dots,s_m)$ vanishes generically nowhere, whence $[Z(s_1,\dots,s_m)]=c_1(M)\cap[X]$. Finally, it is clear that 
$Z(s_1,\dots,s_m)=Z(s_1)\cap\cdots\cap Z(s_m)$.
\end{proof}

\begin{Def} Let $X$ be a reduced scheme of finite type over $k$. 
A \emph{linked net of linear series} on $X$ over a $\mathbb Z^n$-quiver $Q$ is the data $\mathfrak g$ of a maximal linked net $\mathfrak L$ over $Q$ of invertible sheaves on $X$ and a finitely generated pure subnet $\g$ of $H^0(X,\mathfrak L)$. It is said to have rank $r$ if $\g$ has dimension $r+1$. Also, we say $\mathfrak g$ has sections in $\mathfrak L$, and write $\mathfrak g=(Q,\mathfrak L,\g)$. 
\end{Def}

\begin{proposition}\label{gH=P(H)} Let $X$ be a reduced scheme of finite type over $k$ and $\mathfrak g=(Q,\mathfrak L,\g)$ be a linked net of linear series on $X$. Let $H$ be the intersection of all collections of vertices 1-generating $\g$. Then $\g$ is $1$-generated by the finite set $H$. Furthermore, if  $\g$ exact then 
$P(H)=H$ and for each vertex $u$ of $Q$, there is a section of $L_u$ in $V_u$ not vanishing generically anywhere on $X$ if and only if $u\in H$. Finally, $\vp^{u_1}_{u_2}(V_{u_1})\neq 0$ for any $u_1,u_2\in H$
\end{proposition}

\begin{proof} The first statement follows immediately from Lemma~\ref{minsup1}, since $\g$ is finitely generated. As for the second statement, let $u$ be a vertex of $Q$. There is an admissible path $\nu$ connecting a vertex $w$ of $H$ to $u$ such that $\vp_\nu(V_w)=V_u$. If $u\not\in H$, then $\nu$ is nontrivial, and thus Proposition~\ref{maxlink} yields that all sections of $L_u$ in $V_u$ vanish on $X^{\mathfrak L}_{\mathfrak a}$ for each $\mathfrak a$ appearing as the type of an arrow in $\nu$.  

Now, assume $\mathfrak V$ is exact. Assume that all sections of $L_u$ in $V_u$ vanish completely on a component of $X$. Since $H\subseteq P(H)$, we will finish the proof of the second statement by showing that $u\not\in P(H)$. That is the case indeed, since Proposition~\ref{maxlink} yields that all sections of $L_u$ in $V_u$ vanish on $X^{\mathfrak L}_{\mathfrak a}$ for a certain arrow type $\mathfrak a$. Let $a$ be the arrow arriving at $u$ with type $\mathfrak a$. Then $\vp_\gamma(V_u)=0$ for any reverse path $\gamma$ by Proposition~\ref{maxlink}. Since $\mathfrak V$ is exact, $\vp_a(V_x)=V_u$, where $x$ is the initial vertex of $a$. Then $H':=(H-\{u\})\cup\{x\}$ would also $1$-generate $\g$. But since $H$ is minimum, $H\subseteq H'$, and thus $u\not\in H$. 

Actually, $u\not\in P(H)$. Indeed, let $z\in H$. We have just seen that there is a section $s$ of $L_z$ in $V_z$ that does not vanish completely on any component of $X$. Let $\mu$ be an admissible path connecting $z$ to $u$. Since $\vp_{\gamma}\vp_\mu(s)=0$, Proposition~\ref{maxlink} yields that the concatenation $\gamma\mu$ is not admissible. Thus $\mu$ must contain an arrow of type $\mathfrak a$. Since this is true for each $z\in H$, we have that $u\not\in P(H)$.

As for the third statement, let $\nu$ be an admissible path connecting $u_1$ to $u_2$, both in $H$. If $\vp_\nu(V_{u_1})=0$, then it follows from Proposition~\ref{maxlink} that all the sections of $L_{u_1}$ in $V_{u_1}$ vanish on $X^{\mathfrak L}_{\mathfrak a}$ for each arrow type $\mathfrak a$ not appearing in $\nu$. But this contradicts the second statement. 
\end{proof}

%We refrain from directly calling $\mathfrak g$ a ``limit linear series," as the term has been used already in different contexts, and it may not be suitable to add a new one.

\begin{proposition}\label{dual} Let $X$ be a reduced scheme of finite type over $k$ and $\mathfrak g=(Q,\mathfrak L,\g)$ be a linked net of linear series on $X$ with $\g$ exact. Let $H$ be the minimum collection of vertices 1-generating $\mathfrak V$. Then there is a unique rational map $X\dashrightarrow\lp(\g^*)$ which assigns to each $P\in X$ the unique subrepresentation $\h\subseteq\g^*$ associating to each vertex $u\in H$ the class $[\varepsilon^u_P]$ of the evaluation map $\varepsilon^u_P\colon V_u\to L_u|_P$.
\end{proposition}

(There might vertices $v$ of $Q$ not in $H$ for which the evaluation map $V_v\to L_v|_P$ vanishes for $P$ on a whole irreducible component of $X$.) 

\begin{proof} By Proposition~\ref{gH=P(H)}, for each $u\in H$ there is a section of $L_u$ in $V_u$ that does not vanish generically anywhere on $X$, whence there is an open dense subset of $X$ parameterizing $P\in X$ with nonzero evaluation map $\varepsilon^u_P$. As $H$ is finite, there is an open dense subset $U$ of $X$ such that $\varepsilon^u_P\neq 0$ for each $P\in U$ and $u\in H$.

Let $P\in U$. For each vertex $v$ of $Q$, put  $\varepsilon^v_P:=\varepsilon^u_P(\vp^u_v)^{-1}\colon V_v\to L_u|_P$, where $u$ is the shadow of $v$ in $H$. It is well-defined because $P(H)=H$ by Proposition~\ref{gH=P(H)} and $\vp^u_v(V_u)=V_v$. Clearly, $\varepsilon^v_P\neq 0$.

If there is a subrepresentation $\h\subseteq\g^*$ of pure dimension 1 associating to each vertex $u\in H$ the class $[\varepsilon^u_P]$, compatibility yields that $\h$ associates to each vertex $v$ of $Q$ the class  $[\varepsilon^v_P]$. 

Conversely, the assignment of the nonzero class $[\varepsilon^v_P]$ to each vertex $v$ of $Q$ is a subrepresentation $\h\subseteq\g^*$. Indeed, given an arrow $a$ connecting a vertex $v_1$ to a vertex $v_2$ of $Q$, let $u_1$ and $u_2$ be their respective shadows in $H$. Then $\vp^{u_1}_{u_2}(V_{u_1})\neq 0$ by Proposition~\ref{gH=P(H)}, and thus $\vp^{u_2}_{v_2}\vp^{u_1}_{u_2}(V_{u_1})\neq 0$. Hence there is an admissible path from $u_1$ to $v_2$ through $u_2$. Then either $\vp^{v_1}_{v_2}(V_{v_1})=0$ or $\vp^{v_1}_{v_2}\vp^{u_1}_{v_1}=\vp^{u_1}_{v_2}=\vp^{u_2}_{v_2}\vp^{u_1}_{u_2}$. In the first case, $[\varepsilon^{v_2}_P]\vp^{v_1}_{v_2}(V_{v_1})=0$, whereas in the second case $[\varepsilon^{v_2}_P]\vp^{v_1}_{v_2}|_{V_{v_1}}=[\varepsilon^{v_1}_P]$. In any case, $\h$ is a subrepresentation of $\g^*$. Clearly, $\h$ is of pure dimension 1, so $\h\in\lp(\g^*)$. 
\end{proof}

For each scheme $X$ projective over $k$, let $\text{Hilb}_X$ denote the Hilbert scheme of $X$, parameterizing closed subschemes.

\begin{proposition}\label{maplpg} Let $X$ be a reduced projective scheme over $k$. Let $\mathfrak g=(Q,\mathfrak L,\g)$ be a linked net of linear series on $X$ such that $\mathfrak L$ and $\g$ are exact. Then the assignment of $Z(\h)$ to each pure subnet $\h\subseteq\g$ of dimension~1 is the underlying function of a scheme morphism $\lp(\g)\to\text{Hilb}_X$.
\end{proposition}

\begin{proof} That the function is well-defined follows from Proposition~\ref{Zw}, since 
$\mathfrak L$ is exact and maximal, and each subnet $\h\subseteq\g$ is finitely generated because so is $\g$. 

The function $\lp(\g)\to\text{Hilb}_X$ is a morphism of schemes if it is locally so. Let $\h\in\lp(\g)$. Then 
$\h$ is finitely generated by \cite{Esteves_et_al_2021}, Prop.~6.8. Then, by Theorem~\ref{12n}, there are vertices $v_1,\dots,v_m$ forming an oriented polygon minimally generating $\h$. In fact, there is an open neighborhood $U\subseteq\lp(\g)$ parameterizing subrepresentations generated by 
$\{v_1,\dots,v_m\}$. On $U$ the function is given by taking $\h$ to the intersection $Z(s_1)\cap\cdots\cap Z(s_m)$, where $s_i$ is a nonzero element of $V^\h_{v_i}$, thus a section of the invertible sheaf $L_i$ associated by $\mathfrak L$ to $v_i$ for each $i$. A family of $\h$ over $U$ corresponds thus to a family over $U$ of nonzero sections $s_i$ of $L_i$ for each $i$, and thus to a family of intersections $Z(s_1)\cap\cdots\cap Z(s_m)$ over $U$, which is flat over $U$, because the sheaves of ideals of its fibers have the same class in the Grothendieck group of $X$, and because $\lp(\g)$ is reduced by Theorem~\ref{main_result}. Hence the restriction of the function 
$\lp(\g)\to\text{Hilb}_X$ to $U$ is a scheme morphism. Since $\h$ was arbitrary, the statement of the proposition follows.
\end{proof}

Assume now that we are given a \emph{regular smoothing} of a connected reduced projective scheme $X$ over $k$, that is, the data of a discrete valuation ring $R$ with residue field $k$, a flat projective map $\pi\colon\mathcal X\to B$ from a regular scheme $\mathcal X$ to $B:=\text{Spec}(R)$, and an isomorphism from the special fiber to $X$. Let $\text{Hilb}_{\mathcal X/B}$ be the relative Hilbert scheme, parameterizing closed subschemes. 
%It is $B$-flat and i
Its fibers over $B$ are the corresponding 
Hilbert schemes of the fibers of $\pi$.

Assume as well that we are given a linear series $(L_\eta,V_\eta)$ of rank $r$ on the generic fiber of $\pi$. Let $Q$ be the arising $\mathbb Z^n$-quiver. Let $\mathfrak L$ be the arising maximal exact linked net over $Q$ of invertible sheaves on $X$ and 
$\g$ the arising pure exact finitely generated subnet of 
$H^0(X,\mathfrak L)$ of dimension $r+1$; see \cite{Esteves_et_al_2021}, \S 3. The data $\mathfrak g=(Q,\mathfrak L,\g)$ is thus a linked net of linear series on $X$ of rank $r$. 

\begin{Def} Call $\mathfrak g=(Q,\mathfrak L,\g)$ as above the \emph{limit} of $(L_\eta,V_\eta)$ along $\pi$, or simply a \emph{limit linked net of linear series}.
\end{Def}

\begin{theorem} Let $X$ be a reduced scheme projective over $k$. Let $\mathfrak g=(Q,\mathfrak L,\g)$ be a linked net of linear series on $X$. Assume $\mathfrak g$ is a limit. Then $\g$ is smoothable. In addition, $\lp_H(\g)$ is a degeneration of the small diagonal in $\prod_{v\in H}\mathbb P(V^{\g}_v)$ for each finite set of vertices $H$ of $Q$ with $P(H)=H$.
\end{theorem}

\begin{proof} As seen in \cite{Esteves_et_al_2021}, \S 3, the linked net $\g$ is smoothable. The remaining is a consequence of Theorem~\ref{smoothg}.
\end{proof}

\begin{theorem}\label{final} Let $X$ be a connected reduced scheme projective over $k$. Let $\mathfrak g=(Q,\mathfrak L,\g)$ be a linked net of linear series on $X$. If $\mathfrak g$ is the limit of $(L_\eta,V_\eta)$ along a regular smoothing $\pi\colon\mathcal X\to B$ of $X$, then the image of $\lp(\g)$ in 
$\text{Hilb}_X$ is the associated reduced subscheme of the 
limit of the image of 
$\mathbb P(V_{\eta})$ in the generic fiber of 
$\text{Hilb}_{\mathcal X/B}$.
\end{theorem}

\begin{proof} Here $B$ is the spectrum of a discrete valuation ring $R$ with residue field $k$. As seen in \cite{Esteves_et_al_2021}, \S 3, there is a representation of $Q$ in the category of invertible sheaves on $\mathcal X$ restricting to $\mathfrak L$ on $X$ and a subrepresentation $\mathfrak M$ of the associated representation of global sections in the category of $R$-modules which is a smoothing of $\g$. Let $H$ be a finite set of vertices of $Q$ generating $\g$ with $P(H)=H$. The generic fiber of 
$\lp_H(\mathfrak M)$ is $\mathbb P(V_{\eta})$ and the special fiber is $\lp_H(\g)$, which is naturally isomorphic to $\lp(\g)$. By Theorem~\ref{smoothg}, the scheme $\lp_H(\mathfrak M)$ is reduced and flat over $B$. 

Arguing as in the proof of Proposition~\ref{maplpg}, using that $\lp_H(\mathfrak M)$ is reduced, there is a natural associated $B$-morphism of schemes $\lp_H(\mathfrak M)\to\text{Hilb}_{\mathcal X/B}$ restricting to the scheme morphism 
$\lp(\g)\to\text{Hilb}_X$ on the special fiber and to the natural embedding $\mathbb P(V_{\eta})\to \text{Hilb}_{\mathcal X/B}$ into the generic fiber. Since $B:=\text{Spec}(R)$, the scheme-theoretic image 
$Y\subseteq\text{Hilb}_{\mathcal X/B}$ of $\lp_H(\mathfrak M)$ is a $B$-flat closed subscheme. Since $\lp_H(\mathfrak M)\to\text{Hilb}_{\mathcal X/B}$ restricts to a closed embedding over the general point of $B$, the fiber of $Y$ over the general point is the image of $\mathbb P(V_\eta)$ in the generic fiber of 
$\text{Hilb}_{\mathcal X/B}$. And the fiber of $Y$ over the special point is a certain closed subscheme of $\text{Hilb}_X$ whose associated reduced subscheme is the image of $\lp(\g)$ in $\text{Hilb}_X$.
%, which is a closed embedding because so is on each fiber over $B$.
\end{proof}

\section{Example}\label{secexample}

\begin{Exa}\label{example} Let $X$ be the reduced union of $n+1$ distinct lines $M_0,\dots,M_n$ on the plane $\p^2_k$ over the field $k$ for $n>0$. Picking coordinates $Y,Z,W$ for $\p^2_k$, we have $M_i=y_iY+z_iZ+w_iW$ for $y_i,z_i,w_i\in k$ for $i=0,\dots,n$. Assume no three of the $M_i$ intersect. Thus $X$ is reduced and its singularities are nodes.

Let $F$ be a plane curve of degree $n+1$. Let 
$\mathcal X\subset \p^2_k\times_k B$ be the surface given by $M_0\cdots M_n+TF=0$, where $B:=\text{Spec}(k[[T]])$. Assume $F$ does not contain any node of $X$. Then $\mathcal X$ is regular. Denote by $X$ its special fiber over $B$ and by $\mathcal X_\eta$ its generic fiber.

Consider the invertible sheaf $\mathcal L:=\mathcal O_{\mathcal X}(1)$. The coordinates $Y,Z,W$ can be thought of as sections of $\mathcal O_{\p^2_k}(1)$. Consider the linear system of sections $V_\eta$ of $L_\eta:=\mathcal L|_{\mathcal X_\eta}$ generated by (the pullbacks of) $Y,Z,W$.

Let 
$$
\mathbb Z^{n+1}_{n+1}:=\{(d_0,\dots,d_n)\in\mathbb Z^{n+1}\,|\, \sum d_i=n+1\}.
$$
Put $v:=(1,\dots,1)$, the multidegree of $\mathcal L|_X$,  and $s_i:=(1,\dots,-n,1,\dots,1)$, the multidegree of $\mathcal O_{\mathcal X}(M_i)|_X$, for each $i=0,\dots,n$. 
Recall the associated $\mathbb Z^n$-quiver $Q:=Q(v,s_0,\dots,s_n)$, with vertex set $Q_0:=v+\mathbb Z s_0+\cdots+\mathbb Z s_n$ and arrow set $Q_1:=A_0\cup\cdots\cup A_n$ with $A_i:=\{(u,u+s_i)\,|\,u\in Q_0\}$ for $i=0,\dots,n$; see \cite{Esteves_et_al_2021}, \S~2. 

Let $\mathfrak L$ be the representation of $Q$ induced by $(L_\eta,V_\eta)$ in the category of invertible sheaves on $X$ of degree $n+1$. It is a maximal exact linked net; see \cite{Esteves_et_al_2021}, Prop.~3.1. The sheaf $L_u$ associated to $u\in Q_0$ has multidegree $u$. Let 
$\mathfrak V\subseteq H^0(X,\mathfrak L)$ be the subrepresentation induced by $(L_\eta,V_\eta)$. It is a pure exact linked net 
$1$-generated by the set $H$ of effective multidegrees in $\mathbb Z^{n+1}_{n+1}$, that is, by 
$$
H:=\{v,(n+1)e_0,\dots,(n+1)e_n\},
$$
where $e_0,\dots,e_n$ is the canonical basis of $\mathbb Z^{n+1}$; see \cite{Esteves_et_al_2021}, Prop.~3.2. For simplicity, put $v_i:=(n+1)e_i$ for each $i$. 

Notice that $H$ is a ``star." The arrows of $Q$ connecting vertices of $H$ are just the pairs $a_i:=(v_i,v)$ for $i=0,\dots,n$. It follows from Proposition~\ref{cap} and Theorem~\ref{Main_Theorem} that $\lp(\g)$ has at most $n+2$ irreducible components, the nonempty among   $\lp(\g)_v,\lp(\g)_{v_0},\dots,\lp(\g)_{v_n}$. Furthermore, it follows as well from Proposition~\ref{cap} that $\lp(\g)_{v_i}$ intersects at most only $\lp(\g)_v$ for each $i=0,\dots,n$. We will see below that we can remove ``at most" from the last two sentences. 

We may assume for simplicity that $y_i\neq 0$ for each $i$. Clearly, the subspace $V_v\subseteq H^0(X,L_v)$ is that induced by the coordinates of $\p^2_k$. Here is how to obtain $V_{w_0}$: Add to $(L_\eta,V_\eta)$ the base points of the pencil $M_0\cdots M_n+tF=0$ given by $M_1\cdots M_n=F=0$. This is obtained multiplying $Y,Z,W$ by 
$M_1\cdots M_n$. Now, $M_0,Z,W$ generate the same system $V_\eta$ since $y_0\neq 0$. Also, $M_0M_1\cdots M_n=-TF$ on $X_\eta$, which is a nonzero multiple of $F$, and can thus be replaced by $F$. Restrict the net of hypersurfaces generated by $F,ZM_1\cdots M_n,WM_1\cdots M_n$ to $X$; observe it has base points given by $M_1\cdots M_n=F=0$. Subtracting them we obtain $(L_{v_0},V_{v_0})$. The $(L_{v_i},V_{v_i})$ for $i=1,\dots,n$ are obtained similarly.

Notice that there is a linear combination of $Y,Z,W$ that does not vanish totally on any $M_i$, hence a section of $V_v$ spanning a subrepresentation of $\mathfrak V$ of pure dimension 1. So 
$\lp(\g)_v^*\neq\emptyset$. Likewise, since the curve $F$ does not contain the line $M_i$, we have $\lp(\g)_{v_i}^*\neq\emptyset$ for each $i$. Thus 
$\lp(\g)_v,\lp(\g)_{v_0},\dots,\lp(\g)_{v_n}$ are the irreducible components of $\lp(\g)$. Furthermore, 
$\lp(\g)_{v_0}$ intersects $\lp(\g)_v$ by Proposition~\ref{lemma: non_exact_point_is_in_frontier}, as the subnet generated by the section corresponding to $M_0$ in $V_v$ and that corresponding to $YM_1\cdots M_n$ in $V_{v_0}$ is of pure dimension 1 and is minimally generated by $\{v,v_0\}$. By the symmetry, $\lp(\g)_{v_i}$ intersects $\lp(\g)_v$ for each $i=1,\dots,n$ as well.

Notice that $H$ is equal to its hull $P(H)$. Indeed, given $u\in Q_0-H$, consider an admissible path $\gamma$ in $Q$ connecting $v$ to $u$. If it contains arrows of a certain type at least twice, then all admissible paths connecting a vertex of $H$ to $u$ contain an arrow of that type. Thus $u\in P(H)$ only if $\gamma$ is simple. Suppose $\gamma$ is simple. Then $\gamma$ cannot have maximum length $n$, as otherwise $u=v_i$ for some $i$. But then each vertex of $H$ can be connected to $u$ by an admissible path passing through $v$. Since $\gamma$ is nontrivial, $u\not\in P(H)$.

It thus follows from Proposition~\ref{HHPPHH} that $\lp(\g)=\lp(\g)_H$. So 
$\lp(\g)$ can be described as the quiver Grassmannian of pure 1-dimensional subrepresentations of a representation by vector spaces of the quiver with $n+2$ vertices and $2(n+1)$ arrows connecting one of the vertices, called ``central", to the other $n+1$ vertices, called ``outer", and back. Indeed, 
$\vp^{v_i}_{v_j}=\vp^v_{v_j}\vp^{v_i}_v$ for $i\neq j$ for $\mathfrak L$, hence we need only specify $\vp^v_{v_i}$ and $\vp^{v_i}_v$ for each $i=0,\dots,n$. 

Now, $Y,Z,W$ gives us a basis for $V_v$, which we fix, thus identifying $V_v$ with $k^3$. Then, for each $i=0,\dots,n$, the map $\vp^v_{v_i}$ has a one-dimensional kernel, generated by $(y_i,z_i,w_i)$. By the exactness of $\g$, this kernel is the image of $\vp^{v_i}_v$. Thus, by choosing a basis for $V_{v_i}$ appropriately, we have that $\vp^{v_i}_v$ and $\vp^v_{v_i}$ can be respectively represented, up to multiplication by a nonzero scalar, by the matrices
$$
\begin{bmatrix} y_i & 0 & 0 \\ 
z_i & 0 & 0 \\ 
w_i & 0 & 0
\end{bmatrix}\quad\text{and}\quad 
\begin{bmatrix} 0 & 0 & 0 \\ 
-z_i & y_i & 0 \\ 
-w_i & 0 & y_i
\end{bmatrix}.
$$

\begin{figure}[ht]
    \centering
    \begin{subfigure}{.48\textwidth}
        \centering
        \includegraphics[scale = .83,trim = 133 462 187 200 , clip = true]{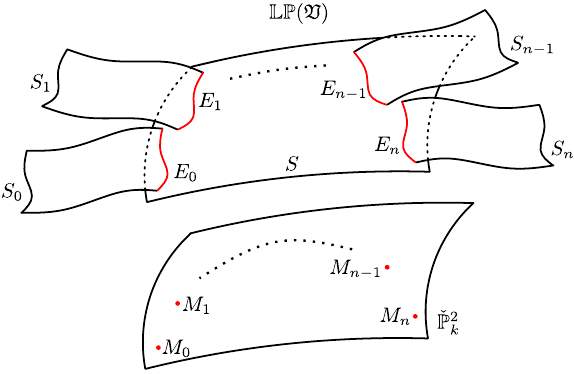}
        \subcaption{$\lp(\g)$.}
        \label{fig:LP(V)}
    \end{subfigure}
    \quad
    \begin{subfigure}{.48\textwidth}
        \centering
        \includegraphics[scale = 0.83,trim = 140 510 180 150 , clip = true]{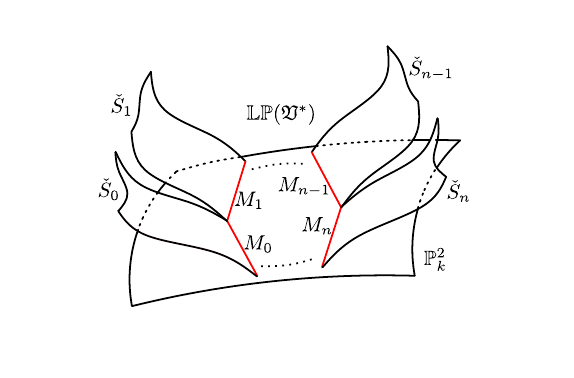}
        \subcaption{$\lp(\g^*)$.}
        \label{fig:LP(V_dual)}
    \end{subfigure}
    \caption{Geometric descriptions of $\lp(\g)$ and $\lp(\g^*)$.}
\end{figure}

We gave a precise description of $\mathfrak V$, yielding one for $\lp(\g)$: It can be obtained by taking the union of the blowup $S$ of the dual $\check{\mathbb P}^2_k$ along the points $(y_i:z_i:w_i)$ corresponding to $M_i$ for each $i$, and $n+1$ copies $S_0,\dots,S_n$ of $\mathbb P^2_k$, identifying a line on $S_i$ with the exceptional divisor $E_i$ on $S$ over $(y_i:z_i:w_i)$ for 
each $i$. 

We have decided not to describe in the current article when the scheme morphism $\lp(\g)\to\text{Hilb}_X$ is an embedding. The reader may check himself that in our example this is the case.

It is straightforward that $(L_v,V_v)$ has no base points. Also, the intersection of the hypersurfaces $F,ZM_1\cdots M_n,WM_1\cdots M_n,M_0\cdots M_n$ is the intersection of $F$ and $M_1\cdots M_n$ because $y_0\neq 0$, whence $(L_{v_0},V_{v_0})$ has no base points either. By symmetry, $(L_{v_i},V_{v_i})$ has no base points for any $i$. It follows that the rational map 
$\psi\colon X\dashrightarrow\lp(\g^*)$ described in Proposition~\ref{dual} is defined everywhere. 

We have a precise description of $\mathfrak V$, thus of $\mathfrak V^*$ as well, which yields one for $\lp(\g^*)$: It can be obtained by taking the union of $\check S:=\mathbb P^2_k$ with the blowups $\check S_i$ of $\mathbb P^2_k$ at $(y_i:z_i:w_i)$ for $i=0,\dots,n$, identifying the exceptional divisor $\check E_i$ on $\check S_i$ with the line $M_i$ on $\check S$ for 
$i=0,\dots,n$. 

Notice that even though $\lp(\g^*)$ and $\lp(\g)$ have the same number of components, they are not isomorphic, since $\lp(\g^*)$ admits a triple intersection of irreducible components but $\lp(\g)$ does not. Also, the image of the composition of $\psi|_{M_i}$ with the projection 
$\lp(\g^*)\to\mathbb P(V^*_{v_i})$ spans the whole space, thus $\psi(X)$ intersects $\check S$ in finitely many points, and is not in particular equal to $M_0\cdots M_n$. The curve $\psi(X)$ is far from being a union of lines, the way it flexes captures for instance the limits of the flexes along the pencil $M_0\cdots M_n+TF=0$. 
\end{Exa}

\bibliographystyle{plain}
\bibliography{mybibtex}

\vspace{0.5cm}

{\smallsc Instituto de Matem\'atica Pura e Aplicada, 
Estrada Dona Castorina 110, 22460-320 Rio de Janeiro RJ, Brazil}

{\smallsl E-mail address: \small\verb?esteves@impa.br?}

\vspace{0.2cm}

{\smallsc Instituto de Matem\'atica Pura e Aplicada, 
Estrada Dona Castorina 110, 22460-320 Rio de Janeiro RJ, Brazil}

{\smallsl E-mail address: \small\verb? silvase@impa.br?}

\vspace{0.2cm}

{\smallsc Instituto de Matem\'atica Pura e Aplicada, 
Estrada Dona Castorina 110, 22460-320 Rio de Janeiro RJ, Brazil}

{\smallsl E-mail address: \small\verb?renanmath@ufc.br?}

\end{document}